\documentclass[12pt,a4paper,reqno]{amsart}
\usepackage{latexsym, amsfonts, amsmath, amsthm, amssymb, amscd,
color}
\usepackage{pstricks}
\usepackage{caption}
\usepackage{listings}
\usepackage{graphics}    
\usepackage{graphicx}
\usepackage{color}
\usepackage{mathrsfs}
\usepackage{xparse}
\usepackage{xpatch}
\usepackage{enumerate}
\usepackage{epstopdf}
\usepackage{hyperref}
\usepackage{verbatim}
\usepackage{mathtools}
\usepackage{subeqnarray}
\usepackage{multirow}
\usepackage{amsaddr}

\makeatletter
\newsavebox\myboxA
\newsavebox\myboxB
\newlength\mylenA

\newcommand*\xoverline[2][0.75]{%
\sbox{\myboxA}{$\m@th#2$}%
\setbox\myboxB\null
\ht\myboxB=\ht\myboxA%
\dp\myboxB=\dp\myboxA%
\wd\myboxB=#1\wd\myboxA
\sbox\myboxB{$\m@th\overline{\copy\myboxB}$}
\setlength\mylenA{\the\wd\myboxA}
\addtolength\mylenA{-\the\wd\myboxB}%
\ifdim\wd\myboxB<\wd\myboxA%
\rlap{\hskip 0.5\mylenA\usebox\myboxB}{\usebox\myboxA}%
\else
\hskip -0.5\mylenA\rlap{\usebox\myboxA}{\hskip
	0.5\mylenA\usebox\myboxB}%
\fi}
\makeatother

\makeatletter

\makeatletter
\def\subsection{\@startsection{subsection}{3}%
  \z@{.5\linespacing\@plus.7\linespacing}{.1\linespacing}%
  {\bf}}
\makeatother

\newtheorem{theorem}{Theorem}[section]
\newtheorem{proposition}[theorem]{Proposition}
\newtheorem{lemma}[theorem]{Lemma}
\newtheorem{corollary}[theorem]{Corollary}

\newtheorem{conjecture}[theorem]{Conjecture}

\newtheorem{problem}{Problem}

\newtheorem*{lemma*}{Lemma}
\newtheorem*{theorem*}{Theorem}
\newtheorem*{conjecture*}{Conjecture}

\theoremstyle{remark}

\numberwithin{equation}{section}
\numberwithin{theorem}{section}
 \def\noninv{\operatorname{nin}}
\def\Rp{\mathbb{R}^{+}}
\newcommand{\ceq}{\coloneqq}
\newcommand{\qn}[1]{\ensuremath{[#1]_q}}
\newcommand{\pqn}[2]{\ensuremath{[#1]_{#2}}}

\newcommand{\qpow}[3]{\ensuremath{\left(\frac{\qn{#1}}{\qn{#2}}\right)^{#3}}}
\newcommand{\qbin}[2]{\genfrac{(}{)}{0pt}{}{#1}{#2}_q}
\newcommand{\pqbin}[3]{\genfrac{(}{)}{0pt}{}{#1}{#2}_{#3}}

\newcommand{\balph}{\ensuremath{\bsym\alpha}}
\newcommand{\bbet}{\ensuremath{\bsym\beta}}
\newcommand{\tp}{T_{\infty}}
\newcommand{\tm}{\sfT}
\newcommand{\fm}{\mathsf{F}}
\newcommand{\mm}{\mathsf{M}}

\newcommand{\bsym}{\ensuremath{\boldsymbol}}

\newcommand{\nibl}{\mathop{\rm niblings}}
\newcommand{\scind}{\mathop{\rm scindex}}
\newcommand{\parents}{\mathop{\rm parents}}

\def\eop{\hbox{\kern1pt\vrule height6pt width4pt
depth1pt\kern1pt}\medskip}

\newcommand{\be}{\begin{equation}}
\newcommand{\ee}{\end{equation}}

\newcommand{\parent}[1]{\mathop{\rm parent}\nolimits(#1)}

\newcommand{\sgn}{\mathop{\rm sgn}\nolimits}

\newcommand{\diag}{\mathop{\rm diag}\nolimits}

\newcommand{\wt}{\mathop{\rm wt}\nolimits}

\newcommand{\scrd}{{\mathcal{D}}}

\newcommand{\scrf}{{\mathcal{F}}}

\newcommand{\scrm}{{\mathcal{M}}}
\newcommand{\scrn}{{\mathcal{N}}}
\newcommand{\scro}{{\mathcal{O}}}

\newcommand{\scrr}{{\mathcal{R}}}

\newcommand{\bfa}{{\mathbf{a}}}
\newcommand{\bfb}{{\mathbf{b}}}
\newcommand{\bfc}{{\mathbf{c}}}

\newcommand{\bfe}{{\mathbf{e}}}
\newcommand{\bfh}{{\mathbf{h}}}

\newcommand{\bfx}{{\mathbf{x}}}

\newcommand{\Z}{{\mathbb Z}}
\newcommand{\N}{{\mathbb N}}
\newcommand{\Q}{{\mathbb Q}}

\newcommand{\R}{{\mathbb R}}

\newcommand{\myge}{\succeq}

\newcommand{\sfT}{{\sf T}}

\newcommand{\inv}{{\rm inv}}

\newcommand{\cross}{{\rm cr}}

\newcommand{\nest}{{\rm ne}}

\newmuskip\pFqskip
\mathchardef\pFcomma=\mathcode`, 

\begin{document}

\title[Trees, forests, and total positivity: I. \texorpdfstring{\MakeLowercase{$q$}}{TEXT}-trees and \texorpdfstring{\MakeLowercase{$q$}}{TEXT}-forests matrices]{Trees, forests, and total positivity: \\I. \texorpdfstring{\MakeLowercase{$q$}}{TEXT}-trees and \texorpdfstring{\MakeLowercase{$q$}}{TEXT}-forests matrices}
\author[T.~Gilmore]{
     Tomack Gilmore}

     \address{{\small Department of Mathematics, University College London,\\
                    London WC1E 6BT, UK\\\vskip0.5cm\tt{t.gilmore@ucl.ac.uk}}}


\begin{abstract}We consider matrices with entries that are polynomials in $q$ arising from natural $q$-generalisations of two well-known formulas that count: forests on $n$ vertices with $k$ components; and trees on $n+1$ vertices where $k$ children of the root are smaller than the root. We give a combinatorial interpretation of the corresponding statistic on forests and trees and show, via the construction of various planar networks and the Lindstr\"om-Gessel-Viennot lemma, that these matrices are coefficientwise totally positive. We also exhibit generalisations of the entries of these matrices to polynomials in \emph{eight} indeterminates, and present some conjectures concerning the coefficientwise Hankel-total positivity of their row-generating polynomials.\end{abstract}

\maketitle{}
\tableofcontents

\section{Introduction}\label{sec.intro}

Two well-known enumerative formulas that often arise in the study of trees and forests are\footnote{See~\cite{Clarke58},~\cite[pp.~26-27]{Moon70},~\cite[p.~70]{Comtet_74},~\cite[pp.~25-28]{Stanley_86},~\cite{Chauve99},~\cite{Sokal_20} or~\cite{Avron2016}. See also~\cite{Francon_74,Pitman_02,Sagan_83,Guo17,Riordan68} and~\cite[pp.~235-240]{Aigner18} for related information.}
\be f_n\ceq\;(n+1)^{n-1}\ee
and
\be\label{eq.treen} t_n\ceq\;(n+1)^n.\ee

The first formula, $f_n$, gives the number of forests of rooted labelled trees on the vertex set $[n]:=\{1,2,\ldots,n\}$ and has the refinement:
\be\label{eq.forest.refine}f_n\;=\;\sum_{k=0}^{n}f_{n,k}\ee
where each summand 
\be\label{eq.forests.refine1} f_{n,k}\ceq\binom{n-1}{k-1}n^{n-k}\ee is the number of forests on $[n]$ comprised of $k$ components (that is, $k$ rooted labelled trees). The first few $f_n$ and $f_{n,k}$ are
\be\label{eq.fnkArr}
\begin{array}{c|rrrrrrr|r}
\downarrow n, \rightarrow k & 0 & 1 & 2 & 3 & 4 & 5 & 6&f_n\\\hline 
0 & 1 &  &  &  &  &  & &1\\
1 & 0 & 1 & & & & & &1\\
2 & 0 & 2 & 1 &  &  &  & &3\\
3 & 0 & 9 & 6 & 1 &  &  & &16\\
4 & 0 & 64 & 48 & 12 & 1 &  & &125\\
5 & 0 & 625 & 500 & 150 & 20 & 1 & &1296\\
6 & 0 & 7776 & 6480 & 2160 & 360 & 30 & 1&16807
\end{array}
\ee
(see~\cite[A061356/A137452 and A000272]{OEIS}). By adding a new vertex (labelled $0$) and connecting it to the root of each component in a forest, we see that $f_n$ is also the number of trees on vertices labelled $0,1,\ldots,n$, rooted at $0$, where the root has precisely $k$ children.

The formula in~\eqref{eq.treen}, $t_n$, gives the number of rooted labelled trees on the vertex set $[n+1]$ and has the following refinement similar to~\eqref{eq.forest.refine}:
\be\label{eq.tree.refine}
t_n\;=\;\sum_{k=0}^{n}t_{n,k},\ee
where each summand
\be\label{eq.tree.refine1}t_{n,k}\ceq\binom{n}{k}n^{n-k}\ee
is the number of rooted labelled trees on the vertex set $[n+1]$ in which precisely $k$ children of the root are lower-numbered than the root (see~\cite{Chauve99,Chauve00,Sokal_20}). The first few $t_n$ and $t_{n,k}$ are:
\be\label{eq.treetable}\begin{array}{c|rrrrrrr|r}
\downarrow n, \rightarrow k & 0 & 1 & 2 & 3 & 4 & 5 & 6&t_n\\\hline 
0 & 1 &  &  &  &  &  & &1\\
1 & 1 & 1 & & & & & &2\\
2 & 4 & 4 & 1 &  &  &  & &9\\
3 & 27 & 27 & 9 & 1 &  &  & &64\\
4 & 256 & 256 & 96 & 16 & 1 &  & &625\\
5 & 3125 & 3125 & 1250 & 250 & 25 & 1 & &7776\\
6 & 46656 & 46656 & 19440 & 4320 & 540 & 36 & 1&117649
\end{array}
\ee
(see~\cite[A071207]{OEIS}). We point out that this array has an alternative combinatorial interpretation; in~\cite{Chauve99,Chauve00} the authors gave a bijection between trees on $n+1$ vertices where $k$ children of the root are lower-numbered than the root and forests on $n+1$ vertices with $k+1$ components where the vertex $n+1$ is a leaf.

The main goal of this paper is to prove the total positivity of some matrices that arise from generalisations of $f_{n,k}$ and $t_{n,k}$. Recall that a matrix is \emph{totally positive} (\emph{strictly totally positive}) if all of its minors are nonnegative (strictly positive, respectively).\footnote{The terms ``totally nonnegative'' and ``totally positive'' have been used by previous authors (see, for example,~\cite{Gantmacher_02,Fomin_00,Fallat_11}) in place of what we have called here ``totally positive'' and ``strictly totally positive'' respectively. In studying the literature it is important to clarify which sense of total positivity is being used, since many theorems are valid only for strictly totally positive matrices.} Total positivity has many applications across combinatorics, statistical physics, representation theory, and background material on the topic can be found in~\cite{Karlin_68,Gantmacher_02,Pinkus_10,Fallat_11}.

Our first result concerns the \emph{forests (of rooted labelled trees) matrix}
\be \fm\ceq(f_{n,k})_{n,k\geq 0}\;=\;\left(\binom{n-1}{k-1}n^{n-k}\right)_{n,k\geq 0},\ee
(the first few rows of which appear in~\eqref{eq.fnkArr}) and the 
\emph{(rooted labelled) trees matrix}
\be \tm\ceq\left(t_{n,k}\right)_{n,k\geq 0}\;=\;\left(\binom{n}{k}n^{n-k}\right)_{n,k\geq 0},\ee
(the first few entries of which appear in~\eqref{eq.treetable}). We have:

\begin{theorem}\label{thm.tree}
The matrices $\fm$ and $\tm$ are totally positive.
\end{theorem}

The total positivity of the forests and trees matrices has been proven very recently in~\cite{Sokal_21f,Sokal_21t} using different methods to the ones we employ in the present paper. Sokal observes in~\cite{Sokal_21f} that $\fm$ is the exponential Riordan array $\scrr[F,G]$ with $F(t):=1$ and $G(t)$ the \emph{tree function}~\cite{Corless96}
\be\label{eq.treefunc} G(t)\ceq\sum_{n=1}^{\infty}n^{n-1}\frac{t^n}{n!},\ee
and similarly in~\cite{Sokal_21t} Sokal and Chen observe that $\tm$ is the exponential Riordan array $\scrr[F,G]$ with
\be F(t)\ceq \sum_{n=0}^{\infty}n^n\frac{t^n}{n!}\;=\;\frac{1}{1-G(t)}\ee
and $G(t)$ as above. In~\cite{Sokal_21f} and~\cite{Sokal_21t} the author(s) exploit these intepretations of $\fm$ and $\tm$ and make use of the production-matrix method to prove these matrices are totally positive.\footnote{We remark that in~\cite{Sokal_21f,Sokal_21t} the author(s) prove that generalisations of $\fm$ and $\tm$ (that are different to those considered here) are coefficientwise totally positive. We discuss their results further in Section~\ref{sec.furthercomments}.} The proof of Theorem~\ref{thm.tree} we offer below instead takes a leaf from~\cite{Brenti_95,Fomin_00,Fallat_11}, making use of the Lindstr\"om-Gessel-Viennot (LGV) lemma and planar networks (see Section~\ref{subsec.planarNetworks}). 

The matrices $\fm$ and $\tm$ are closely related; by studying the entries one can, in effect, see the woods for the trees since for $n>0$
\be f_{n,k}\; = \;\binom{n-1}{k-1}n^{n-k}\;=\;\frac{k}{n}\binom{n}{k}n^{n-k}=\frac{k}{n}t_{n,k},\ee
that is, $\fm$ and $\tm$ are related via a diagonal similarity transform (see Lemma~\ref{lem.transF} in Section~\ref{subsec.identities}):
\be\label{eq.FTtrans} \fm\;=\;\lim_{\epsilon\to0}\;\diag((n+\epsilon)_{n\geq 0})^{-1}\tm\diag((k+\epsilon)_{k\geq 0})\ee
where $\diag(\bfa)$ denotes (for the sequence $\bfa=(a_n)_{n\geq0}$) the diagonal matrix with $(n,n)$-entry equal to $a_n$ and all other entries $0$ (note we must take a limit in the above to avoid division by zero in the first column of $\tm$). Total positivity is preserved under matrix multiplication (this follows from the Cauchy-Binet formula), so~\eqref{eq.FTtrans} shows that the total positivity of $\tm$ implies that of $\fm$ (a diagonal matrix with nonnegative entries on the diagonal is clearly totally positive).

In fact there are a number of identities we present in Section~\ref{subsec.identities} that allow one to go from $\tm$ to $\fm$ (and vice versa) via simple matrix multiplications that preserve total positivity. By specialising $q$ to $q=1$ in Lemma~\ref{lem.fqtqT} below we obtain the identity
\be\label{eq.fM} \tm\;=\;\fm\cdot M\ee
where $M\ceq(m_{n,k})_{n,k\geq 0}$ is the unit-lower-triangular matrix with $(n,k)$-entry
\be m_{n,k}\;=\;\frac{n!}{k!}\ee
for $n\geq k$ and all other entries $0$. The matrix $M$ is easily shown to be totally positive (see Section~\ref{subsec.identities}), so it follows from~\eqref{eq.fM} that the total positivity of $\fm$ implies that of $\tm$. In order to prove Theorem~\ref{thm.tree} it therefore suffices to prove that the forests matrix $\fm$ is totally positive; we do this in Section~\ref{sec.pointwise} by constructing a planar network with nonnegative rational weights and showing that the corresponding path matrix is $\fm$. The total positivity of $\fm$ then follows from the LGV lemma (see Section~\ref{subsec.LGV}). 

In this paper, however, we are chiefly concerned with $q$-generalisations of the forests and trees matrices. The formula $f_{n,k}$ has a perfectly natural $q$-analogue, thus we define the \emph{$q$-forests matrix} $\fm(q)\ceq(f_{n,k}(q))_{n,k\geq 0}$ to be the matrix with $(n,k)$-entry 
\be\label{eq.fnkqdef} f_{n,k}(q)\ceq\begin{cases}\delta_{nk} & \textrm{if }k=0,\\
                           \displaystyle\qbin{n-1}{k-1}(\qn{n})^{n-k} &\textrm{if } n\geq k\geq 1,\\
                           0&\textrm{otherwise.}
                           \end{cases}\ee
where 
\be\qn{n}\ceq \frac{1-q^n}{1-q}\;=\;\begin{cases}1+q+\cdots+q^{n-2}+q^{n-1}&\textrm{if }n>0,\\
0 & \textrm{if }n=0,
\end{cases}\ee
$\qbin{n}{k}$ is the $q$-binomial coefficient:
\be \qbin{n}{k} \ceq\frac{\qn{n}!}{\qn{k}!\qn{n-k}!},\ee
in which
\be \qn{n}!\ceq \begin{cases}\displaystyle\prod_{j=1}^n\qn{j} &\textrm{if }n>0,\\
1& \textrm{if }n=0,
\end{cases}\ee
and $\delta_{nk}$ is the Kronecker delta:
\be \delta_{nk}\ceq\begin{cases}1&\textrm{if } n=k,\\
                    0&\textrm{if }n\neq k\end{cases}\ee
It follows easily from~\eqref{eq.fnkqdef} that $f_{n,k}(q)$ is a monic self-reciprocal polynomial of degree $(n-1)^2-(k-1)^2$, and the first few rows of $\fm(q)$ are:
\be\left[
\begin{array}{ccccc}
 1 &  &  &  &  \\
 0 & 1 & &  &  \\
 0 & q+1 & 1 &  &  \\
 0 & q^4+2 q^3+3 q^2+2 q+1 & q^3+2 q^2+2 q+1 & 1 &  \\
 \vdots&\vdots&\vdots&\vdots&\ddots
\end{array}
\right]\ee

The $q$-forests matrix counts forests on $n$ vertices with $k$ components with respect to some statistic, which we interpret combinatorially in Section~\ref{sec.combin.interp} (see Corollary~\ref{cor.niblingsq}). Please note that in~\cite{Sokal_21f} Sokal considers a \emph{different} generalisation of $\fm$ (he introduces indeterminates into $\fm$ that count forests of rooted labelled trees in terms of \emph{proper} and \emph{improper} edges), whereas our results grew from studying natural $q$-generalisations of the matrix entries. In Section~\ref{sec.furthercomments} we discuss how our generalisation gives rise to some curious conjectures concerning generalisations of well-known polynomial sequences including: the Schl\"afli-Gessel-Seo polynomials; the general Abel polynomials; $(p,q)$-Stirling cycle polynomials; and the reverse Bessel polynomials.

Each entry $t_{n,k}$ of the trees matrix has a similar natural $q$-analogue; we thus define the \emph{$q$-trees matrix} to be $\tm(q)\ceq(t_{n,k}(q))_{n,k\geq 0}$ where
\be t_{n,k}(q)\ceq\qbin{n}{k}(\qn{n})^{n-k}.\ee
Once more it is easy to see that $t_{n,k}(q)$ is a monic self-reciprocal polynomial of degree $n(n-1)-k(k-1)$, and the first few rows of $\tm(q)$ are:
\be \left[{\tiny
\begin{array}{ccccc}
 1 &  &  &  &  \\
 1 & 1 &  &  &  \\
 q^2+2 q+1 & q^2+2 q+1 & 1 &  &  \\
 q^6+3 q^5+6 q^4+7 q^3+6 q^2+3 q+1 & q^6+3 q^5+6 q^4+7 q^3+6 q^2+3 q+1 & q^4+2 q^3+3 q^2+2 q+1 & 1 &  \\
 \vdots&\vdots&\vdots&\vdots&\ddots\\
\end{array}}
\right]\ee
It follows from~\cite{Chauve99,Chauve00} that the entries of the matrix $\tm(q)$ count forests on $n+1$ vertices with $k+1$ components where the vertex $n+1$ is a leaf with respect to some statistic, and we give a combinatorial interpretation of $q$ in Section~\ref{sec.combin.interp}. We note that our $q$-generalisation of the tree matrix and the corresponding combinatorial interpretation differs from the generalisation of $\tm$ studied by Chen and Sokal in~\cite{Sokal_21t}.

The entries of $\fm(q)$ and $\tm(q)$ are polynomials in $q$ with integer coefficients, and there is a natural extension of total positivity to matrices whose entries are polynomials in one or more indeterminates $\bfx$. We equip the polynomial ring $\R[\bfx]$ with the \emph{coefficientwise
partial order}, that is, we say that $P$ is nonnegative
(and write $P \myge 0$)
in case $P$ is a polynomial with nonnegative coefficients.
We then say that a matrix with entries in $\R[\bfx]$ is
\emph{coefficientwise totally positive} (TP)
if all of its minors are polynomials with nonnegative coefficients, and if all of its minors of size $\leq r$ are polynomials with nonnegative coefficients we say that the matrix is \emph{totally positive of order $r$} (TP$_r$). Our main result is the following:
\begin{theorem}\label{thm.fqtqCoeff}
The matrices $\fm(q)$ and $\tm(q)$ with entries in the polynomial ring $\Z[q]$ are coefficientwise totally positive.\end{theorem}

It is plain to see that the above theorem reduces to Theorem~\ref{thm.tree} when $q$ is specialised to $q=1$. More generally, we say that a matrix $M$ with polynomial entries belonging to $\R[\bfx]$ is \emph{pointwise totally positive} on some domain $\scrd\subseteq\Rp$ if $M$ is totally positive for all $\bfx\in\scrd$ (by which we mean all the indeterminates $\bfx$ are specialised to values in $\scrd$). Coefficientwise total positivity of $M$ thus implies pointwise total positivity for all $\bfx \in\Rp$, but the converse is not true.

The main goal of this paper is to prove Theorem~\ref{thm.fqtqCoeff} and the structure of our proof is as follows: after reviewing some basic concepts in total positivity from the perspective of planar networks in Section~\ref{sec.prelim}, we then provide combinatorial interpretations of the entries of $\fm(q)$ and $\tm(q)$ in Section~\ref{sec.combin.interp} and establish some identities that show the coefficientwise total positivity of $\fm(q)$ implies that of $\tm(q)$. In Section~\ref{sec.pointwise} we construct a planar network with weights that are \emph{rational} and \emph{pointwise nonnegative} functions of $q$, and prove that the path matrix corresponding to this planar network agrees with $\fm(q)$; thanks to the LGV lemma this proves \emph{pointwise total positivity} of $\fm(q)$ and $\tm(q)$ for $q\in\Rp$ (Theorem~\ref{thm.tree} is then obtained by specialising $q$ to $q=1$). In Section~\ref{sec.coeffpos} we show how the network from Section~\ref{sec.pointwise} can be transformed into a planar network with weights that are \emph{polynomials} in $q$, thereby proving Theorem~\ref{thm.fqtqCoeff}. We conclude in Section~\ref{sec.furthercomments} with some further generalisations and open problems, some of which will be the subject of future work.

\section{Preliminaries}\label{sec.prelim}
Here we review some fundamental definitions and facts regarding total positivity. Please note that since our proofs rely on Lindstr\"om-Gessel-Viennot lemma, many of the results in this section are motivated by performing operations on planar networks (see~\cite{Fallat_11,Fomin_00,Brenti_95,Brenti_96}), from which various matrix identites follow.

\subsection{Total positivity and the Lindstr\"om-Gessel-Viennot lemma}\label{subsec.LGV}

One fundamental tool in the study of total positivity is the Lindstr\"om-Gessel-Viennot (LGV) lemma. Suppose $G$ is a locally finite acyclic digraph with source vertices $U\ceq\{u_0,u_1,\ldots,u_n\}$ and sink vertices $V\ceq\{v_0,v_1,\ldots,v_n\}$, where the weight $w_e$ of an edge $e$ is an element of some commutative ring $R$. Let $w(P(u_i,v_j))$ denote the product of the weights of the edges of a path $P(u_i,v_j)$ starting at $u_i$ and ending at $v_j$, and define
\be P_G(u_i\to v_j)\ceq\sum_{P(u_i,v_j)}w(P(u_i,v_j))\ee
to be the sum over all weighted paths between $u_i$ and $v_j$. The \emph{path matrix} corresponding to $G$ is then the matrix
\be P_G\ceq(P_G(u_i\to v_j))_{0\leq i,j\leq n}.\ee

We say that a \emph{family of nonintersecting paths} from $U$ to $V$ is an $n$-tuple $(P_0,P_1,\ldots,P_{n})$ of paths in $G$ such that:
\begin{enumerate}[(i)]
  \item{There exists a permutation $\sigma$ of $\{0,1,2,\ldots,n\}$ such that for each $i$, $P_i$ is a path from $u_i$ to $v_{\sigma(i)}$.}
  \item{Whenever $i\neq j$, paths $P_i$ and $P_j$ have no vertices in common (this includes endpoints).}
\end{enumerate}
The LGV lemma states that
\be \det(P_G)\;=\;\sum_{(P_0,P_1,\ldots,P_n):U\to V}\sgn(\sigma)\prod_{i=0}^nw(P_i),\ee
where $\sgn(\sigma)=(-1)^{\inv(\sigma)}$ is the signature of the permutation $\sigma$ arising from $P_i$ mapping $u_i$ to $v_{\sigma(i)}$ (note that $\inv(\sigma)$ denotes the number of \emph{inversions} of $\sigma$, that is, pairs $(i,j)$ in $\{0,1,\ldots,n\}$ such that $i<j$ and $\sigma(i)>\sigma(j)$). In particular, if the only permutation giving rise to nonempty families of nonintersecting paths is the identity then $\det(P_G)$ gives the sum over all weighted families of nonintersecting paths in $G$ that begin at $U$ and end at $V$, where the weight of each family is the product of the weights of the paths $(P_0,\ldots,P_n)$.

Suppose now the source and sink vertices of $G$ are \emph{fully compatible}, by which we mean that for any subset of sources $u_{n_1},u_{n_2},\ldots,u_{n_r}$ (where $n_1<n_2<\cdots<n_r$) and sinks $v_{k_1},v_{k_2},\ldots,v_{k_r}$ (where $k_1<k_2<\ldots<k_r$), the only permutation $\sigma\in\mathfrak{S}_r$ mapping each source $u_{n_i}$ to the sink $v_{k_{\sigma(i)}}$ that gives rises to nonempty families of paths is the identity. The LGV lemma then implies that \emph{every} minor of $P_G$ is a sum over families of nonintersecting paths between specified subsets of $U$ and $V$, where each family has weight $\prod w(P_i)$. If every edge has a weight that is a positive real number then $P_G$ is totally positive; if the weights of the network belong to the field $\Q(\bfx)$ of rational functions of $\bfx$ that are pointwise nonnegative on some domain $\scrd$ then $P_G$ is pointwise totally positive on $\scrd$; and if the weight of each edge is a polynomial in the indeterminates $\bfx$ with nonnegative integer coefficients then $P_G$ is totally positive in $\Z[\bfx]$ equipped with the coefficientwise order.

Determining whether a general locally finite acyclic digraph is fully compatible is non-trivial, however, if $G$ is embedded in the plane and the source and sink vertices $U$ and $V$ lie on the boundary of $G$ in the order ``first $U$ in reverse order, then $V$ in order'' then the topology of $G$ clearly implies that $U$ and $V$ are fully compatible (see, for example, Lemma~9.18 in~\cite{petreolle2020}). From now on we will refer to locally acyclic digraphs embedded in the plane with fully compatible sources and sinks as \emph{planar networks}, in the spirit of~\cite{Fomin_00}.

\subsection{The binomial-like planar network}\label{subsec.planarNetworks}

\begin{figure}\centering
\includegraphics[scale=0.5]{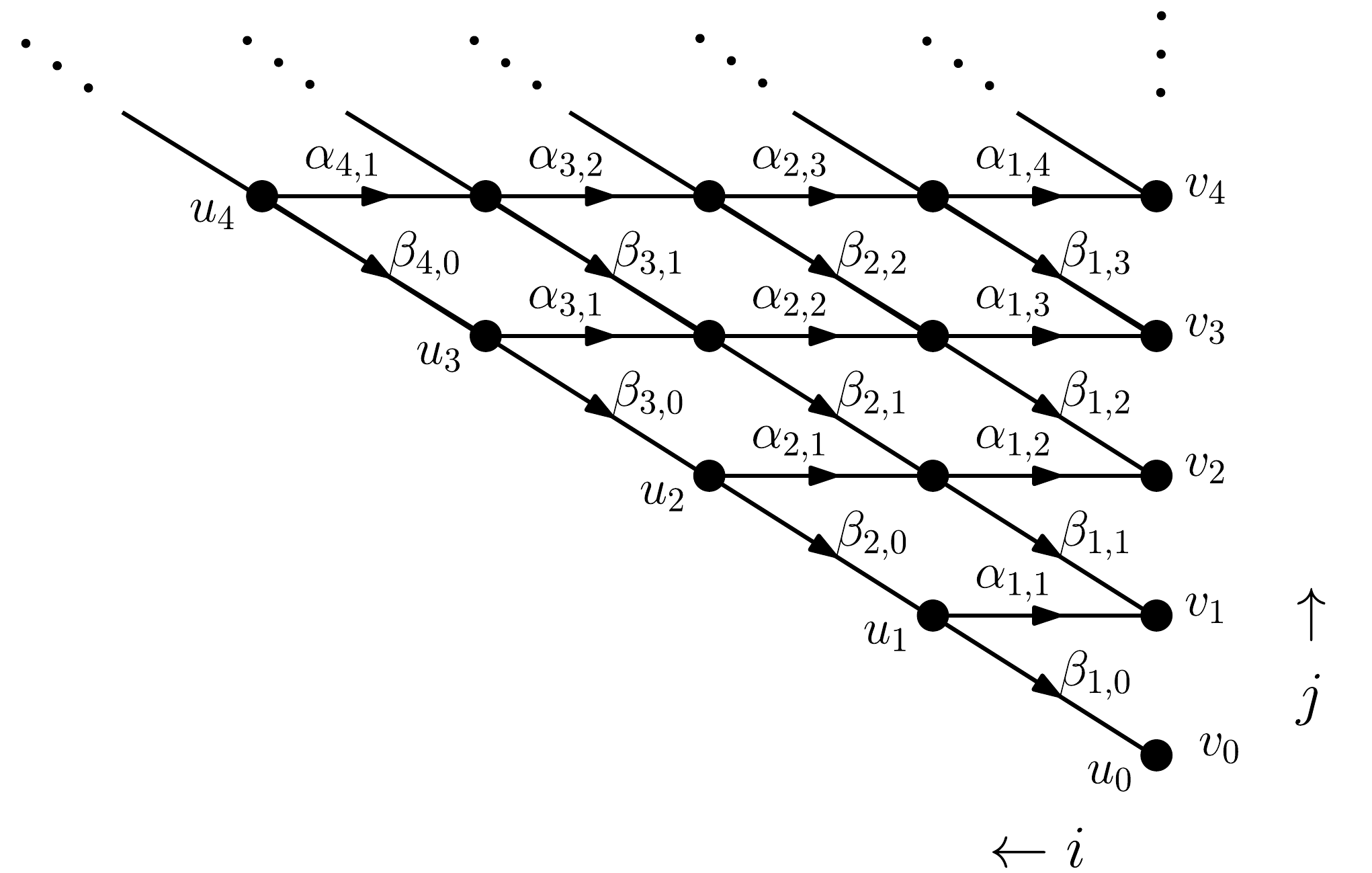}\caption{The binomial-like planar network $\scrn$ for vertices $u_0,u_1,\ldots,u_5$ and $v_0,v_1,\ldots,v_5$.}\label{fig.standardTriPN}
\end{figure}
The matrices $\tm(q)$ and $\fm(q)$ are unit-lower-triangular, so in the sequel we will be principally concerned with what we call the \emph{binomial-like planar network} (depicted in Figure~\ref{fig.standardTriPN}). This planar network consists of vertices indexed by $i$ ($j$) increasing horizontally from the right (vertically upwards, respectively) where $0\leq i\leq j$. Each source $u_{i}$ is the vertex with index $(i,i)$ while each sink $v_{j}$ is the vertex $(0,j)$. Each horizontal edge directed from vertex $(i,j)$ to vertex $(i-1,j)$ has weight $\alpha_{i,j-i+1}$, while the weight of a diagonal edge directed from vertex $(i,j)$ to $(i-1,j-1)$ has weight $\beta_{i,j-i}$. We also assume there exists a directed edge from $u_0$ to $v_0$ with weight 1.

The binomial-like planar network with weights $\balph=\{\alpha_{i,l}\}_{i,l\geq 1}$ and $\bbet=\{\beta_{i,l}\}_{i\geq1,l\geq 0}$ described above is denoted $\scrn$. Please note that the source and sink vertices of $\scrn$ are fully compatible, so it follows immediately from the LGV lemma that the path matrix $P_{\scrn}$ is automatically coefficientwise totally positive over $\Z[\balph,\bbet]$ in all of these indeterminates!

This is really the beauty of the LGV lemma (in the context of total positivity): given a lower-triangular matrix $M$ that appears to be coefficientwise totally positive in $\R[\bfx]$, one can try to prove total positivity by specialising $\alpha_{i,l}$, $\beta_{i,l}$ to suitable nonnegative elements of $\R[\bfx]$ and showing that under such a specialisation 
\be P_{\scrn}\;=\;M.\ee

Indeed, this is how we concoct our proof of Theorem~\ref{thm.fqtqCoeff}. We first show that specialising the weights of the binomial-like network to suitably chosen rational functions of $q$ yields a path matrix that agrees with $\fm(q)$ in Section~\ref{sec.pointwise}; since these rational functions are all nonnegative for $q\in\Rp$ we conclude that $\fm(q)$ is pointwise totally positive for $q\in\Rp$. In Section~\ref{sec.coeffpos} we then transform this binomial-like network to obtain a different planar network with weights that are \emph{polynomials} in $q$ in order to show that $\fm(q)$ is \emph{coefficientwise} totally positive. 

The planar network construction described above originally goes back to Brenti~\cite{Brenti_95} who also observed that if the binomial-like network $\scrn$ has weights $\balph$ and $\bbet$ that depend purely on the first index (that is, $\alpha_{i,l}=\alpha_{i,*}$ and $\beta_{i,l}=\beta_{i,*}$ for all $i$) then the entries of the path matrix $P_{\scrn}$ satisfy the $n$-dependent recurrence
\be P_{\scrn}(u_n\to v_k)\;=\; \alpha_{n,*} \, P_{\scrn}(u_{n-1}\to v_{k-1}) \:+\: \beta_{n,*} \, P_{\scrn}(u_{n-1}\to v_{k})\;.
 \label{eq.brenti.intro.n-dependent}
\ee
Similarly if the weights depend solely on the second index (that is, $\alpha_{i,l}=\alpha_{*,l}$ and $\beta_{i,l}=\beta_{*,l}$ for all $l$) then the entries of the path matrix satisfy the $k$-dependent recurrence
\be P_{\scrn}(u_{n}\to v_{k})\;=\; \alpha_{*,k} \, P_{\scrn}(u_{n-1}\to v_{k-1}) \:+\: \beta_{*,k} \, P_{\scrn}(u_{n-1}\to v_{k})\;.
 \label{eq.brenti.intro.k-dependent}
\ee

The most straightforward (and thus, eponymous) example that illustrates this connection between matrices with entries satisfying purely $n$- (or $k$-) dependent recurrences and binomial-like planar networks is the \emph{weighted binomial matrix}, 
\be B_{x,y}\ceq\left(\binom{n}{k}x^{n-k}y^k\right)_{n,k\geq 0},\ee the entries of which satisfy the $n$-dependent linear recurrence
\be (B_{x,y})_{n,k}\;=\;x(B_{x,y})_{n-1,k}+y(B_{x,y})_{n-1,k-1}\ee
for $n>0$ with initial condition $(B_{x,y})_{0,k}=\delta_{0k}$. The corresponding planar network is the network $\scrn$ with $\alpha_{i,l} =y$ and $\beta_{i,l}=x$, and it follows immediately from the LGV lemma that
\begin{corollary}
The weighted binomial matrix $B_{x,y}$ is coefficientwise totally positive in $\Z[x,y]$.
\end{corollary}
Matrices with entries given by purely $n$-dependent or purely $k$-dependent recurrences have relatively straightforward planar network representations, and often in these cases it is easy to deduce a straightforward planar network by studying the recurrence. 

Matrices with entries that satisfy recurrences dependent on both $n$ \emph{and} $k$, however, give rise to seemingly much more complex planar networks (see, for example,~\cite{Gilmore20}). For example, the entries of the forests matrix 
\be\fm\;=\;(f_{n,k})_{n,k\geq 0}\;=\;\left(\binom{n-1}{k-1}n^{n-k}\right)_{n,k\geq 0},\ee
satisfy the recurrence:
\be f_{n,k}\;=\;\left(\frac{n}{n-1}\right)^{n-1}\left[(n-1)f_{n-1,k}+\frac{1}{k-1}\left(n\left(\frac{n-1}{n}\right)^{k}+(k-n)\right)f_{n-1,k-1}\right]\ee
for $n>1$~ with initial conditions $f_{n,0}=\delta_{n0}$, and $f_{1,1}=1$. We note that we were able to find the correct weights for the binomial-like planar network corresponding to $\fm$ in \emph{spite} of this complicated recurrence, and not \emph{because} of it.

In the sequel we will also make use of the following fact. Suppose $\alpha_{i,l}=1$ and consider the sum over all weighted paths from $u_n$ to $v_k$ in $\scrn$ (denoted $P_{\scrn}(u_n\to v_k)$). By studying the network it is relatively easy to see that $P_{\scrn}(u_n\to v_k)$ can be expressed as a nested sum
\begin{eqnarray} P_{\scrn}(u_n\to v_k)&=&\sum_{i_1=1}^{k+1}\beta_{i_1,k+1-i_1}\sum_{i_2=i_1+1}^{k+2}\beta_{i_2,k+2-i_2}\cdots\sum_{i_{n-k}=i_{n-k-1}+1}^{n}\beta_{i_{n-k},n-i_{n-k}}\\\label{eq.symm}
&=&\sum_{1\leq i_1<\cdots<i_{n-k}\leq n}\beta_{i_1,k+1-i_1}\beta_{i_2,k+2-i_2}\cdots \beta_{i_{n-k},n-i_{n-k}}.\end{eqnarray}
Observe that if the weights $\beta_{i,l}$ are dependent only on the first index then $P_{\scrn}(u_n\to v_k)$ can be realised as the \emph{elementary symmetric polynomial}
\be\label{eq.eSymm} \bfe_{n-k}(X_1,X_2,\ldots,X_n)\ceq\sum_{1\leq i_1<\cdots<i_{n-k}\leq n}X_{i_1}X_{i_2}\cdots X_{i_{n-k}}\ee
where $X_n\ceq \beta_{n,*}$. If instead the weights $\beta_{i,l}$ depend solely on the second index then~\eqref{eq.symm} can be realised as the \emph{complete homogeneous symmetric polynomial}
\be\label{eq.hSymm} \bfh_{n-k}(X_1,X_2,\ldots,X_{k+1})\ceq\sum_{1\leq i_1\leq\cdots\leq i_{n-k}\leq k+1}X_{i_1}X_{i_2}\cdots X_{i_{n-k}}\ee
where $X_k=\beta_{*,k-1}$ (this observation again goes back to Brenti~\cite{Brenti_95}).

In this article we study total positivity primarily through the lens of planar networks; the next section dicusses how planar networks can be intepreted algebraically as matrix factorisations.

\subsection{Matrix factorisations}\label{subsec.MatFac}

The planar network approach outlined above allows us to easily write down various factorisations of the path matrix $P_{\scrn}$. Before discussing these factorisations we first clarify some notation. We use $\scrn$ to denote the binomial-like network described in the previous subsection with the set of horizontal weights $\balph$ and diagonal weights $\bbet$. Conversely, given a matrix $M$, in this section we will often use $\scrn_M$ to denote a planar network with corresponding path matrix $P_{\scrn_M}$ satisfying
\be P_{\scrn_M}\;=\;M.\ee
In this case $\scrn_M$ is referred to as a \emph{planar network representation} of the matrix $M$.

We will make much use of the following definition. For sequences $\bfa\ceq(a_n)_{n\geq 0}$ and $\bfb\ceq(b_n)_{n\geq0}$, let $L(\bfa,\bfb)$ denote the \emph{lower-bidiagonal matrix} with $(n,k)$-entry
\be \left(L(\bfa,\bfb)\right)_{n,k\geq 0}\;=\;\begin{cases}a_k & \textrm{if } n=k,\\
                                                          b_k & \textrm{if }n=k+1,\\
                                                          0&\textrm{otherwise.}
                                                          \end{cases}\ee
The first few rows and columns of $L(\bfa,\bfb)$ are:
\be L(\bfa,\bfb)\;=\;\left[
\begin{array}{cccccccc}
 a_0 &  &  &  &  &  &  &\\
 b_0 & a_1 &  &  &  &  &  &\\
  & b_1 & a_2 &  &  &  &  &\\
  &  & b_2 & a_3 &  &  &  &\\
  &  &  & b_3 & a_4 &  &  &\\
  &  &  &  & b_4 & a_5 &  &\\
  &  &  &  &  & b_5 & a_6 &\\
  &  &  &  &  & & \ddots&\ddots
\end{array}
\right].\ee
In case $a_n=1$ for all $n$ we abbreviate $L(\bfa,\bfb)$ to $L(\bfb)$, while if $b_n=0$ for all $n$ the matrix $L(\bfa,\bfb)$ is simply the diagonal matrix $\diag(\bfa)$ with $(n,n)$-entry $a_n$ and all other entries $0$. Lastly, if $a_n=1$ and the sequence $\bfb$ is constant, that is, $b_n=x$ for all $n$, we abbreviate $L(\bfb)$ to $L(x)$, and if we have $a_n=x$ and $b_n=y$ for all $n$ with $y\neq1$ we abbreviate $L(\bfa,\bfb)$ to $L(x,y)$.

The lower-bidiagonal matrix $L(\bfa,\bfb)$ has a planar network representation $\scrn_{L(\bfa,\bfb)}$ in which
\be P_{\scrn_{L(\bfa,\bfb)}}(u_n\to v_n)\;=\;a_n\ee
and
\be P_{\scrn_{L(\bfa,\bfb)}}(u_n\to v_{n-1})\;=\;b_{n-1}\ee
(see the diagram on the left in Figure~\ref{fig.bidiag}). Since we view these matrices through the lens of planar networks we will refer to the sequences $\bfa$ and $\bfb$ as \emph{edge-sequences corresponding to $L(\bfa,\bfb)$}. According to the LGV lemma we have:
\begin{figure}
\centering
\includegraphics[scale=0.5]{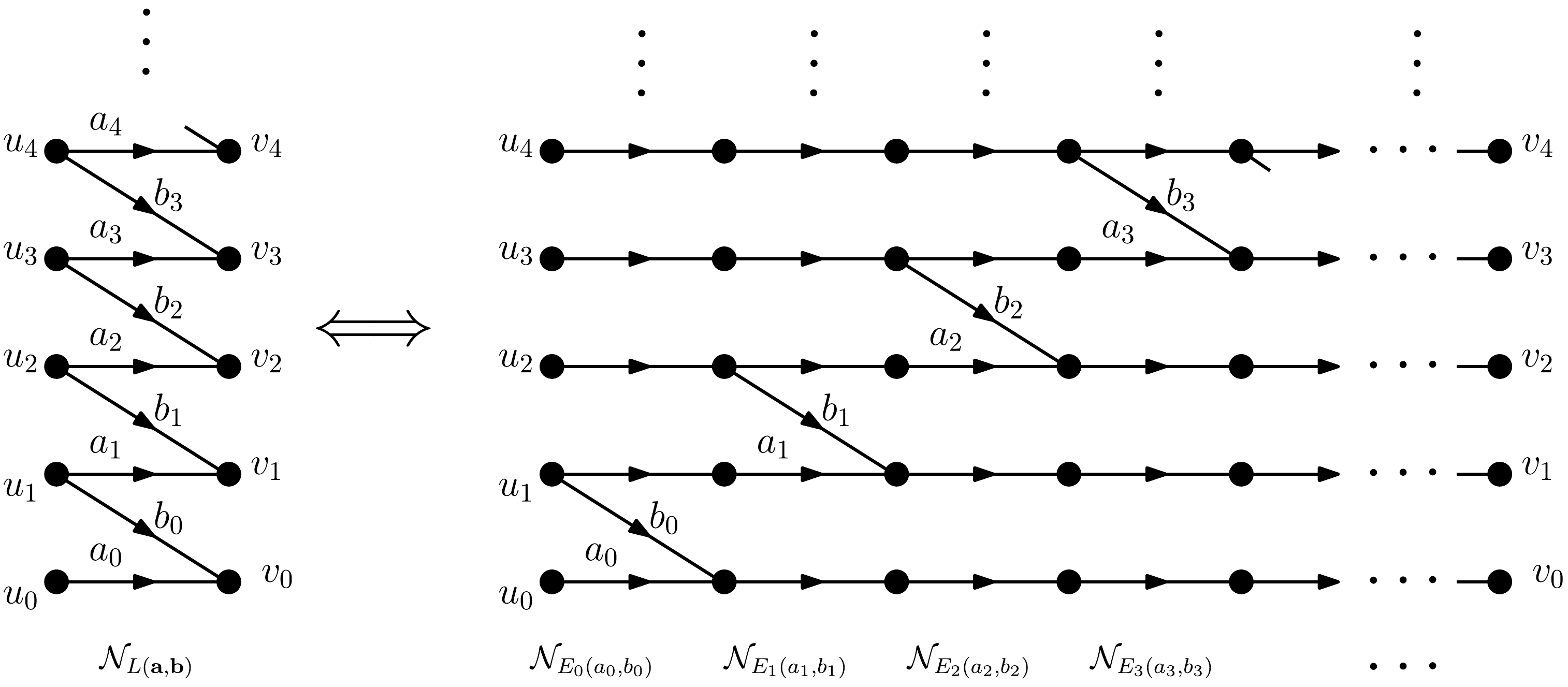}\caption{The planar network $\scrn_{L(\bfa,\bfb)}$ (left), and the planar network that arises from expressing $L(\bfa,\bfb)$ as a product of elementary lower-bidiagonal transfer matrices and concatenating the networks $N_{E_0(a_0,b_0)},N_{E_1(a_1,b_1)},\dots$ (right). Note that unlabelled directed edges are assumed to have weight 1.}\label{fig.bidiag}
\end{figure}

\begin{corollary}[Lower-bidiagonal matrices]\label{cor.bidiag}
The matrix $L(\bfa,\bfb)$ is totally positive in $\Z[\bfa,\bfb]$ equipped with the coefficientwise order. 
\end{corollary}

It is straightforward to verify that 
\be\label{eq.elemBidiag} L(\bfa,\bfb)\;=\;E_0(a_0,b_0)E_1(a_1,b_1)E_2(a_2,b_2)\cdots\ee
where $E_{n}(a,b)$ denotes the \emph{elementary (lower) bidiagonal matrix} with $(n,n)$-entry $a$, $(n+1,n)$-entry $b$, all other diagonal entries 1, and all other entries $0$. As with bidiagonal matrices, in case $a=1$ we abbreviate $E_n(1,b)$ to $E_n(b)$. We now present a useful lemma relating matrix products and concatenated planar networks:

\begin{lemma}[Concatenating planar networks]\label{lem.concat}
Suppose $M_1$ and $M_2$ are path matrices corresponding to two planar networks $\scrn_{M_1}$ and $\scrn_{M_2}$, where $\scrn_{M_1}$ has source vertices 
\be U\ceq\{u_0,u_1,\ldots,u_n\}\ee and sink vertices 
\be V\ceq\{v_0,v_1,\ldots,v_m\},\ee and $\scrn_{M_2}$ has source vertices 
\be U'\ceq\{u'_0,u'_1,\ldots,u'_m\}\ee and sink vertices 
\be V'=\{v'_0,v'_1,\ldots,v'_l\}.\ee Then the matrix product $M_1M_2$ is the path matrix corresponding to the planar network $\scrn_{M_3}$ obtained by concatenating $M_1$ and $M_2$, with source vertices $U$ and sink vertices $V'$, where each vertex $v_i\in V$ is identified with $u'_i\in U'$ for all $i$.
\end{lemma}
\begin{proof}
Since $v_i$ is identified with $u'_i$ we can write the sum over over weighted paths from $u_n$ to $v'_k$ in $\scrn_{M_3}$ as
\be P_{\scrn_{M_3}}(u_n\to v'_k)=\sum_{s=0}^{m}P_{\scrn_{M_1}}(u_{n}\to v_s)P_{\scrn_{M_2}}(u'_s\to v_k')\;=\;(M_1M_2)_{n,k}.\ee
\end{proof}

Concatenating planar networks thus corresponds to multiplying path matrices, and the matrices $M_1$ and $M_2$ in the above are referred to as \emph{transfer matrices} of the planar network $\scrn_{M_3}$. The elementary bidiagonal matrices on the right-hand side of~\eqref{eq.elemBidiag} are thus transfer matrices of $\scrn_{L(\bfa,\bfb)}$,\footnote{Note that the elementary bidiagonals described here are, in fact, referred to as \emph{column} transfer matrices, we will shortly also consider what we call \emph{diagonal transfer matrices}.} and concatenating the planar networks for $E_0(a_0,b_0),E_{1}(a_1,b_1),\ldots$ yields the planar network on the right in Figure~\ref{fig.bidiag}.
\begin{figure}\begin{center}
\includegraphics[scale=0.6]{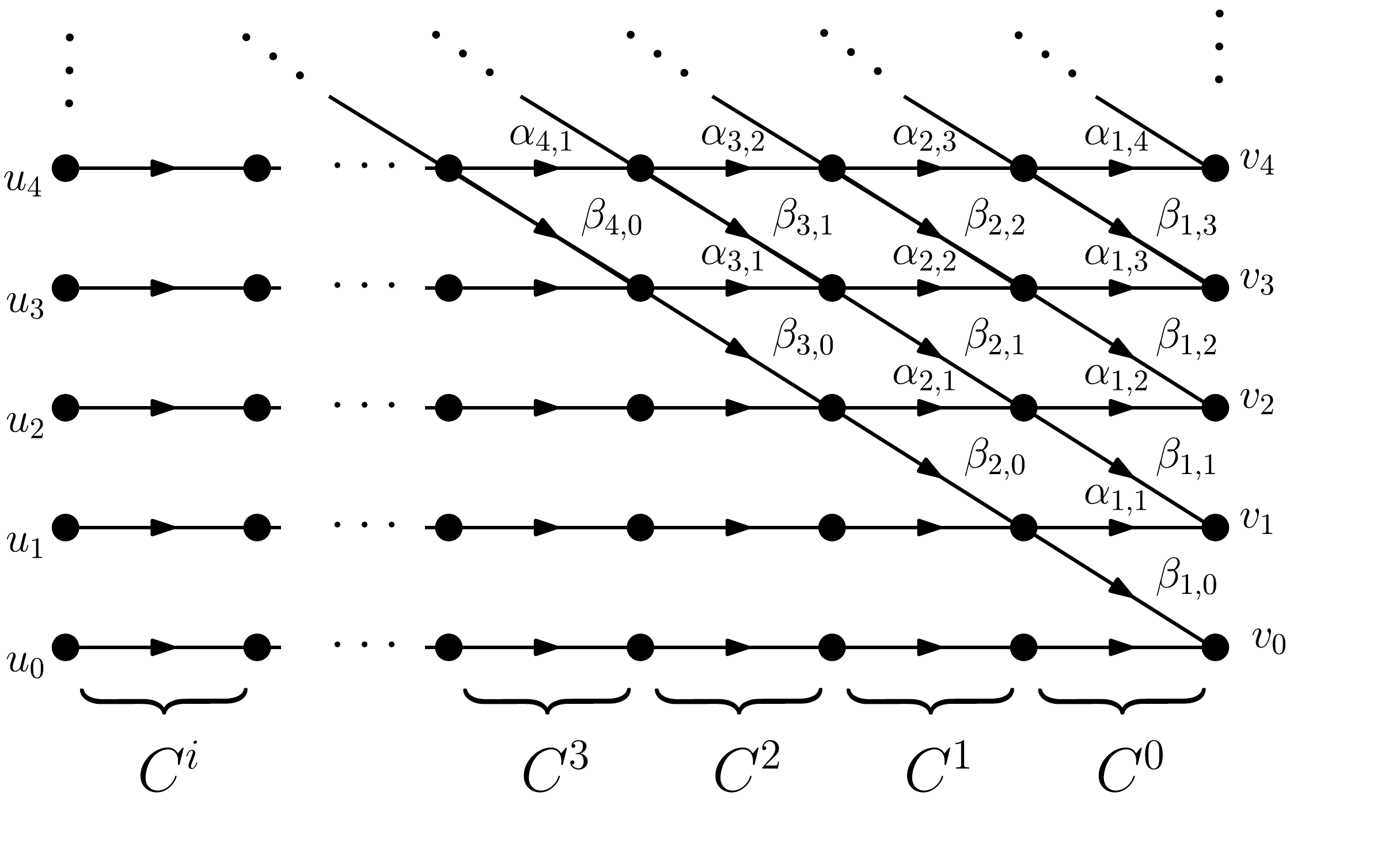}\caption{The planar network $\scrn_1$ up to sources $u_0,u_1,\ldots,u_4$ and $v_0,v_1,\ldots,v_4$ (unlabelled directed edges have weight $1$).}\label{fig.columntransfer}\end{center}
\end{figure}

We now interpret the binomial-like planar network $\scrn$ and its path matrix $P_{\scrn}$ from Subsection~\ref{subsec.planarNetworks} in light of these definitions. By extending the source vertices of $\scrn$ to the left so that they lie on the same vertical line we obtain a planar network $\scrn_1$ that is isomorphic to $\scrn$ (see Figure~\ref{fig.columntransfer}), in particular we have
\be P_{\scrn}(u_n\to v_k)\;=\;P_{\scrn_1}(u_n\to v_k).\ee 

Consider a subnetwork $C^i$ of $\scrn_1$ consisting of source vertices $U\ceq\{(i+1,j):j\in\N\}$ and sink vertices $V\ceq\{(i,j):j\in\N\}$. Let $\bfa_i$ denote the sequence of weights of horizontally directed edges emanating from vertex $(i+1,i+n)$ for increasing $n$, that is, $\bfa_i\ceq(a_n)_{n\geq 0}$ where
\be a_n\;=\;\begin{cases}\alpha_{i+1,n} & \textrm{if }n>0,\\
                            1 & \textrm{if }n=0\end{cases}.\ee
Similarly let $\bfb_i$ denote the sequence of weights of the diagonal edges emanating from vertex $(i+1,i+n+1)$ for increasing $n$, that is, $\bfb_i\ceq(\beta_{i+1,n})_{n\geq0}$.
It is easy to see that the path matrix corresponding to $C^i$ is the \emph{column transfer matrix}
\be P_{C^i}\;=\;\left[\begin{array}{c|c}
            I_{i} & 0\\\hline
            0 & L(\bfa_i,\bfb_i)\end{array}\right]\ee
where $I_i$ is the identity matrix of size $i$, and $L(\bfa_i,\bfb_i)$ is the lower-bidiagonal matrix with corresponding edge sequences $\bfa_i$ and $\bfb_i$. Obviously the matrix $P_{C^i}$ is totally positive in $\Z[\bfa,\bfb]$ equipped with the coefficientwise order.

Since $\scrn_1$ is the concatenation of column transfer matrices we immediately obtain what we refer to as a \emph{production-like factorisation} of $P_{\scrn}$:
\begin{corollary}[Production-like factorisation]\label{cor.prodlikeFactor}
The path matrix $P_{\scrn}$ corresponding to the binomial-like planar network $\scrn$ has the factorisation:
\be P_{\scrn} \;=\; \cdots\left[\begin{array}{c|c}
            I_{2} & 0\\\hline
            0 & L(\bfa_2,\bfb_2)\end{array}\right]\cdot\left[\begin{array}{c|c}
            I_{1} & 0\\\hline
            0 & L(\bfa_1,\bfb_1)\end{array}\right]\cdot L(\bfa_0,\bfb_0)\ee
where $\bfa_i=(\alpha_{i+1,n})_{n\geq 0}$ (with $\alpha_{i+1,0}=1$) and $\bfb_i=(\beta_{i+1,n})_{n\geq0}$.       
\end{corollary}

Recall from the previous subsection that the weighted binomial matrix $B_{x,y}$ is the path matrix of the binomial-like network $\scrn$ with $\alpha_{i,l}=y
$ and $\beta_{i,l}=x$. The above corollary thus yields the following factorisation of $B_{x,y}$:
\begin{corollary}[Production-like factorsation of $B_{x,y}$]\label{cor.binomprodLike}
The matrix $B_{x,y}$ has the production-like factorisation
\be B_{x,y} \;=\; \cdots\left[\begin{array}{c|c}
            I_{2} & 0\\\hline
            0 & L({\bf y},x)\end{array}\right]\cdot\left[\begin{array}{c|c}
            I_{1} & 0\\\hline
            0 & L({\bf y},x)\end{array}\right]\cdot L({\bf y},x)\ee
where $L({\bf y},x)$ denotes the lower-bidiagonal matrix with $(0,0)$-entry 1, all other diagonal entries $y$, all subdiagonal entries $x$, and all other entries $0$.
\end{corollary}

We note that the above corollary imples that $B_{x,y}$ satisfies
\be B_{x,y}\;=\;\left[\begin{array}{c|c}
            1 & 0\\\hline
            0 & B_{x,y}\end{array}\right]L({\bf y},x).\ee
\begin{figure}\centering
\includegraphics[scale=0.5]{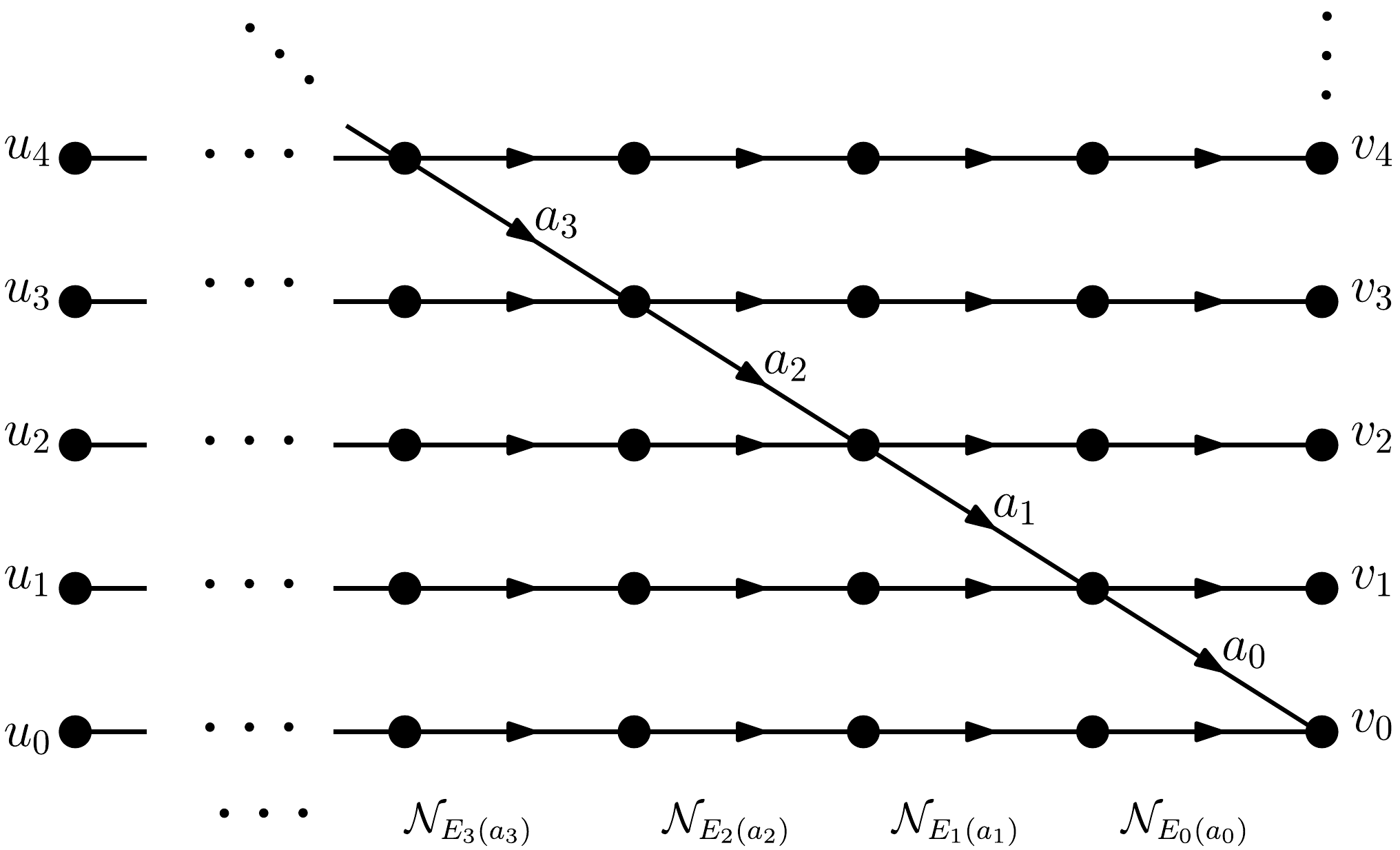}\caption{The planar network $\scrn_{T(\bfa)}$ that has corresponding path matrix $T(\bfa)$.}\label{fig.invbidiag}
\end{figure}

There are other ways to factorise $P_{\scrn}$ that rely on the following helpful definition. Given a sequence $\bfa=(a_n)_{n\geq0}$ let $T(\bfa)\ceq(a^*_{n,k})_{n,k\geq 0}$ denote the lower-triangular matrix with $(n,k)$-entry
\be a^*_{n,k}\;=\;\prod_{j=k}^{n-1}a_j\ee
for $n\geq k$ and 0 in all other cases (we consider the empty product arising from $n=k$ to be $1$). The first few rows and columns of $T(\bfa)$ are thus:
\be T(\bfa)\;=\;\left[
\begin{array}{ccccccc}
 1&  &  &  &  &  &  \\
 a^*_{1,0} & 1&  &  &  &  &  \\
 a^*_{2,0} & a^*_{2,1} & 1&  &  &  &  \\
 a^*_{3,0} & a^*_{3,1} & a^*_{3,2} & 1&  &  &  \\
 a^*_{4,0} & a^*_{4,1} & a^*_{4,2} & a^*_{4,3} & 1&  &  \\
 a^*_{5,0} & a^*_{5,1} & a^*_{5,2} & a^*_{5,3} & a^*_{5,4} & 1&  \\
 \vdots & \vdots & \vdots & \vdots & \vdots & \vdots & \ddots \\
\end{array}
\right].\ee

The matrix $T(\bfa)$ is the path matrix corresponding to the planar network $\scrn_{T(\bfa)}$ in which
\be P_{\scrn_{T(\bfa)}}(u_n\to v_k)\;=\;a^*_{n,k}\ee
(see Figure~\ref{fig.invbidiag}). Clearly $T(\bfa)$ is totally positive in $\Z[\bfa]$ equipped with the coefficientwise order, and by comparing Figure~\ref{fig.bidiag} with Figure~\ref{fig.invbidiag} and~\eqref{eq.elemBidiag} it is easy to see that
\be\label{eq.Taelembi} T(\bfa)\;=\;\cdots E_2(a_2)E_1(a_1)E_0(a_0) = (E_0(-a_0)^{-1}E_1(-a_1)^{-1}E_2(-a_2)^{-1}\cdots)^{-1},\ee
since $E_n(a)^{-1}=E_{n}(-a)$. We thus obtain the following corollary:
\begin{corollary}[Inverse lower-bidiagonal matrices]\label{cor.invbidiag}
For a sequence $\bfa=(a_n)_{n\geq0}$ we have
\be T(\bfa)\;=\;L(-\bfa)^{-1}\ee
where $-\bfa = (-a_n)_{n\geq 0}$.
\end{corollary}
\noindent In light of Corollary~\ref{cor.invbidiag} we refer to the matrix $T(\bfa)$ as an \emph{inverse (lower-)bidiagonal matrix} with corresponding edge sequence $\bfa$. The case where $a_n=x$ for all $n$ often arises in the study of total positivity, and $T(\bfa)$ is then referred to as a \emph{Toeplitz matrix of powers of $x$} since each $(n,k)$ entry is $x^{n-k}$, and we denote it $T_{\infty}(x)$. It follows from the LGV lemma that $T_{\infty}(x)$ is coefficientwise totally positive in $\Z[x]$. More generally, given a sequence $\bfa=(a_n)_{n\geq 0}$ we call the infinite lower-triangular matrix $T_{\infty}(\bfa)=(a_{n-k})_{n,k\geq 0}$ where $a_{l}=0$ for $l<0$ the (infinite) \emph{Toeplitz matrix} associated to $\bfa$.\footnote{We observe here that a sequence $\bfa=(a_n)_{n\geq0}$ of real numbers is \emph{Toeplitz-totally positive} if its associated Toeplitz matrix $\tp(\bfa)$ is totally positive, and such sequences are commonly referred to as \emph{P\'olya frequency sequences}. A sufficient condition for Toeplitz-total positivity of a sequence $\bfa$ of real numbers is given by the celebrated Aissen--Schoenberg--Whitney--Edrei theorem~\cite[Theorem~5.3, p.~412]{Karlin_68}. Similarly, a sequence $\bfa$ of elements belonging to $\R[\bfx]$ is \emph{coefficientwise} Toeplitz-totally positive if its associated Toeplitz matrix $\tp(\bfa)$ is coefficientwise totally positive in $\R[\bfx]$; an extension of the Aissen--Schoenberg--Whitney--Edrei theorem to this more general setting can be found in~\cite[Lemmas~2.4 and~2.5]{Sokal_21f}.}
\begin{figure}\centering
\includegraphics[scale=0.4]{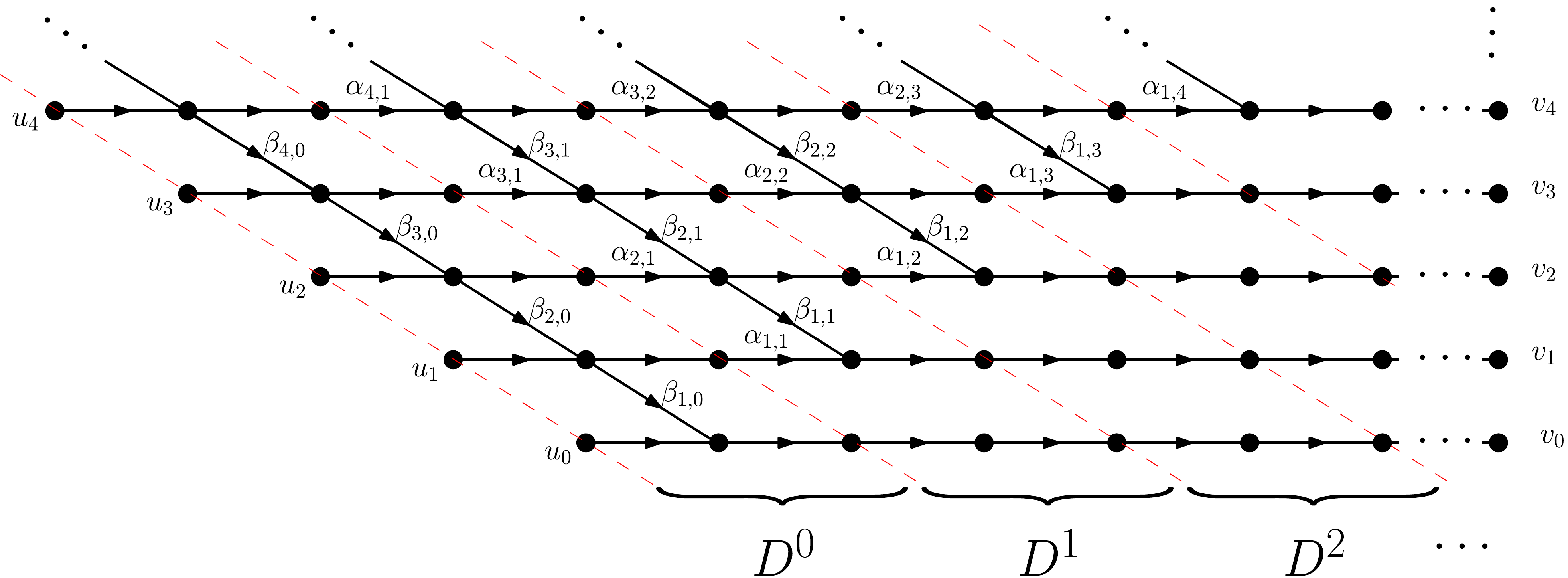}\caption{The planar network $\scrn_2$ up to source $u_4$ and sink $v_4$ (unlabelled edges have weight 1). Please note the red dashed lines are not part of the planar network, they delineate the sources and sinks of each diagonal transfer matrix $D^l$.}\label{fig.diagonaltransfer}\end{figure}

Let us return again to the binomial-like planar network $\scrn$. Observe that $\scrn$ is isomorphic to the planar network $\scrn_2$ in Figure~\ref{fig.diagonaltransfer}, in particular we have
\be P_{\scrn}(u_n\to v_k)\;=\;P_{\scrn_2}(u_n\to v_k)\ee
for all $n,k\geq 0$. 

Consider a subnetwork $D^l$ of $\scrn_2$ contained between a pair of consecutive dashed diagonal lines. The path matrix for each such subnetwork is the \emph{diagonal transfer matrix}:
\be\label{eq.PDi} P_{D^l}\;=\;\left[\begin{array}{c|c}
            I_l & 0\\\hline
            0 & \diag(\bfa_l)T(\bfb_l)\end{array}\right]\ee
where $\bfa_0=(1)_{n\geq 0}$, $\bfa_l=(\alpha_{n+1,l})_{n\geq 0}$ for $l>0$, and $\bfb_l=(\beta_{n+1,l})_{n\geq0}$ for all $l$.

The concatenation of subnetworks $D^0,D^1,\ldots$ from left to right yields a different factorisation of $P_{\scrn},$ namely its \emph{quasi-production-like factorisation}:
\begin{corollary}[Quasi-production-like factorisation]\label{cor.quasiprodlike}
The path matrix $P_{\scrn}$ has the factorisation
\be P_{\scrn}\;=\; T(\bfb_0)\cdot\left[\begin{array}{c|c}
            I_1 & 0\\\hline
            0 & \diag(\bfa_1)T(\bfb_1)\end{array}\right]\cdot\left[\begin{array}{c|c}
            I_2 & 0\\\hline
            0 & \diag(\bfa_2)T(\bfb_2)\end{array}\right]\cdots\ee
where $\bfa_l=(\alpha_{n+1,l})_{n\geq 0}$, and $\bfb_l=(\beta_{n,l})_{n\geq0}$.
\end{corollary}

The weighted binomial matrix $B_{x,y}$ thus has the following quasi-production-like factorisation which can also be found in~\cite[Lemma~2.8]{Gilmore22}:
\begin{corollary}[Quasi-Production-like factorisation of $B_{x,y}$]\label{cor.binomquasiProd}
The weighted binomial matrix $B_{x,y}$ has the quasi-production-like factorisation
\be B_{x,y}\;=\;T(x)\cdot\left[\begin{array}{c|c}
            I_1 & 0\\\hline
            0 & yT(x)\end{array}\right]\cdot\left[\begin{array}{c|c}
            I_2 & 0\\\hline
            0 & yT(x)\end{array}\right]\cdots\ee
where
\be yT(x)\;=\; \diag(y)T(x).\ee
\end{corollary}
Observe that we can write the factorisation of $B_{x,y}$ above as
\be B_{x,y}\;=\;T(x)\diag({\bf y})\cdot\left[\begin{array}{c|c}
            I_1 & 0\\\hline
            0 & T(x)\diag({\bf y})\end{array}\right]\cdot\left[\begin{array}{c|c}
            I_2 & 0\\\hline
            0 & T(x)\diag({\bf y})\end{array}\right]\cdots\ee
where ${\bf y}=(1,y,y,\ldots)$, from which it immediately follows that
\be B_{x,y}\;=\;T(x)\diag({\bf y})\cdot\left[\begin{array}{c|c}
            1 & 0\\\hline
            0 & B_{x,y}\end{array}\right].\ee

We conclude this subsection with one final lemma relating bidiagonal and inverse bidiagonal matrices.

\begin{lemma}\label{lem.switch}
Suppose $\bfa=(a_n)_{n\geq 0}$ and $\bfb=(b_n)_{n\geq 0}$ are sequences of elements belonging to a field $F$. Then
\begin{enumerate}[(i)]
\item If $a_n=-b_n$ for all $n$ then $L(\bfa)T(\bfb)=I$;
\item If $a_n+b_n\neq 0$ for all $n$ then
\be\label{eq.switch} L({\bfa})T({\bfb}) = \left[\begin{array}{c|c}
          1 & 0 \\
          \hline
          0 & T({\bfb'})
          \end{array}\right]\cdot L({\bfa'})\ee
where $\bfa' = (a'_n)_{n\geq 0}$ is the edge sequence in which 
\be a'_n\ceq \begin{cases}a_0+b_0 &\textrm{if } n=0,\\ \frac{a_{n-1}(a_n+b_n)}{a_{n-1}+b_{n-1}} &\textrm{if }n>0,\end{cases}\ee
and $\bfb'=(b'_n)_{n\geq 0}$ is the edge sequence in which
\be b'_n = \frac{b_{n}(a_{n+1}+b_{n+1})}{a_n+b_n}.\ee
\end{enumerate}
\end{lemma}
\begin{proof}
$(i)$ If $a_n=-b_n$ then we trivially have
\be L(\bfa)T(\bfb)\;=\; L(-\bfb)T(\bfb)\;=\;L(-\bfb)L(-\bfb)^{-1}\;=\;I,\ee

$(ii)$ Suppose $a_n+b_n\neq0$ for all $n$, and
let $M=L(\bfa)T(\bfb)$ and 
\be N\;=\;\left[\begin{array}{c|c}
          1 & 0 \\
          \hline
          0 & T(\bfb')
          \end{array}\right]\cdot L({\bfa'}).\ee
Clearly $M_{n,n}=1=N_{n,n}$, while for $n>k\geq 0$ the matrix entries of $M_{n,k}$ can be written as a telescoping product
\be M_{n,k}\;=\;(a_{n-1}+b_{n-1})\prod_{i=k}^{n-2}b_i\;=\;(a_{k}+b_{k})\prod_{i=k}^{n-2}\frac{b_i(a_{i+1}+b_{i+1})}{(a_i+b_i)}\;=\;N_{n,k}\ee
(the empty product is taken to be $1$).

\end{proof}
\begin{figure}\centering
\includegraphics[scale=0.36]{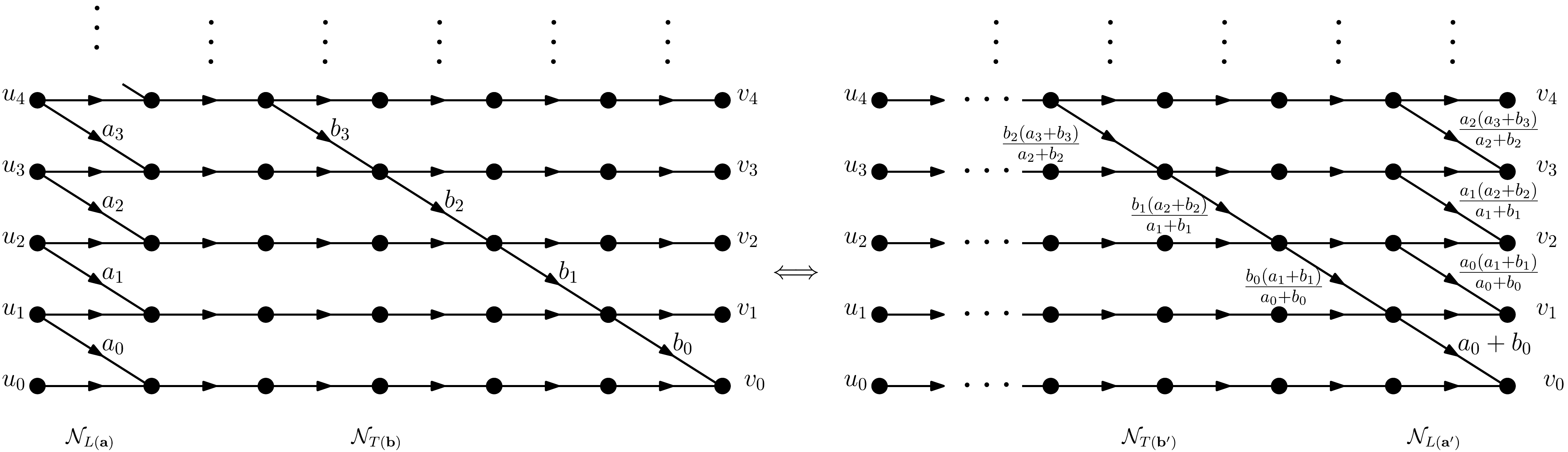}\caption{The planar network representation of Lemma~\ref{lem.switch}}\label{fig.switch}\end{figure}

Figure~\ref{fig.switch} is a planar network representation of Lemma~\ref{lem.switch}. As simple as Lemma~\ref{lem.switch} is, it will prove to be a fundamental tool in our proof of Theorem~\ref{thm.fqtqCoeff}. In applications in this paper the field $F$ will be the field $\Q(q)$ of rational functions of $q$, which is the fraction field of the ring $\Z[q]$ of polynomials in $q$. The planar network we present in Section~\ref{sec.pointwise} arises from specialising the weights $\balph$ and $\bbet$ of the binomial-like planar network $\scrn$ to certain rational functions of $q$, and Lemma~\ref{lem.switch} allows us to transform this network into a different network with weights that are polynomials in $q$.

\subsection{A remark on Neville elimination}
Given a totally positive matrix $M$ there are myriad ways in which to express $M$ as a product of totally positive matrices, however, if the entries of $M$ belong to $\R$ then Gasca and Pe\~na~\cite{Gasca92} provide an algorithm for systematically determining a ``canonical'' factorisation. The process is based on the Neville-Aitken technique~\cite{Gasca87} which, in the context of solutions of linear systems, gives rise to \emph{Neville elimination}~\cite{Gasca87a,Gasca87b}. Neville elimination can be a powerful tool for studying totally positive matrices in general; here we sketch a brief overview of it purely for unit-lower-triangular matrices. For a full treatment please see~\cite{Gasca92}.

Let $M\ceq(m_{n,k})_{0\leq n,k\leq N}$ be a unit-lower-triangular matrix of size $N+1$. If $M$ contains any zeroes in its initial $k=0$ column then form the matrix $M_0\ceq(m^0_{n,k})_{0\leq n,k\leq N}$ by moving the offending rows of $M$ to the bottom in such a way that the relative order among them is preserved (if $M$ does not contain any zeroes in its initial column then set $M_0\ceq(m^0_{n,k})_{0\leq n,k\leq N}$). The entries $m^0_{i,0}$ for $i\geq 0$ are referred to as the $(i,0)$ \emph{pivots} of $M$.

Suppose row $n_1$ is the bottom-most nonzero entry in column $k=0$. For $i$ decreasing from $n_1$ to 1 subtract
\be \frac{m^0_{i,0}}{m^0_{i-1,0}}(\textrm{row  } (i-1))\ee
from row $i$ (thereby reducing $m^0_{i,0}$ to $0$ from bottom to top as $i$ decreases). This operation is equivalent to left-multiplying $M$ by a product of elementary lower-bidiagonal matrices, and can thus be expressed algebraically as
\be E_{N+1,1}(-\rho_{1,0})E_{N+1,2}(-\rho_{2,0})\cdots E_{N+1,N}(-\rho_{N,0})R_{N+1}M\;=\;\left[\begin{array}{c|c}1&0\\\hline0&M_1\end{array}\right],\ee
where: $E_{N+1,n}(x)$ is the elementary bidiagonal matrix of size $N+1$ with $(n,n-1)$-entry $x$, $1$s on the diagonal, and $0$s everywhere else; $R_{N+1}$ is some matrix of size $N+1$ that encodes a permutation of rows; and 
\be \rho_{i,0}\;=\;\begin{cases}\frac{m^0_{i,0}}{m^0_{i-1,0}}&\textrm{if }m^0_{i,0}\neq 0,\\
                                0 &\textrm{if }m^0_{i,0}=0\end{cases}\ee
is a ratio of pivots of $M$.

The entries $m^1_{i,0}$ in the $k=0$ column of $M_1=(m^1_{n,k})_{0\leq n,k\leq N-1}$ are the $(i,1)$ pivots of $M$. Reducing the zeroth column of $M_1$ to $0$s everywhere except on the diagonal in the manner described above (that is, rearranging rows and successively subtracting them from each other) yields
\begin{multline} \left[\begin{array}{c|c}
            1 & 0\\\hline
            0 & E_{N,1}(-\rho_{1,1})E_{N,2}(-\rho_{2,1})\cdots E_{N,N-1}(-\rho_{N-1,1})R_N\end{array}\right]\\\cdot E_{N+1,1}(-\rho_{1,0})E_{N+1,2}(-\rho_{2,0})\cdots E_{N+1,N}(-\rho_{N,0})R_{N+1}M\;=\;\left[\begin{array}{c|c}I_2&0\\\hline0&M_2\end{array}\right]\end{multline}
where $M_2=(m^2_{n,k})_{0\leq n,k\leq N-2}$ is a matrix of size $N-1$ and
\be \rho_{i,1}\;=\;\begin{cases}\frac{m^1_{i,0}}{m^1_{i-1,0}}&\textrm{if }m^1_{i,0}\neq0,\\
                                0 &\textrm{if }m^1_{i,0}=0\end{cases}\ee
is a ratio of pivots of $M_1$.

Proceeding iteratively on smaller and smaller matrices we obtain, after a finite number of steps, the identity
\begin{multline}
\left[\begin{array}{c|c}
            I_{N-1} & 0\\\hline
            0 & E_{2,1}(-\rho_{1,N-1})\end{array}\right]\cdot\left[\begin{array}{c|c}
            I_{N-2} & 0\\\hline
            0 & E_{3,1}(-\rho_{1,N-2})E_{3,2}(-\rho_{2,N-2})R_{3}\end{array}\right]\\\cdots E_{N+1,1}(-\rho_{1,0})E_{N+1,2}(-\rho_{2,0})\cdots E_{N+1,N}(-\rho_{N,0})R_{N+1}M\;=\;I_{N+1},\end{multline}
equivalently,
\begin{multline}\label{eq.nevFac1}
R_{N+1}^{-1}E_{N+1,N}(\rho_{N,0})\cdots E_{N+1,2}(\rho_{2,0})E_{N+1,1}(\rho_{1,0})\cdot \left[\begin{array}{c|c}
            1 & 0\\\hline
            0 & R_N^{-1}E_{N,N-1}(\rho_{N-1,1})\cdots E_{N,1}(\rho_{1,1})\end{array}\right]\\\cdots\left[\begin{array}{c|c}
            I_{N-2} & 0\\\hline
            0 & R_3^{-1}E_{3,2}(\rho_{2,N-2})E_{3,1}(\rho_{1,N-2})\end{array}\right]\left[\begin{array}{c|c}
            I_{N-1} & 0\\\hline
            0 & E_{2,1}(\rho_{1,N-1})\end{array}\right]\;=\;M.
\end{multline}
where each $\rho_{i,j}$ is a ratio of pivots obtained at each step of the algorithm:
\be\rho_{i,j}\;=\;\frac{m^j_{i,0}}{m^j_{i-1,0}}.\ee
We call the factorisation of $M$ obtained in this way the \emph{Neville factorisation} of $M$.

Gasca and Pe\~na showed in~\cite{Gasca92} that the matrix $M$ is totally positive if and only if $R_j=I_j$ for all $j$ (that is, Neville elimination can be applied to $M$ without ever interchanging rows at any stage of the algorithm), and the pivots of $M$ are all nonnegative (the same result can also be found in Chapter~6 of Pinkus' book~\cite{Pinkus_10}, although Neville elimination is not treated explicitly there).

\begin{figure}
\includegraphics[scale=0.5]{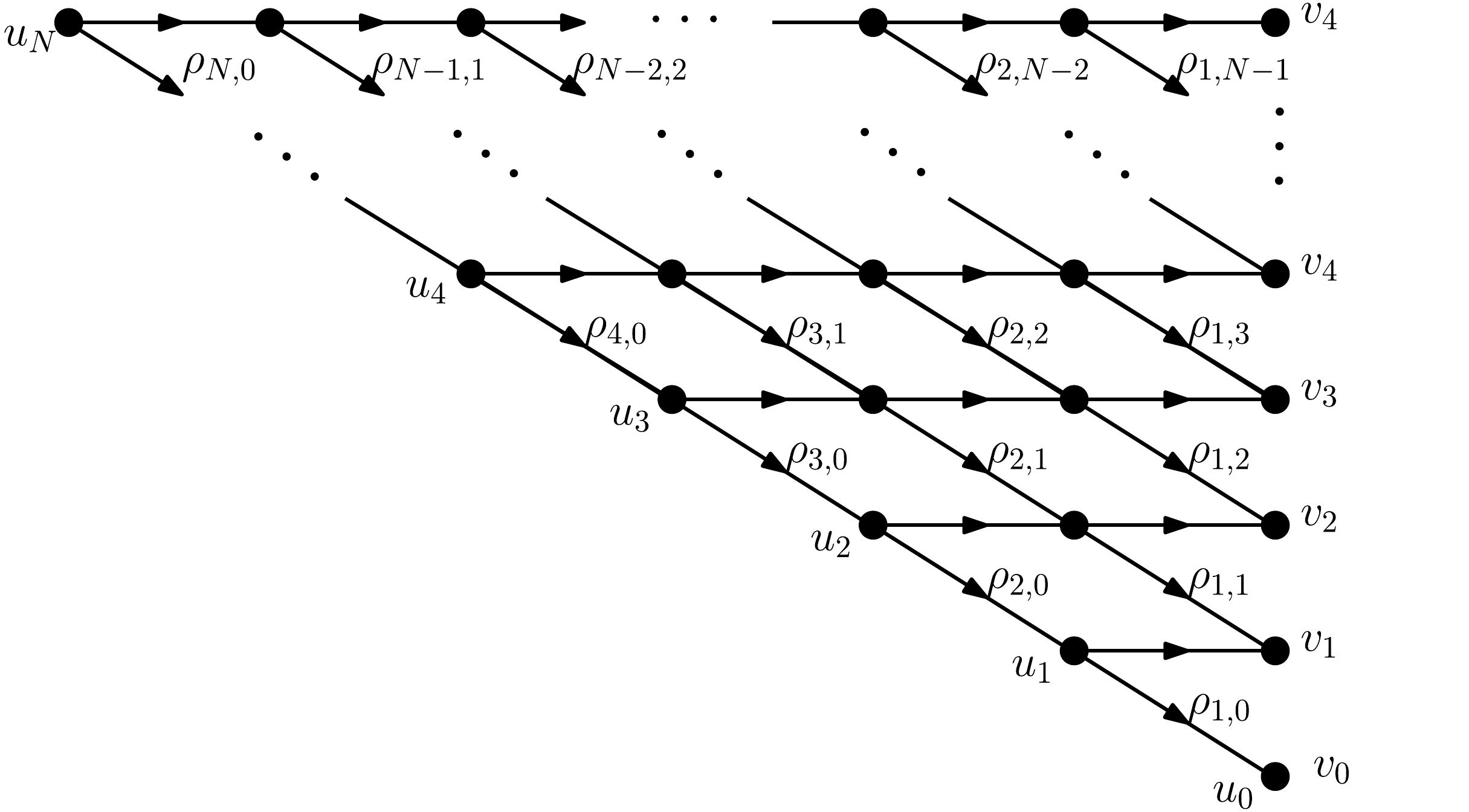}\caption{The binomial-like planar network obtained using Neville elimination.}\label{fig.binomNev}\end{figure}

We can translate the Neville factorisation of $M$ into the language of planar networks. If $M$ is totally positive then it follows from~\eqref{eq.nevFac1} that
\begin{multline}\label{eq.nevFac2}
E_{N+1,N}(\rho_{N,0})\cdots E_{N+1,2}(\rho_{2,0})E_{N+1,1}(\rho_{1,0})\cdot \left[\begin{array}{c|c}
            1 & 0\\\hline
            0 & E_{N,N-1}(\rho_{N-1,1})\cdots E_{N,1}(\rho_{1,1})\end{array}\right]\\\cdots\left[\begin{array}{c|c}
            I_{N-2} & 0\\\hline
            0 & E_{3,2}(\rho_{2,N-2})E_{3,1}(\rho_{1,N-2})\end{array}\right]\left[\begin{array}{c|c}
            I_{N-1} & 0\\\hline
            0 & E_{2,1}(\rho_{1,N-1})\end{array}\right]\;=\;M,
\end{multline}
where every pivot $\rho_{i,j}$ is nonnegative. We can easily construct a planar network representation $\scrn_M$ of this factorisation (see Figure~\ref{fig.binomNev}); the network we obtain is simply the standard binomial-like planar network with $\alpha_{n,k}=1$ and $\beta_{n,k}=\rho_{n,k}$. We call the planar network $\scrn_M$ obtained using Neville elimination the \emph{Neville network} corresponding to $M$. 

We have already seen how to interpret such planar networks as production-like and quasi-production-like factorisations. By reading the network as a product of column transfer matrices we obtain the production-like factorisation
\be\label{eq.nevfacprod} M\;=\;\left[\begin{array}{c|c}
            I_{N-1} & 0\\\hline
            0 & L_2(\bfb'_{N-1})\end{array}\right]\cdots\left[\begin{array}{c|c}
            1 & 0\\\hline
            0 & L_N(\bfb'_1)\end{array}\right]L_{N+1}(\bfb'_0)\ee
where $L_{N+1-i}(\bfb'_{i})$ is the finite lower-bidiagonal matrix of size $N+1-i$ with corresponding finite edge sequence
\be \bfb'_i\;=\;(\rho_{i+1,n})_{0\leq n\leq N-1-i}.\ee
Similarly, by reading the network as a product of diagonal transfer matrices we obtain the quasi-production-like factorisation
\be\label{eq.nevFac3} M\;=\; T_{N+1}(\bfb_0)\left[\begin{array}{c|c}
            1 & 0\\\hline
            0 & T_N(\bfb_1)\end{array}\right]\cdots \left[\begin{array}{c|c}
            I_{N-1} & 0\\\hline
            0 & T_2(\bfb_{N-1})\end{array}\right].\ee
where
\be \bfb_j\;=\;\left(\rho_{n,j}\right)_{1\leq n\leq N-j}.\ee

We note that the quasi-production-like factoristion above can be obtained directly from the products of elementary bidiagonal matrices in~\eqref{eq.nevFac2}, since according to~\eqref{eq.Taelembi} we have
\be\left[\begin{array}{c|c}
            I_j & 0\\\hline
            0 & E_{N-j+1,N-j}(\rho_{N-j,j})\cdots E_{N-j+1,1}(\rho_{1,j})\end{array}\right]\;=\;\left[\begin{array}{c|c}
            I_j & 0\\\hline
            0 & T_{N-j+1}(\bfb_j)\end{array}\right].\ee
The production-like factorisation, on the other hand, can be obtained from~\eqref{eq.nevFac2} by commuting elementary bidiagonal matrices. We have
\be\left[\begin{array}{c|c}
            I_i & 0\\\hline
            0 & E_{N-i+1,N-i}(\rho_{N-i,i})\cdots E_{N-i+1,1}(\rho_{1,i})\end{array}\right]\;=\;E_{N+1,N}(\rho_{N-i,i})\cdots E_{N+1,i+1}(\rho_{1,i}),\ee
and since $E_{N,n_1}(x)E_{N,n_2}(x)=E_{N,n_2}(x)E_{N,n_1}(x)$ for $n_2>n_1+1$ it follows that we can commute elementary bidiagonal matrices in~\eqref{eq.nevFac2}, thereby obtaining
\begin{multline}\label{eq.nevFac4}
[E_{N+1,N}(\rho_{N,0})\cdots E_{N+1,2}(\rho_{2,0})E_{N+1,1}(\rho_{1,0})]\cdot [E_{N+1,N}(\rho_{N-1,1})\cdots E_{N+1,2}(\rho_{1,1})]\cdots\\ [E_{N+1,N}(\rho_{2,N-2})E_{N+1,N-1}(\rho_{1,N-2})][E_{N+1,N}(\rho_{1,N-1})]\\\;=\;[E_{N+1,N}(\rho_{N,0})][(E_{N+1,N-1}(\rho_{N-1,0})E_{N+1,N}(\rho_{N-1,1})]\\\cdots[E_{N+1,1}(\rho_{1,0})E_{N+1,2}(\rho_{1,1})\cdots E_{N+1,N}(\rho_{1,N-1})].
\end{multline}
The right-hand side agrees with that of~\eqref{eq.nevfacprod} since
\be [E_{N+1,i+1}(\rho_{i+1,0})E_{N+1,i+2}(\rho_{i+1,1})\cdots E_{N+1,N}(\rho_{i+1,N-i-1})]\;=\;\left[\begin{array}{c|c}
            I_{i} & 0\\\hline
            0 & L_{N-i+1}(\bfb'_{i})\end{array}\right].\ee

Neville elimination can easily be applied to a unit-lower-triangular matrix $M$ with entries that belong to the polynomial ring $\Z[\bfx]$, though the resulting factorisation may well consist of matrices with entries that belong to the field $\Q(\bfx)$ of rational functions of $\bfx$. The planar network $\scrn_1'$ described in Section~\ref{sec.pointwise} below (see Figure~\ref{fig.binomialqrat}) that we use to prove the pointwise total positivity of the $q$-forests matrix is the Neville network for the matrix $\fm'(q)\;=\;(f_{n,k}(q))_{n,k\geq 1}$. By transforming this Neville network into a different planar network with weights that are polynomials in $q$ we eventually show in Section~\ref{sec.coeffpos} that $\fm(q)$ is \emph{coefficientwise} totally positive in $\Z[q]$.

Determining binomial-like planar networks via Neville elimination can be incredibly useful if the set of pivots are of a ``nice'' form from which a general pattern is easy to guess; however, often when attempting to construct planar networks for suspected totally positive matrices\footnote{Such as, for example, the maddeningly stubborn \emph{Eulerian triangle}, which was conjectured by Brenti~\cite{Brenti_96} to be totally positive over a quarter of a century ago.} one finds oneself with a set of rational expressions for the pivots that are difficult to make any sense of at all. Luckily, in the case of the forests and the trees matrices the pivots were of a general form that was easy to understand.

We have now established all of the fundamental concepts required to prove our main result (Theorem~\ref{thm.fqtqCoeff}). In the next section we will interpret the entries of the $q$-forests matrix and the $q$-trees matrix combinatorially, and prove some matrix identities relating them.

\section{Combinatorial interpretations of the entries of the \texorpdfstring{$q$}{TEXT}-forests and \texorpdfstring{$q$}{TEXT}-trees matrices, and some identites relating them}\label{sec.combin.interp}{}
The entries of the matrices $\fm(q)$ and $\tm(q)$ are polynomials in $q$ that count trees or forests according to some statistic, and it is natural to try to interpret what that statistic might be. We present one such interpretation in the following two subsections, before considering some identities relating $\fm(q)$ and $\tm(q)$.

\subsection{A combinatorial interpretation of the entries of \texorpdfstring{$\fm(q)$}{TEXT}}
We begin with the $q$-forests matrix
\be \fm(q)\;=\;\left(\qbin{n-1}{k-1}(\qn{n})^{n-k}\right)_{n,k\geq 0}\ee
which counts forests on $n$ vertices with $k$-components.
Let $\scrf_{n,k}$ denote the set of such forests for given $n,k$ and consider a forest $F\in\scrf_{n,k}$. The vertices of $F$ can be partitioned into three subsets:
\be V_r(F)\cup V_s(F)\cup V_t(F)\ee
where 
\begin{eqnarray*} V_r(F)&\ceq&\{v\in F: v \textrm{ is a root of a component}\},\\
                V_s(F)&\ceq&\{v\in F: v \textrm{ is a lowest-numbered child of a root}\},\\
                V_t(F)&\ceq&\{v\in F:v \textrm{ is neither a root, nor a lowest-numbered child of a root}\}.\end{eqnarray*}
We now describe one way in which to assign weights to the vertices in $V_r(F),V_s(F)$, and $V_t(F)$.

For each vertex $v$ in $V_r(F)$ define
\be V_r(F,v)\ceq\{v':v'\notin\scrf_r,v'<v\}\ee
to be the set of nonroot vertices in $F$ that are lower-numbered than a root vertex $v$.
Each root $v\in V_r(F)$ is then given the weight
\be\wt(v)\;=\;r^{|V_r(F,v)|}\ee
and for a forest $F\in\scrf_{n,k}$ define\footnote{The portmanteau ``niblings" is a combination of nephews/nieces and siblings.}
\be\nibl(F)\ceq \sum_{v\in V_r(F)}|V_r(F,v)|.\ee

For the vertex $v_{\min}$ (if it exists), we consider the roots $v_0,v_1,\ldots,v_{k-1}\in V_r(F)$ where
\be v_{0}<v_{1}<\cdots<v_{k-1}.\ee
The vertex $v_{\min}$ is the child of one of $v_{0},\ldots,v_{k-1}$, so we assign $v_{\min}$ the weight
\be \wt(v_{\min})\;=\;s^{j},\ee
where $v_{j}$ is the parent of $v_{\min}$ (we say that $j$ is the \emph{smallest child index} of $F$ and denote it $\scind(F)$).

Lastly, for each vertex $v$ in $V_t(F)$ set
\be \wt(v)\;=\;t^{\parent{v}-1},\ee
where $\parent{v}$ denotes the label of the parent of $v$, and define
\be \parents(F)\ceq\sum_{v\in V_t(F)}\parent{v}.\ee

Putting the above together in one place we have for a vertex $v$ in $F$:
\be\label{eq.vertexwt} \wt(v)\ceq\begin{cases}r^{|V_r(F,v)|} & v\in V_r(F),\\
                           s^{\scind(F)} & v=v_{\min},\\
                           t^{\parent{v}-1} & v\in V_t(F).\end{cases}\ee
We then define the weight of a forest to be the product of the weights of its vertices:
\be\label{eq.forestwt}\wt(F)\ceq\prod_{v\in F}\wt(v)\;=\;r^{\nibl(F)}s^{\scind(F)}t^{\parents(F)-n+k},\ee
and note that if $n=k$ we have $\wt(F)=1$.

\begin{proposition}\label{prop.niblings}
For the set $\scrf_{n,k}$ of forests on $n$ vertices comprised of $k$ components with $n>k$ we have
\be\sum_{F\in\scrf_{n,k}}r^{\nibl(F)}s^{\scind(F)}t^{\parents(F)-n+k}\;=\;\pqbin{n}{k}{r}\pqn{k}{s}(\pqn{n}{t})^{n-k-1}.\ee
\end{proposition}
Our proof below is an extension of the first proof of Proposition~5.3.2 given by Stanley in~\cite{Stanley_86}, and we are very grateful to Bishal Deb for pointing it out.
\begin{proof}
A forest $F\in\scrf_{n,k}$ consists of a root set $V_r(F)$, and a subforest attached to the roots consisting of the vertices $V_t(F)\cup V_{s}(F)$. We construct a sequence $\sigma_1,\sigma_2,\ldots,\sigma_{n-k+1}$ of subforests of $F$ (all with root set $V_{r}(F)$) in the following way: set $\sigma_1=F$. If $i<n-k+1$ and $\sigma_i$ is defined then let $\sigma_{i+1}$ be the subforest obtained from $\sigma_i$ by removing its largest nonroot endpoint $v_i$ (together with the edge incident to it). Let $p_i$ be the unique vertex of $\sigma_i$ adjacent to $v_i$ and consider the \emph{Pr\"ufer sequence} (or \emph{Pr\"ufer code}, see~\cite{Prufer18}) of the subforest of $F$ that arises from removing all vertices in $V_t(F)$ in this way, 
\be\gamma(\scrf)\ceq(p_1,p_2,\ldots,p_{n-k}).\ee
Note that for $i<n-k$, $p_i\in[n]$, while $p_{n-k}\in V_r(F)$, so the number of such sequences is $kn^{n-k-1}$.

Now let $\scrf_{n,V_r}$ denote the set of forests with a specified root set $V_r$ of size $k$. The map from forests to Pr\"ufer sequences described above
\be \gamma:\scrf_{n,V_r}\to[n]^{n-k-1}\times V_r\ee
is a bijection (see~\cite[Proposition~5.3.2]{Stanley_86}).\footnote{Please note that here $[n]^{n-k-1}$ denotes the cartesian product of the set $[n]$ with itself $n-k-1$ times.} Our aim is to understand the contribution of the weights from the vertices in terms of these Pr\"ufer sequences and root sets.

Each subforest on the vertices $V_t(F)\cup V_s(F)$ has a unique Pr\"ufer sequence 
\be(p_1, p_2, \ldots, p_{n-k-1}, p_{n-k}),\ee and the first $n-k-1$ elements of the sequence correspond to vertices with weights $(t^{p_1-1},t^{p_2-1},\ldots,t^{p_{n-k-1}-1})$. The weight of all sequences $(p_1,p_2,\ldots,p_{n-k-1})$ is thus $(\pqn{n}{t})^{n-k-1}$. It is easy to see that the element $p_{n-k}$ is the parent of $v_{\min}$, and since $v_{\min}$ has weight $s^l$ for $0\leq l\leq k-1$ (according to the definition of the weight function), the weight of the set of subforests with root set $V_r$ is
\be\pqn{k}{s}(\pqn{n}{t})^{n-k-1}.\ee  What remains is to show that the weight of all possible root sets is 
\be\pqbin{n}{k}{r}.\ee
Given a set of $k$ root vertices $V_r$ chosen from $[n]$, form the word $w=w_1w_2\ldots w_n$ of length $n$ on the alphabet $\{0,1\}$, containing 0s at positions $v_1,\ldots,v_k\in V_r$ and 1s everywhere else. The weight we assign to each vertex $v_i\in V_r$ then corresponds to the number of $1$s preceding the $0$ at position $w_{v_i}$, so counting root sets with the weighting specified above is equivalent to counting all words $w$ with respect to the number of inversions in $w$ (an inversion is a pair $(i,j)$ such that $w_i=1$, $w_j=0$, and $i<j$). This is a well-known interpretation of the $r$-binomial coefficient:
\be\pqbin{n}{k}{r}.\ee
\end{proof}

The above proposition might invite one to consider the more general matrix \be\fm(r,s,t)\ceq(f_{n,k}(r,s,t))_{n,k\geq 0}\ee with $(n,k)$-entry given by
\be f_{n,k}(r,s,t)\;=\;\pqbin{n}{k}{r}\pqn{k}{s}(\pqn{n}{t})^{n-k-1},\ee
but alas, the first few rows of this matrix are:
\be\left[
\begin{array}{ccccc}
 1 &  &  &  &\\
 0 & 1 &  &  &\\
 0 & r+1 & 1 &  &\\
 0 & \left(r^2+r+1\right) \left(t^2+t+1\right) & \left(r^2+r+1\right) (s+1) & 1& \\
 \vdots & \vdots& \vdots & \vdots &\ddots \\
\end{array}
\right]\ee
which is not even coefficientwise TP$_2$ in $\Z[r,s,t]$ since
\begin{multline} f_{2,1}(r,s,t)f_{3,2}(r,s,t)-f_{3,1}(r,s,t)f_{2,2}(r,s,t)\\=\;(1 + r + r^2) (s + r (1 + s) - t (1 + t))\not\myge0\end{multline}
In fact it seems the only way to ensure that $\fm(r,s,t)$ is coefficientwise totally positive is to set $r=s=t=q$, thereby reducing $\fm(r,s,t)$ to the $q$-forests matrix $\fm(q)$. We conclude:

\begin{corollary}\label{cor.niblingsq}
For the $q$-forests matrix $\fm(q)=(f_{n,k}(q))_{n,k\geq 0}$ we have
\be f_{n,k}(q)\;=\;\sum_{F\in\scrf_{n,k}}q^{\nibl(F)+\scind(F)+\parents(F)-n+k},\ee
where $\scrf_{n,k}$ denotes the set of forests on $n$ vertices with $k$ components.
\end{corollary}

Having established this combinatorial interpretation of the entries of $\fm(q)$ we now turn to the $q$-trees matrix.

\subsection{A combinatorial interpretation of the entries of \texorpdfstring{$\tm(q)$}{TEXT}}
The entries of the $q$-trees matrix
\be \tm(q)\;=\;\left(\qbin{n}{k}(\qn{n})^{n-k}\right)_{n,k\geq 0}\ee
count rooted trees on $n+1$ vertices with $k$ children smaller than the root with respect to some statistic. In~\cite[Lemma~2]{Chauve99}, however, the authors give a bijection between the set of all such trees and the set of forests $\scrf_{n+1,k+1}^*$ on the vertex set $[n+1]$ comprised of $k+1$ components where $n+1$ is a leaf (note that a root with no children is also a leaf); we therefore interpret the entries of $\tm(q)$ as enumerating forests in $\scrf_{n+1,k+1}^*$ with respect to certain statistics.

We can can construct the set $\scrf_{n+1,k+1}^*$ in the following way: either take the set $\scrf_{n,k+1}$ of forests on $n$ vertices with $k+1$ components, and for each forest $F\in\scrf_{n,k+1}$ attach the vertex $n+1$ as a leaf to any of the vertices in $F$ (there are $n$ ways to do this); or take the set of forests $\scrf_{n,k}$ on $n$ vertices with $k$ components and to each $F\in\scrf_{n,k}$ attach the vertex $n+1$ as a singleton component. It follows that
\be|\scrf_{n+1,k+1}^*|\;=\;n|\scrf_{n,k+1}|+|\scrf_{n,k}|\;=\;n\binom{n-1}{k}n^{n-k-1}+\binom{n-1}{k-1}n^{n-k}\;=\;\binom{n}{k}n^{n-k}.\ee
If we weight the forests in $\scrf_{n+1,k+1}^*$ according to~\eqref{eq.forestwt} in the previous subsection then it is not difficult to see that
\begin{multline}\sum_{F\in\scrf_{n+1,k+1}^*}r^{\nibl(F)}s^{\scind(F)}t^{\parents(F)-n+k}\;=\;\pqbin{n}{k+1}{r}\pqn{k+1}{s}(\pqn{n}{t})^{n-k-1}\\+r^{n-k}\pqbin{n}{k}{r}\pqn{k}{s}(\pqn{n}{t})^{n-k-1}\end{multline}
since the weight of the set of forests where $n+1$ is a leaf, but not a component, is 
\be [n]_t\sum_{F\in\scrf_{n,k+1}}r^{\nibl(F)}s^{\scind(F)}t^{\parents(F)-n+k+1}\;=\;\pqbin{n}{k+1}{r}\pqn{k+1}{s}(\pqn{n}{t})^{n-k-1}\ee
(we obtain a multiplicative factor of $t^{\parent{n+1}-1}$ each time we attach $n+1$ as a leaf to a leaf of a component, and each possible parent of $n+1$ is an element of $[n]$ in a forest $F\in\scrf_{n,k+1}$)
and the weight of the set of forests where $n+1$ is a singleton component is 
\be r^{n-k}\sum_{F\in\scrf_{n,k}}r^{\nibl(F)}s^{\scind(F)}t^{\parents(F)-n+k}\;=\;r^{n-k}\pqbin{n}{k}{r}\pqn{k}{s}(\pqn{n}{t})^{n-k-1}\ee
(since all nonroot vertices of $F\in\scrf_{n,k}$ are smaller than $n+1$, and this corresponds to appending a $0$ to the word obtained from the root set $V_r(F)$ described in the proof of Proposition~\ref{prop.niblings}). Letting 
\be \tm(r,s,t)\ceq(t_{n,k}(r,s,t))_{n,k\geq 0}\;=\;([n]_t)^{n-k}\left(r^{n-k}\pqbin{n}{k}{r}[k]_s+[k+1]_s\pqbin{n}{k+1}{r}\right)\ee
we have:
\begin{lemma}
The entries of $\tm(r,s,t)$ satisfy
\be t_{n,k}(r,s,t)\;=\;\sum_{F\in\scrf_{n+1,k+1}^*}r^{\nibl(F)}s^{\scind(F)}t^{\parents(F)-n+k}.\ee
\end{lemma}
The first few rows of $\tm(r,s,t)$ are:
\be{\tiny\left[
\begin{array}{ccccc}
 1 &  &  &  &  \\
 1 & 1 &  &  &  \\
 (r+1) (t+1) & r^2+r+s+1 & 1 &  &  \\
 \left(r^2+r+1\right) \left(t^2+t+1\right)^2 & \left(r^2+r+1\right) \left(t^2+t+1\right) \left(r^2+s+1\right) & r^3 s+r^2 s+r^3+r^2+r s+r+s^2+s+1 & 1 &  \\
 \vdots & \vdots & \vdots & \vdots & \ddots 
\end{array}
\right]}\ee
which contains the $2\times 2$ minor
\be t_{1,0}(r,s,t)t_{2,1}(r,s,t)-t_{2,0}(r,s,t)t_{1,1}(r,s,t)\;=\;r^2 + s - t - r t\not\myge0,\ee
so $\tm(r,s,t)$ is not even coefficientwise TP$_2$ in $\Z[r,s,t]$. Again, it seems the only way to restore total positivity is to specialise $r=s=t=q$, in which case
\begin{multline} t_{n,k}(q,q,q)\;=\;([n]_q)^{n-k-1}\left(q^{n-k}\pqbin{n}{k}{q}[k]_q+[k+1]_q\pqbin{n}{k+1}{q}\right)\\\;=\qbin{n}{k}([n]_q)^{n-k-1}(q^{n-k}[k]_q+[n-k]_q)\;=\;([n]_q)^{n-k}\qbin{n}{k}\\\;=\;t_{n,k}(q).\end{multline} We thus have:
\begin{corollary}\label{cor.niblingsqT}
The entries of the $q$-trees matrix $\tm(q) = (t_{n,k}(q))_{n,k\geq 0}$ satisfy
\be t_{n,k}(q)\;=\;\sum_{F\in\scrf_{n+1,k+1}^*}q^{\nibl(F)+\scind(F)+\parents(F)-n+k}\;=\;\qbin{n}{k}([n]_q)^{n-k}\ee
where $\scrf_{n+1,k+1}^*$ denotes the set of forests on $n+1$ vertices with $k+1$ components, in which the vertex $n+1$ is a leaf.
\end{corollary}

It is safe to say that the statistics defined above are quite unconventional, arising from inserting $q$-weights into the Pr\"ufer sequence for trees and forests. In~\cite{Egecioglu86} the authors provide an alternative method for encoding weighted trees that is more refined than that of Pr\"ufer; a closer examination of this weight-preserving bijection in relation to $q$-forests and $q$-trees matrices will form part of future work on these matrices. For the time being we will turn our attention to some identities relating the matrices $\fm(q)$ and $\tm(q)$.

\subsection{Some identities involving \texorpdfstring{$\fm(q)$}{TEXT} and \texorpdfstring{$\tm(q)$}{TEXT}}\label{subsec.identities}
We have already seen in Section~\ref{sec.intro} that the entries of the forests matrix
\be \fm\ceq(f_{n,k})_{n,k\geq 0}\;=\;\left(\binom{n-1}{k-1}n^{n-k}\right)_{n,k\geq 0}\ee
and the entries of the trees matrix
\be \tm\ceq(t_{n,k})_{n,k\geq 0}\;=\;\left(\binom{n}{k}n^{n-k}\right)_{n,k\geq 0}\ee
satisfy
\be f_{n,k}\;=\;\frac{k}{n}t_{n,k}\ee
for $n>0$. We therefore have:
\begin{lemma}\label{lem.transF}
The matrices $\fm$ and $\tm$ satisfy
\be \fm\;=\;\lim_{\epsilon\to0}\diag((n+\epsilon)_{n\geq0})^{-1}\tm\diag((n+\epsilon)_{n\geq0}).\ee
\end{lemma}
The above lemma shows that the total positivity of $\tm$ implies that of $\fm$, however, by observing that $t_{n,0}=t_{n,1}$ we also obtain the following identity:

\begin{lemma}\label{lem.fe}
The matrices $\fm$ and $\tm$ satisfy
\be \tm\;=\;\lim_{\epsilon\to0}\diag((n+\epsilon)_{n\geq0})\fm\diag((n+\epsilon)_{n\geq 0})^{-1}E_0(1)\ee
where $E_{0}(1)$ is the elementary bidiagonal matrix.
\end{lemma}
\begin{proof}
This follows from Lemma~\ref{lem.transF} and observing that right multiplying a matrix with $E_0(1)$ corresponds to adding column $1$ to column $0$.
\end{proof}
We conclude that total positivity of $\fm$ implies that of $\tm$, but there are yet more identities to uncover. By observing that
\begin{multline} f_{n+1,k+1}\;=\binom{n}{k}(n+1)^{n-k}\;=\;\binom{n}{k}\sum_{l=0}^{n-k}\binom{n-k}{l}n^{n-k-l}\\\;=\;\sum_{l=k}^n\binom{n}{l}\binom{l}{k}n^{n-l}\;=\;\sum_{l=k}^n\binom{l}{k}t_{n,l}\end{multline}
we obtain:
\begin{lemma}\label{lem.ftbi}
The matrices $\fm$ and $\tm$ satisfy
\be \fm \;=\;\left[\begin{array}{c|c}
          1 & 0 \\
          \hline
          0 & \tm B
          \end{array}\right]\ee
          where $B$ is the weighted binomial matrix $B_{x,y}$ with $x=y=1$,
          \be B\ceq \left(\binom{n}{k}\right)_{n,k\geq 0}.\ee
\end{lemma}

Combining Lemmas~\ref{lem.transF}--\ref{lem.ftbi} yields the curious dual identities below:
\begin{corollary}\label{cor.cur}
The matrices $\fm$ and $\tm$ satisfy
\be\label{eq.ffd} \fm\;=\;\left[\begin{array}{c|c}
          1 & 0 \\
          \hline
          0 & \lim_{\epsilon\to0}\diag((n+\epsilon)_{n\geq0})\fm\diag((n+\epsilon)_{n\geq 0})^{-1}E_0(1)B
          \end{array}\right]\ee
and
\be\label{eq.ttd}\tm \;=\;\lim_{\epsilon\to0}\diag((n+\epsilon)_{n\geq0})\left[\begin{array}{c|c}
          1 & 0 \\
          \hline
          0 & \tm B
          \end{array}\right]\diag((n+\epsilon)_{n\geq0})^{-1}E_0(1).\ee
\end{corollary}

But what of the more general matrices $\fm(q)$ and $\tm(q)$? In order to understand how the $q$-forests and $q$-trees matrices are related we first present some well-known identities concerning $q$-binomial coefficients.

The $q$-binomial coefficient satisfies the dual recurrences:
\be\label{eq.binomIdentities1} \qbin{n}{k}\;=\;\qbin{n}{n-k}\;=\;\displaystyle\qbin{n-1}{k-1}+q^k\qbin{n-1}{k}\ee
and
\be\label{eq.binomIdentities2} \qbin{n}{k}\;=\;\qbin{n}{n-k}\;=\;q^{n-k}\displaystyle\qbin{n-1}{k-1}+\qbin{n-1}{k}\ee
(see~\cite[equations~(3.3.3) and~(3.3.4)]{andrews84}, for example). Combining~\eqref{eq.binomIdentities1} and~\eqref{eq.binomIdentities2} above we obtain
\be \qbin{n-1}{k-1}\;=\;\frac{1-q^k}{1-q^n}\qbin{n}{k},\ee
so it follows that for the $q$-forest numbers
\be f_{n,k}(q)\;=\;\qbin{n-1}{k-1}(\qn{n})^{n-k}\ee
and the $q$-trees numbers
\be t_{n,k}(q)\;=\;\left(\qbin{n}{k}(\qn{n})^{n-k}\right)_{n,k\geq 0}\ee
we have
\be f_{n,k}(q)\;=\;\frac{\qn{k}}{\qn{n}}t_{n,k}(q)\ee
for $n\geq 1$.
We therefore have a direct $q$-generalisation of Lemma~\ref{lem.transF}:
\begin{lemma}\label{lem.transFq}
The matrices $\fm(q)$ and $\tm(q)$ satisfy
\be \fm(q)\;=\;\lim_{\epsilon\to0}\diag((\qn{n}+\epsilon)_{n\geq 0})^{-1}\tm(q)\diag((\qn{n}+\epsilon)_{n\geq 0})\ee
\end{lemma}

By observing that $t_{n,0}(q)=t_{n,1}(q)$ we also obtain a $q$-generalisation of Lemma~\ref{lem.fe}:
\begin{lemma}\label{lem.feq}
The matrices $\fm(q)$ and $\tm(q)$ satisfy
\be \tm(q)\;=\;\lim_{\epsilon\to0}\diag((\qn{n}+\epsilon)_{n\geq 0})\fm(q)\diag((\qn{n}+\epsilon)_{n\geq 0})^{-1}E_0(1).\ee
\end{lemma}

Since the entries of $\diag((\qn{n}+\epsilon)_{n\geq 0})^{-1}$ are \emph{rational} functions of $q$, Lemma~\ref{lem.transFq} shows that the \emph{pointwise} total positivity of $\fm(q)$ implies that of $\tm(q)$, and vice versa (also, note that Lemma~\ref{lem.feq} shows that the pointwise total positivity of $\fm(q)$ implies that of $\tm(q)$). It follows that to show both of these matrices are pointwise totally positive (Corollary~\ref{cor.pointwise} below) it suffices to prove that one of them is. Neither of the above lemmas, however, are of any help when considering the \emph{coefficientwise} total positivity of our matrices. 

We have one more simple identity that shows the coefficientwise total positivity of $\fm(q)$ implies the coefficientwise total positivity of $\tm(q)$. It turns out that $\tm(q)$ can be obtained by right-multiplying $\fm(q)$ with a simple inverse lower-bidiagonal matrix with entries that are \emph{polynomials} in $q$.

\begin{lemma}\label{lem.fqtqT}
The matrices $\fm(q)$ and $\tm(q)$ satisfy
\be\label{eq.invFq} \tm(q) \;=\; \fm(q)T(\bfa(q)),\ee
where $T(\bfa(q))$ is the inverse bidiagonal matrix with edge sequence 
\be\bfa(q)\;=\;(q^n\qn{n+1})_{n\geq 0}.\ee
\end{lemma}

\begin{proof}
Since $T(\bfa(q))$ is the inverse of the lower-bidiagonal matrix $L(-\bfa(q))$ (Corollary~\ref{cor.invbidiag} above), proving~\eqref{eq.invFq} is equivalent to showing that
\be \fm(q)\;=\;\tm(q)L(-\bfa(q)),\ee
where $L(-\bfa(q))$ is a lower-bidiagonal matrix with 1s on the diagonal, $(n,n-1)$-entry 
\be-q^n\qn{n+1}\ee
and all other entries $0$. We have
\be(\tm(q)L(-\bfa(q)))_{n,k}\;=\;\qbin{n}{k}\qn{n}^{n-k}-q^k\qn{k+1}\qbin{n}{k+1}\qn{n}^{n-k-1}\ee
which clearly reduces to $\delta_{nk}$ for $k=0$. For $k> 1$ the right-hand side reduces, via the dual recurrences~\eqref{eq.binomIdentities1} and~\eqref{eq.binomIdentities2}, to
\be \qn{n}^{n-k}\left(\qbin{n}{k}-q^{k}\qbin{n-1}{k}\right)\;=\;\qbin{n-1}{k-1}\qn{n}^{n-k}\;=\;f_{n,k}(q),\ee
completing the proof.
\end{proof}
The matrix $T(\bfa(q))$ in the above lemma is clearly coefficientwise totally positive, so it follows from Lemma~\ref{lem.fqtqT} that in order to prove that $\tm(q)$ is coefficientwise totally positive it suffices to show that $\fm(q)$ is coefficientwise totally positive. We were unable to find a similar identity that proves the coefficientwise total positivity of $\tm(q)$ implies that of $\fm(q)$, and this is why we focus on the matrix $\fm(q)$ in the sequel.

In the next section we describe how to specialise the weights of the binomial-like network to rational functions of $q$, and show that the corresponding path matrix agrees with $\fm(q)$. In Section~\ref{sec.coeffpos} we show how this network can be transformed into a different network with weights that are \emph{polynomials} in $q$, thus proving that $\fm(q)$ (and therefore $\tm(q)$) is coefficientwise totally positive in $\Z[q]$.

\section{Pointwise total positivity of the \emph{q}-forests and \emph{q}-trees matrices}\label{sec.pointwise}
\begin{figure}\centering
\includegraphics[scale=0.6]{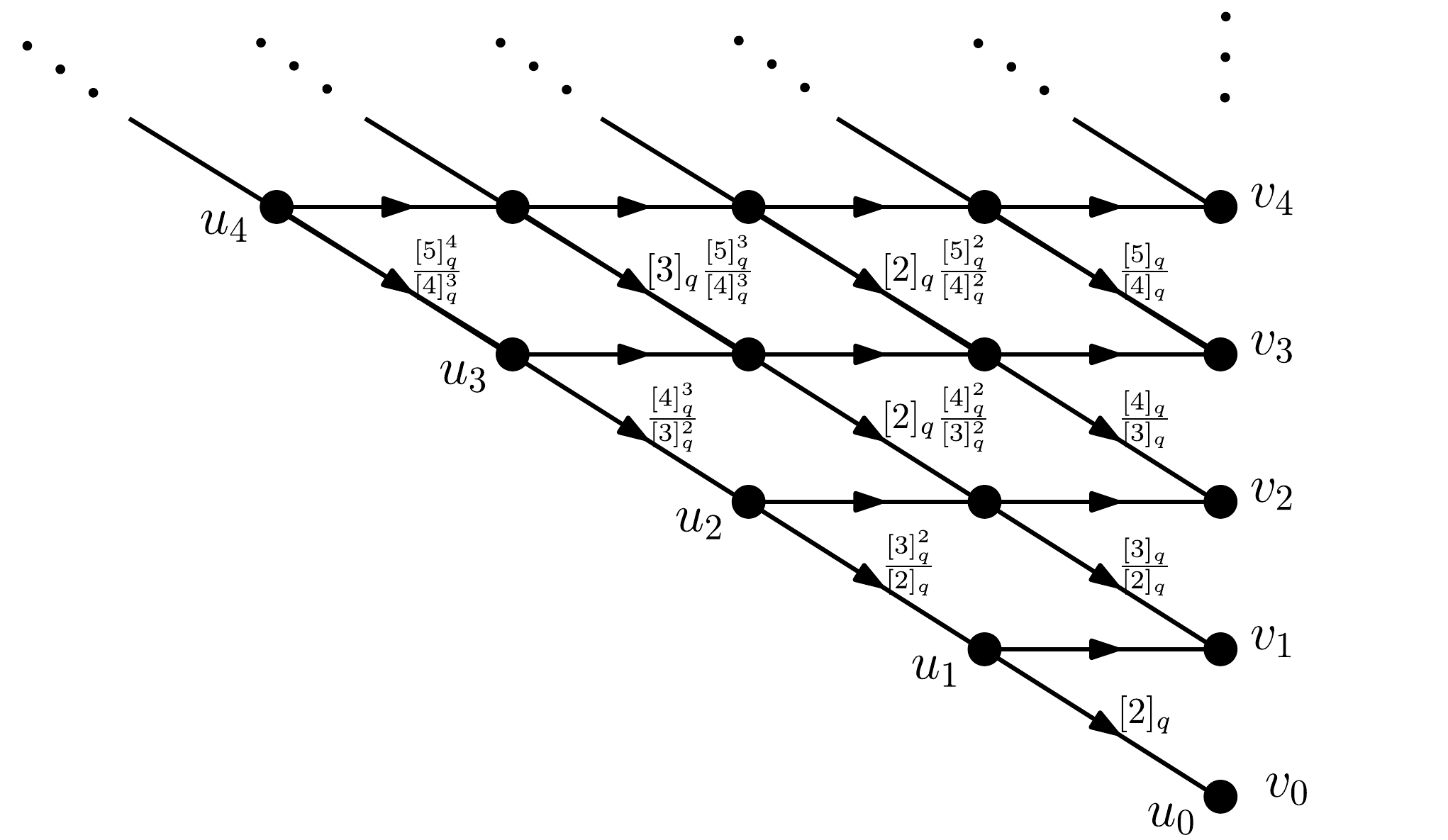}\caption{The binomial-like network $\scrn_1'$ up to $u_4$ and $v_4$.}\label{fig.binomialqrat}\end{figure}
Recall the binomial-like network $\scrn$ with edge weights $\balph$ and $\bbet$ described in Section~\ref{sec.prelim}. This section is devoted to showing how specialising the weights $\balph$ and $\bbet$ to pointwise nonnegative elements of the field $\Q(q)$ (that is, the fraction field of the polynomial ring $\Z[q]$) yields a path matrix that agrees with $\fm(q)$, from which the pointwise total positivity of $\fm(q)$ (and thus, $\tm(q)$) trivially follows. We then present the production-like and quasi-production-like factorisations of $\fm(q)$.

\subsection{The planar network \texorpdfstring{$\scrn_1'$}{TEXT}}
We begin by observing that since
\be \fm(q)\;=\;\left[\begin{array}{c|c}
            1 & 0\\\hline
            0 & \fm'(q)\end{array}\right]\ee
where 
\be \fm'(q)\ceq(f_{n+1,k+1}(q))_{n,k\geq 0}\;=\;\qbin{n}{k}(\qn{n+1})^{n-k},\ee
the pointwise (and indeed, coefficientwise) total positivity of $\fm(q)$ is trivially equivalent to that of $\fm'(q)$.

Now let $\scrn_1'$ denote the binomial-like planar network $\scrn$ in which
\be\label{eq.alphaweights} \alpha_{i,l}\;=\;1\ee
and
\be\label{eq.betaweights} \beta_{i,l}\;=\;q^l\qn{i}\left(\frac{\qn{i+l+1}}{\qn{i+l}}\right)^{i}\ee
(see Figure~\ref{fig.binomialqrat}). Observe that since $\scrn_1'$ is a binomial-like network in which the horizontal edges have weight $1$ it must be a Neville network for some unit-lower-triangular matrix. We make the following proposition:

\begin{proposition}\label{prop.forestQ}
The path matrix matrix $P_{\scrn_1'}$ satisfies
\be P_{\scrn_1'}\;=\;\fm'(q).\ee
\end{proposition}

Please note that by relabelling the source vertices $u_n\to u_{n+1}$ and $v_k\to v_{k+1}$, and inserting a new source vertex $u_0$ at the bottom of $\scrn_1'$ connected via a horizontal edge of weight $1$ to a new sink vertex $v_0$, we obtain the planar network $\scrn_1$ with corresponding path matrix
\be P_{\scrn_1}\;=\;\left[\begin{array}{c|c}
            1 & 0\\\hline
            0 & P_{\scrn_1'}\end{array}\right].\ee
Since the weights of $\scrn_1$ are nonnegative for $q\in\Rp$, the pointwise total positivity of $\fm(q)$ (and therefore also $\tm(q)$ according to Lemma~\ref{lem.fqtqT}) follows from Proposition~\ref{prop.forestQ}. The author would like to point out that he is very grateful to Xi Chen and Shaoshi Chen for discovering a proof of Proposition~\ref{prop.forestQ} with $q$ specialised to $q=1$, which inspired what follows.

Our proof of Proposition~\ref{prop.forestQ} arises from considering directly the sum over weighted paths in the network $\scrn_1'$ from $u_n$ to $v_k$:
\begin{multline}\label{eq.Fpaths}P_{\scrn_1'}(u_{n}\to v_{k})\;=\;\sum_{1\leq i_1<\cdots<i_{n-k}\leq n}\qn{i_1}\left(\frac{\qn{k+2}}{\qn{k+1}}\right)^{i_1}q^{k-i_1+1}\\\cdots\qn{i_{n-k}}\left(\frac{\qn{n+1}}{\qn{n}}\right)^{i_{n-k}}q^{n-i_{n-k}}\end{multline}
(see~\eqref{eq.symm}, where the $\beta_{i,l}$s have been replaced with the weights defined in~\eqref{eq.betaweights}). Proving Proposition~\ref{prop.forestQ} amounts to showing that the right-hand side of~\eqref{eq.Fpaths} equals
\be \qbin{n-1}{k-1}(\qn{n})^{n-k}.\ee
Since the nested sum in~\eqref{eq.Fpaths} is somewhat complicated, in order to make better sense of it we first prove the following lemma:
\begin{lemma}\label{lem.SnkF}
Given integers $n,m,t$ with $0\leq m < n$ and  $0\leq t < n-m$ we have
\be\label{eq.SnkF.ind1} \sum_{s=t+1}^{n-m}q^{n-m-s}\qn{s}\left(\frac{\qn{n+1}}{\qn{n-m}}\right)^{s}\qbin{n-s}{m} \;=\; \qn{n+1}\qbin{n-t}{m+1}\qpow{n+1}{n-m}{t}.\ee
\end{lemma}

\begin{proof}[Proof of Lemma~\ref{lem.SnkF}]It is easy to verify~\eqref{eq.SnkF.ind1} for $t=n-m-1$ so suppose it holds for $t\leq n-m-1$ and consider the case $t-1$. We have
\begin{multline} \sum_{s=t}^{n-m}q^{n-m}\qn{s}\left(\frac{\qn{n+1}}{q\qn{n-m}}\right)^{s}\qbin{n-s}{m} \;=\;\qn{t}q^{n-m}\left(\frac{\qn{n+1}}{q\qn{n-m}}\right)^{t}\qbin{n-t}{m}\\
+ \qn{n+1}\qbin{n-t}{m+1}\qpow{n+1}{n-m}{t},\end{multline}
the right-hand side of which can be written
\begin{multline}\label{eq.qbinCancel} \left(\frac{\qn{n+1}}{\qn{n-m}}\right)^{t}\left(\qn{t}\left(q^{n-m-t}\qbin{n-t}{m}+\qbin{n-t}{m+1}\right)\right.\\\left.\;+\;q^t\qn{n-t-m}\frac{\qn{n+1-t}}{\qn{n-t-m}}\qbin{n-t}{m+1}\right).\end{multline}
The dual recurrences for $q$-binomial coefficients~\eqref{eq.binomIdentities1} and~\eqref{eq.binomIdentities2} then reduce~\eqref{eq.qbinCancel} to 
\be\label{eq.qbinCancel2} \left(\frac{\qn{n+1}}{\qn{n-m}}\right)^{t}\qbin{n-t+1}{m+1}\qn{n-m},\ee
which agrees with the right-hand side of~\eqref{eq.SnkF.ind1} when $t$ is replaced with $t-1$.\end{proof}

With Lemma~\ref{lem.SnkF} in hand we are now ready to prove Proposition~\ref{prop.forestQ}.
\begin{proof}[Proof of Proposition~\ref{prop.forestQ}]
We will show that
\be P_{\scrn_1'}(u_n\to v_k)\;=\;\qbin{n}{k}(\qn{n+1})^{n-k},\ee
where $P_{\scrn_1'}(u_{n}\to v_{k})$ denotes the sum over weighted directed paths in $\scrn_1'$ starting at $u_n$ and ending at $v_k$. According to~\eqref{eq.Fpaths} we have
\begin{multline}\label{eq.PNFpaths}P_{\scrn_1'}(u_{n}\to v_{k})\;=\;\sum_{1\leq i_1<\cdots<i_{n-k}\leq n}\qn{i_1}\left(\frac{\qn{k+2}}{\qn{k+1}}\right)^{i_1}q^{k-i_1+1}\\\cdots\qn{i_{n-k}}\left(\frac{\qn{n+1}}{\qn{n}}\right)^{i_{n-k}}q^{n-i_{n-k}}.\end{multline}
The nested sum~\eqref{eq.PNFpaths} can be expressed recursively. Let
\be\label{eq.snkdef} S_{n-k}\ceq\sum_{i_{n-k}=i_{n-k-1}+1}^{n}\qn{i_{n-k}}\qpow{n+1}{n}{i_{n-k}}q^{n-i_{n-k}},\ee
and for $1\leq m\leq n-k-1$ define
\be S_{n-k-m}\ceq \sum_{i_{n-k-m}=i_{n-k-m-1}+1}^{n-m}\qn{i_{n-k-m}}\left(\frac{\qn{n+1-m}}{q\qn{n-m}}\right)^{i_{n-k-m}}q^{n-m}S_{n-k-m+1},\ee
so that in particular
\be\label{eq.eval} P_{\scrn_1'}(u_{n}\to v_{k}) \;=\;S_{n-k-m}|_{m=n-k-1,i_0=0}\ee
(where $S_{n-k-m}|_{m=n-k-1,i_0=0}$ denotes the right-hand side of $S_{n-k-m}$ where $m$ is replaced with $n-k-1$ and $i_0=0$). We make the following claim: 

{\bf Claim.} For $0\leq m\leq n-k-1$ we have
\be\label{eq.claim} S_{n-k-m} \;=\; ([n+1]_q)^{m+1}\qbin{n-i_{n-k-m-1}}{m+1}\left(\frac{[n+1]_q}{[n-m]_q}\right)^{i_{n-k-m-1}}.\ee
To prove~\eqref{eq.claim} we first observe that since
\be S_{n-k}\;=\;\sum_{i_{n-k}=i_{n-k-1}+1}^{n}\qn{i_{n-k}}\qpow{n+1}{n}{i_{n-k}}q^{n-i_{n-k}}\qbin{n-i_{n-k}}{0}\ee
(remember $\qbin{n}{0}=1$), the claim clearly holds for $m=0$ by virtue of Lemma~\ref{lem.SnkF}. Suppose then~\eqref{eq.claim} holds for $m\geq 0$ and consider the case $m_1 = m+1$. Then
\begin{eqnarray} S_{n-k-m_1}&=&\sum_{s=t+1}^{n-m_1}\qn{s}\left(\frac{\qn{n+1-m_1}}{q\qn{n-m_1}}\right)^{s}q^{n-m_1}S_{n-k-l}\\\label{eq.induction2}
              &=&\qn{n+1}^{m_1}\sum_{s=t+1}^{n-m_1-s}q^{n-m_1}\qn{s}\left(\frac{\qn{n+1}}{\qn{n-m_1}}\right)^{s}\qbin{n-s}{m_1}\end{eqnarray}
by the induction hypothesis (where we have replaced $i_{n-k-m-1}$ and $i_{n-k-m-2}$ with $s$ and $t$ respectively). Applying Lemma~\ref{lem.SnkF} to~\eqref{eq.induction2} then proves the claim~\eqref{eq.claim}.

Substituting $m=n-k-1$ and $i_0=0$ in the right-hand side of~\eqref{eq.claim}) for $n\geq k\geq 0$ yields
\begin{eqnarray} P_{\scrn_1'}(u_{n}\to v_{k})&=&([n+1]_q)^{n-k}\qbin{n}{n-k}\\
                &=&\qbin{n}{k}(\qn{n+1})^{n-k},\end{eqnarray}
completing the proof.

\end{proof}

The above proposition immediately implies the following corollary:
\begin{corollary}\label{cor.pointwise}
The matrices $\fm(q)$ and $\tm(q)$ are pointwise totally positive for $q$ specialised to $q\in\Rp$.
\end{corollary}

\begin{proof}
According to Proposition~\ref{prop.forestQ} the matrix $\fm'(q)$ is the path matrix corresponding to a planar network with weights that are rational functions of $q$. Each weight in $\scrn_1'$ is pointwise totally positive for $q\in\Rp$, so according to the LGV lemma (see Subsection~\ref{subsec.LGV}) $\fm'(q)$ is pointwise totally positive for $q\in\Rp$. Since pointwise total positivity of $\fm(q)$ is trivially equivalent to that of $\fm'(q)$, the matrix $\fm(q)$ is pointwise totally positive. We have already seen in Section~\ref{sec.combin.interp} (see Lemma~\ref{lem.fqtqT}, for example) that pointwise total positivity of $\fm(q)$ implies that of $\tm(q)$, so $\tm(q)$ and $\fm(q)$ are pointwise totally positive.
\end{proof}

Of course for $q$ specialised to $q=1$, Corollary~\ref{cor.pointwise} proves that both the trees matrix $\tm$ and the forests matrix $\fm$ are totally positive (Theorem~\ref{thm.tree}). Moreover, the planar network in Figure~\ref{fig.binomialqrat} implies, via Proposition~\ref{prop.forestQ}, a matrix factorisation of $\fm'(q)$ and hence a factorisation of $\fm(q)$. Indeed, we can obtain two such factorisations depending on whether we view $\scrn_1$ as a concatenation of column transfer matrices (yielding a production-like factorisation, see Corollary~\ref{cor.prodlikeFactor} above), or a concatenation of diagonal column transfer matrices (yielding a quasi-production-like factorisation, see Corollary~\ref{cor.quasiprodlike} above).

\begin{corollary}\label{cor.FqFactor}
The matrix $\fm(q)$ has the production-like factorisation
\be \fm(q)\;=\;\cdots\left[\begin{array}{c|c}
            I_3 & 0\\\hline
            0 & L(\bfa_2(q))\end{array}\right]\cdot\left[\begin{array}{c|c}
            I_2 & 0\\\hline
            0 & L(\bfa_1(q))\end{array}\right]\cdot\left[\begin{array}{c|c}
            I_1 & 0\\\hline
            0 &L(\bfa_0(q))\end{array}\right]\ee
where $L(\bfa_i(q))$ is the lower-bidiagonal matrix with corresponding edge sequence
\be \bfa_i(q)\;=\;\left(q^n\qn{i+1}\left(\frac{\qn{n+i+2}}{\qn{n+i+1}}\right)^{i+1}\right)_{n\geq 0}.\ee
\end{corollary}

We will use the above matrix factorisation in our proof of Theorem~\ref{thm.fqtqCoeff}, however, we also remark that $\fm(q)$ also has a quasi-production-like factorisation (thanks to Corollary~\ref{cor.quasiprodlike}):
\begin{corollary}\label{cor.FqFactor2}
The matrix $\fm(q)$ has the quasi-production-like factorisation 
\be \fm(q)\;=\;\left[\begin{array}{c|c}
            1 & 0\\\hline
            0 & T(\bfa^*_0(q))\end{array}\right]\cdot\left[\begin{array}{c|c}
            I_{2} & 0\\\hline
            0 &T(\bfa^*_{1}(q))\end{array}\right]\cdot\left[\begin{array}{c|c}
            I_3 & 0\\\hline
            0 & T(\bfa^*_2(q))\end{array}\right]\cdots\ee
where $T(\bfa^*_l(q))$ is the inverse bidiagonal matrix with corresponding edge sequence
\be \bfa^*_l(q)\ceq\left(q^{l+1}\qn{n+1}\left(\frac{\qn{n+l+3}}{\qn{n+l+2}}\right)^{n+1}\right)_{n\geq 0}.\ee
\end{corollary}

Having established the pointwise total positivity of the $q$-forests and $q$-trees matrix via the planar network $\scrn_1'$, in the next section we will show how $\scrn_1'$ can be transformed into a different planar network with weights that are polynomial in $q$, from which the coefficientwise total positivity of $\fm(q)$ and $\tm(q)$ follows.

\section{Coefficientwise total positivity of the \texorpdfstring{$q$}{TEXT}-forests and \texorpdfstring{$q$}{TEXT}-trees matrices}\label{sec.coeffpos}

We now present a different planar network with weights that are polynomials in $q$ and show that its corresponding path matrix can be transformed into the production-like factorisation of $\fm(q)$ from the previous section.

Let $\scrn_2$ denote the planar network in Figure~\ref{fig.coeffFq}. Observe that this network does \emph{not} have the structure of the binomial-like network that we have considered up to this point, nonetheless we can still apply the tools from Section~\ref{sec.prelim} to understand the corresponding path matrix. In particular we will repeatedly make use of Lemma~\ref{lem.switch}, so for ease of reading we recall it here:

\begin{lemma*}[Lemma~\ref{lem.switch} above]
Suppose $\bfa=(a_n)_{n\geq 0}$ and $\bfb=(b_n)_{n\geq 0}$ are sequences of elements belonging to a field $F$. Then
\begin{enumerate}[(i)]
\item If $a_n=-b_n$ for all $n$ then $L(\bfa)T(\bfb)=I$;
\item If $a_n+b_n\neq 0$ for all $n$ then
\be L({\bfa})T({\bfb}) = \left[\begin{array}{c|c}
          1 & 0 \\
          \hline
          0 & T({\bfb'})
          \end{array}\right]\cdot L({\bfa'})\ee
where $\bfa' = (a'_n)_{n\geq 0}$ is the edge sequence in which 
\be a'_n\ceq \begin{cases}a_0+b_0 & \textrm{if }n=0,\\ \frac{a_{n-1}(a_n+b_n)}{a_{n-1}+b_{n-1}} & \textrm{if }n>0,\end{cases}\ee
and $\bfb'=(b'_n)_{n\geq 0}$ is the edge sequence in which
\be b'_n = \frac{b_{n}(a_{n+1}+b_{n+1})}{a_n+b_n}.\ee
\end{enumerate}
\end{lemma*}

Consider a subnetwork $C^i$ of $\scrn_2$ depicted in Figure~\ref{fig.coeffFq}. The column transfer matrix corresponding to each subnetwork is
\be\label{eq.Pcicoeff} P_{C^i}\;=\;\left[\begin{array}{c|c}
            I_{i+1} & 0\\\hline
            0 & L(\bfa_i(q))T(\bfb_i(q))\end{array}\right]\ee
where $L(\bfa_i(q))$ and $T(\bfb_i(q))$ are lower-bidiagonal and inverse lower-bidiagonal matrices (respectively) with corresponding edge sequences
\be\bfa_i(q)\;=\;(q^{i+n+1}[i+1]_{q})_{n\geq 0}\ee
and
\be\bfb_i(q)\;=\;(q^{i}[n+1]_{q})_{n\geq0}.\ee
The path matrix corresponding to $\scrn_2$ has the factorisation
\be\label{eq.n2fac} P_{\scrn_2}\;=\;\cdots P_{C^2}P_{C^1}P_{C^0}\ee
We make the following proposition:
\begin{figure}
\includegraphics[scale=0.4]{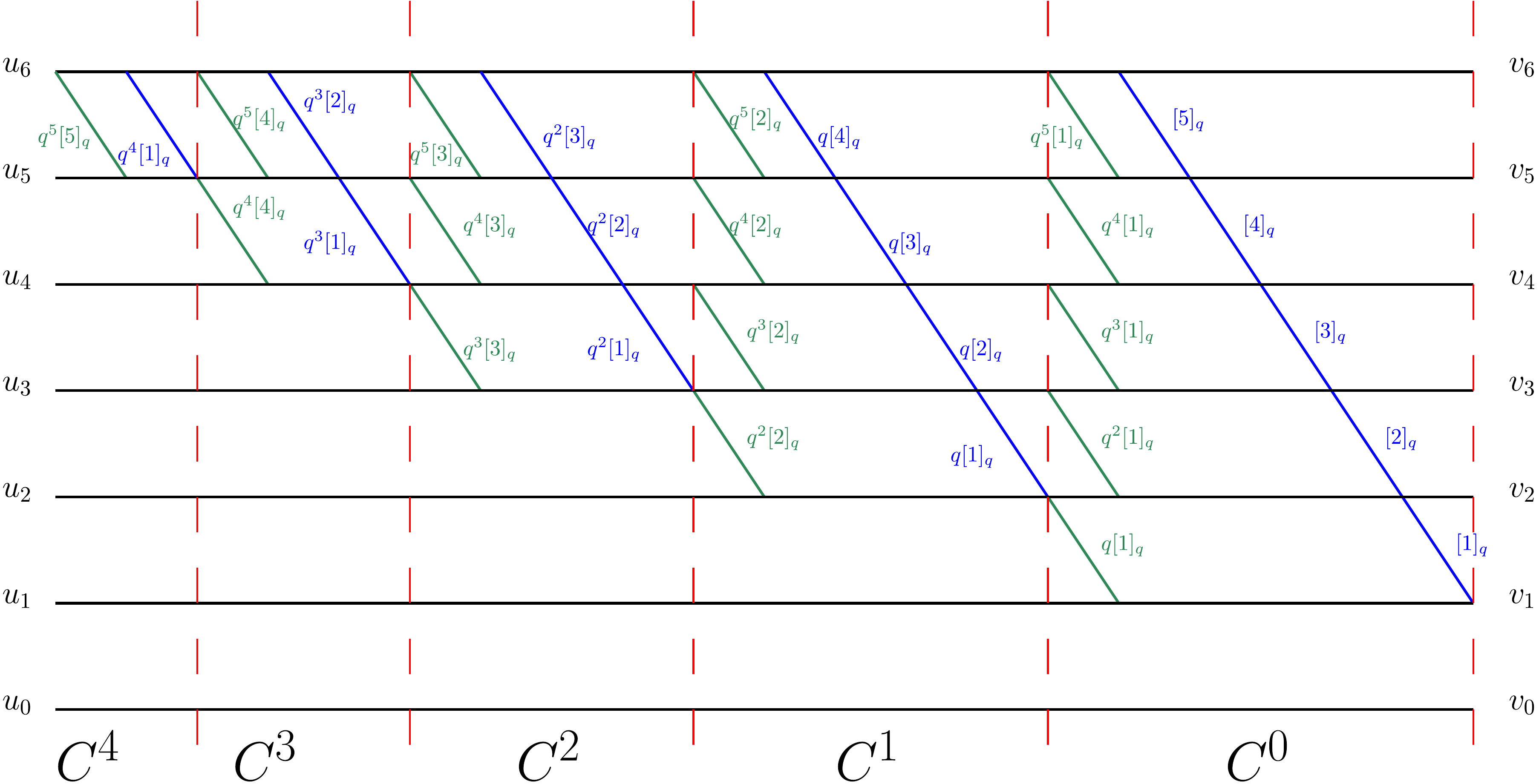}\caption{The planar network $\scrn_2$ up to source $u_6$ and sink $v_6$. Note that the dashed red lines do not form part of the network, they simply delineate the column transfer matrices of $\scrn_2$. Also note that the direction of the edges have been omitted, the reader should assume all edges are directed from left to right.}\label{fig.treeqDigraph}\label{fig.coeffFq}
\end{figure}

\begin{proposition}\label{prop.forestcoeff}
The path matrix for $\scrn_2$ satisfies
\be\label{eq.PNFFQ} P_{\scrn_2}\;=\;\fm(q).\ee
\end{proposition}

Before beginning the proof proper we first outline our approach. The factorisation of $P_{\scrn_2}$ in~\eqref{eq.n2fac} is of the form
\be\label{eq.n2fac2}\cdots (L_4T_4)(L_3T_3)(L_2T_2)(L_1T_1),\ee
where $L_i$ denotes a lower-bidiagonal matrix $L$ \emph{downshifted by $i$}:
\be L_i\;=\left[\begin{array}{c|c}
            I_{i} & 0\\\hline
            0 & L\end{array}\right],\ee
and similarly $T_i$ denotes an inverse lower-bidiagonal matrix $T$ downshifted by $i$:
\be T_i\;=\;\left[\begin{array}{c|c}
            I_{i} & 0\\\hline
            0 & T\end{array}\right].\ee
Applying Lemma~\ref{lem.switch} to the rightmost product $L_1T_1$ we obtain
\be L_1T_1\;=\;T_2'L_1',\ee
where $T_2'$ is an inverse lower-bidiagonal matrix downshifted by $2$ and $L_1'$ a lower-bidiagonal matrix downshifted by $1$. The factorisation in~\eqref{eq.n2fac2} can therefore be expressed as
\be\cdots (L_4T_4)(L_3T_3)(L_2T_2T_2')L_1'.\ee
More generally, Lemma~\ref{lem.switch} enables us to ``pull'' lower-bidiagonal matrices ``through'' inverse lower-bidiagonal matrices to the right (compare Figure~\ref{fig.switch} with Figure~\ref{fig.coeffFq}), so that in particular for a given $i>0$ the factorisation in~\ref{eq.n2fac2} can be written in the form
\be \cdots (L_{i+2}T_{i+2})\left(L_{i+1}T_{i+1}\prod_{j=1}^iT_{i+1,j}\right)L_i'L_{i-1}'\cdots L_1'\ee
where each $T_{i+1,j}$ is an inverse bidiagonal matrix downshifted by $i+1$, and each $L_i'$ is a lower-bidiagonal matrix downshifted by $i$. It is clear that in the limit (as $i\to\infty$) the inverse lower-bidiagonal matrices disappear, and what remains is a factorisation of $P_{\scrn_2}$ into downshifted lower-bidiagonal matrices:
\be P_{\scrn_2}\;=\;\cdots L_4'L_3'L_2'L_1'.\ee
We will show that this resulting factorisation of $P_{\scrn_2}$ agrees with the production-like factorisation of $\fm(q)$ given in Corollary~\ref{cor.FqFactor} of the previous section.

\begin{proof}[Proof of Proposition~\ref{prop.forestcoeff}]
It is easy to see that the matrices agree in the first column, since 
\be P_{\scrn_2}(u_n\to v_0)\;=\;\delta_{n0}\;=\;f_{n,0}(q),\ee so consider again the matrix
\be \fm'(q)\ceq(f'_{n,k}(q))_{n,k\geq 0}\;=\;\left(\qbin{n}{k}([n+1]_q)^{n-k}\right)_{n,k\geq 0}\ee
and let $\scrn_2'$ denote the planar network obtained by removing from $\scrn_2$ vertices $u_0,v_0$ and relabelling vertices $u_i\to u_{i-1}$, $v_i\to v_{i-1}$ for $i>0$ (the corresponding path matrix $P_{\scrn_2'}$ is thus obtained by deleting the first row and column of $P_{\scrn_2}$).

As described above, we will prove~\eqref{eq.PNFFQ} by repeatedly applying Lemma~\ref{lem.switch} to the matrix factorisation of $P_{\scrn_2'}$, thereby expressing $P_{\scrn_2}$ as a product of lower-bidiagonal matrices that agree with those given in Corollary~\ref{cor.FqFactor}. To this end, fix $i\geq 0$ and consider the product
\be P'_{C^i}\cdots P'_{C^1}P'_{C^0}\; =\; \left[\begin{array}{c|c}
            I_i & 0\\\hline
            0 & L(\bfa_i(q))T(\bfb_i(q))\end{array}\right]\cdots\left[\begin{array}{c|c}
            1 & 0\\\hline
            0 & L(\bfa_1(q))T(\bfb_1(q))\end{array}\right]\cdot L(\bfa_0(q))T(\bfb_0(q))\ee
where
\be\bfa_i(q)\;=\;(q^{i+n+1}\qn{i+1})_{n\geq 0}\ee
and
\be\bfb_i(q)\;=\;(q^{i}\qn{n+1})_{n\geq0}\ee
(each matrix $P'_{C^i}$ is obtained by deleting the first row and column from $P_{C^i}$ defined in~\eqref{eq.Pcicoeff} above).

{\bf Claim 1.} We claim that
\be\label{eq.claimPCi} P'_{C^i}\cdots P'_{C^1}P'_{C^0}\;=\;\left[\begin{array}{c|c}
          I_{i+1} & 0 \\
          \hline
          0 & {\bf T}^i(q)
          \end{array}\right]\cdot \left[\begin{array}{c|c}
            I_i & 0\\\hline
            0 & L(\bfc_i(q))\end{array}\right]\cdots\left[\begin{array}{c|c}
            1 & 0\\\hline
            0 & L(\bfc_1(q))\end{array}\right]L(\bfc_0(q))\ee
where 
\be {\bf T}^i(q) = \prod_{j=0}^{i}T(\bfb_{i,j}(q))\ee
is a product of inverse bidiagonal matrices with corresponding edge sequences 
\be \bfb_{i,j}(q)\ceq\left(q^{i-j}\qn{n+1}\left(\frac{\qn{n+i+3}}{\qn{n+i+2}}\right)^{j+1}\right)_{n\geq 0}\ee
and each $L(\bfc_i(q))$ is the lower-bidiagonal matrix with edge sequence
\be \bfc_i(q)\;=\;\left(q^n\qn{i+1}\left(\frac{\qn{n+i+2}}{\qn{n+i+1}}\right)^{i+1}\right)_{n\geq 0}.\ee
Please note that the lower-bidiagonal matrices $L(\bfc_i(q))$ are obtained by deleting the first row and column of the matrices found in the production-like factorisation of $\fm(q)$ in Corollary~\ref{cor.FqFactor} above, so proof of the proposition follows once we have proved the claim (since~\eqref{eq.claimPCi} implies that for $0\leq n\leq i+1$, each $(n,k)$-entry of $P'_{C^i}\cdots P'_{C^1}P'_{C^0}$ agrees with that of $\fm'(q)$).

First observe that we can apply Lemma~\ref{lem.switch} to the product $L(\bfa_i(q))T(\bfb_i(q))$ in $P'_{C^i}$, yielding
\be\label{eq.Cfac1} P'_{C^i}\; =\;\left[\begin{array}{c|c}
          I_{i+1} & 0 \\
          \hline
          0 & T(\bfb'_i(q))
          \end{array}\right]\cdot \left[\begin{array}{c|c}
          I_{i} & 0 \\
          \hline
          0 & L(\bfa'_i(q))
          \end{array}\right]\ee
where
\be \bfa'_i(q)\ceq \left(\frac{q^{i+n}\qn{i+1}\qn{n+i+2}}{\qn{n+i+1}}\right)_{n\geq 0},\ee
and
\be \bfb'_i(q) \ceq \left(\frac{q^{i}\qn{n+1}\qn{n+i+3}}{\qn{n+i+2}}\right)_{n\geq 0}.\ee
Clearly $T(\bfb'_i(q))=T(\bfb_{i,0}(q))$, moreover it is easily verified that $\bfa'_0(q)=\bfc_0(q)$, so suppose~\eqref{eq.claimPCi} holds for $i- 1\geq 0$. By the induction hypothesis and~\eqref{eq.Cfac1} we have
\be P'_{C^i}\cdots P'_{C^0} = \left[\begin{array}{c|c}
          I_{i+1} & 0 \\
          \hline
          0 & T(\bfb_{i,0}(q))
          \end{array}\right]\cdot\left[\begin{array}{c|c}
          I_{i} & 0 \\
          \hline
          0 & L(\bfa'_i(q)){\bf T}^{(i-1)}(q)
          \end{array}\right]\cdot \left[\begin{array}{c|c}
            I_{i-1} & 0\\\hline
            0 & L(\bfc_{i-1}(q))\end{array}\right]\cdots L(\bfc_0(q)).\ee
We now make the following claim:

{\bf Claim 2.} For $0\leq m\leq i-1$,

\be\label{eq.claimGamma} L(\bfa'_i(q))\prod_{j=0}^m T(\bfb_{i-1,j}(q))\;=\; \left[\begin{array}{c|c}
          1 & 0 \\
          \hline
          0 & T(\bfb_{i,1}(q))\cdots T(\bfb_{i,m+1}(q))
          \end{array}\right]\cdot L(\bfc'_{m+1}(q))\ee
where $\bfc'_{m}(q)$ is the edge sequence 
\be \bfc'_{m}(q)\ceq \left(q^{i+n-m}\qn{i+1}\left(\frac{\qn{n+i+2}}{\qn{n+i+1}}\right)^{m+1}\right)_{n\geq 0}.\ee
Claim 2 ~\eqref{eq.claimGamma} clearly holds for $m=0$ by way of Lemma~\ref{lem.switch}:
\be L(\bfa'_i(q))T(\bfb_{i-1,0}(q))=\left[\begin{array}{c|c}
          1 & 0 \\
          \hline
          0 & T(\bfb_{i,1}(q))
          \end{array}\right]\cdot L(\bfc'_{1}(q)),\ee
and more generally we have (again by Lemma~\ref{lem.switch})
\be L(\bfc'_m(q))T(\bfb_{i-1,m}(q)) =\left[\begin{array}{c|c}
          1 & 0 \\
          \hline
          0 & T(\bfb_{i,m+1}(q))
          \end{array}\right]\cdot L(\bfc'_{m+1}(q)).\ee
Letting $m=i-1$ above proves Claim 2~\eqref{eq.claimGamma}, which in turn confirms Claim 1~\eqref{eq.claimPCi} since it is straightforward to check that $\bfc'_{i}(q)=\bfc_i(q)$, thereby completing the proof of the proposition.

\end{proof}
Coefficientwise total positivity of $\fm(q)$ and $\tm(q)$ follows almost immediately from Proposition~\ref{prop.forestcoeff}:
\begin{proof}[Proof of Theorem~\ref{thm.fqtqCoeff}]
Since $\scrn_2$ is a planar network with weights that are polynomials with nonnegative coefficients in $q$ and we have just shown that $P_{\scrn_2}=\fm(q)$, the matrix $\fm(q)$ is coefficientwise totally positive in $\Z[q]$ thanks to the LGV lemma. The coefficientwise total positivity of $\tm(q)$ then follows from the fact that
\be \tm(q)\;=\;\fm(q)T((q^n[n+1]_q)_{n\geq 0})\ee
(see Lemma~\ref{lem.fqtqT} above), since the inverse lower-bidiagonal matrix $T((q^n[n+1]_q)_{n\geq 0})$ is coefficientwise totally positive.
\end{proof}

Proving Theorem~\ref{thm.fqtqCoeff} was the main goal of this paper, however, the planar network $\scrn_2$ can be seen as a specialisation of a more general network that appears to satisfy some interesting properties that we believe warrant further investigation. We discuss this more general network in the following section.

\section{Further comments and some open problems}\label{sec.furthercomments}

In this final section we turn our attention to a generalisation of the planar network $\scrn_2$ that under certain specialisations yield matrices with entries that appear to count a variety of combinatorial objects with respect to different statistics. These matrices are automatically coefficientwise totally positive, since they arise from planar networks with weights that are polynomials with nonnegative coefficients. Furthermore, under certain specialisations, sequences of row-generating polynomials related to these matrices turn out to be sequences of polynomials that are already well-known, and in some cases are also known to be \emph{coefficientwise Hankel-totally positive} (a sequence $\bfa=(P_n(\bfx))_{n\geq 0}$ of polynomials in one or more indeterminates $\bfx$ is \emph{coefficientwise Hankel-totally positive} if its associated \emph{Hankel matrix} $H_{\infty}(\bfa)\ceq (P_{n+k}(\bfx))_{n,k\geq 0}$ is coefficientwise totally positive in all the indeterminates $\bfx$).

We will briefly touch upon the method of \emph{production matrices} (see~\cite{Deutsch09,Deutsch04}), which in recent years has become an important tool in enumerative combinatorics and has its roots in Stieltjes' work on continued fractions\footnote{See~\cite{Stieltjes94,Stieltjes95}, which were reprinted together with an English translation in~\cite[pp.~401--566 and 609--745]{Stieltjes93}}. The theory of production matrices with respect to total positivity is extensively studied in~\cite{Sokal22}, however, since~\cite{Sokal22} is not yet publicly available we direct the reader to Sections~2.2 and~2.3 of~\cite{Sokal_21f} for a fuller treatment. 

For our purposes we will require the following principles: let $\mathsf{P}=(p_{i,j})_{i,j\geq 0}$ be an infinite matrix with entries in $\R[\bfx]$ equipped with the coefficientwise order. In order that powers of $\mathsf{P}$ are well defined we assume that $\mathsf{P}$ is either row finite (that is, has only finitely many nonzero entries in each row) or column finite. Let us now define the infinite matrix $A\ceq(a_{n,k})_{n,k\geq 0}$ where
\be a_{n,k}=(\mathsf{P}^n)_{0,k},\ee
that is, row $n$ of $A$ is the first row of the matrix power $\mathsf{P}^n$ (in particular we set $a_{0,k}=\delta_{0k}$). We call $\mathsf{P}$ the \emph{production matrix} of $A$, while $A$ is referred to as the \emph{output matrix} of $\mathsf{P}$ and we write $\scro(\mathsf{P})=A$. The entries of $A$ are, explicitly,
\be a_{n,k}\;=\;\sum_{i_1,\ldots,i_{n-1}}p_{0,i_1}p_{i_1,i_2}p_{i_2,i_3}\cdots p_{i_{n-2},i_{n-1}}p_{i_{n-1},k},\ee
which can be seen as the total weight of all $n$-step walks in $\N$ from $i_0=0$ to $i_n=k$, in which the weight of a walk is the product of the weights of its steps, and each step from $i$ to $j$ has weight $p_{i,j}$. An equivalent formulation is to define the entries by the recurrence
\be a_{n,k}\;=\;\sum_{i=0}^{\infty}a_{n-1,i}p_{i,k}\ee
for $n\geq 1$ with initial condition $a_{0,k}\;=\;\delta_{0k}$.

Given a production matrix $\mathsf{P}$, we define the \emph{augmented production matrix} to be
\be \tilde{\mathsf{P}}\;\ceq\;\left[\begin{array}{c}1\,\,0\,\,0\,\,0\,\,\cdots\\\hline \mathsf{P}\end{array}\right],\ee
and it is easy to show~\cite[Section~2.2]{Sokal_21f} that the output matrix of $\mathsf{P}$ has the factorisation
\be\label{eq.augProd}\scro(\mathsf{P})\;=\;\cdots\left[\begin{array}{c|c}
            I_3 & 0\\\hline
            0 & \tilde{\mathsf{P}}\end{array}\right]\cdot\left[\begin{array}{c|c}
            I_2 & 0\\\hline
            0 & \tilde{\mathsf{P}}\end{array}\right]\cdot\left[\begin{array}{c|c}
            1 & 0\\\hline
            0 & \tilde{\mathsf{P}}\end{array}\right]\tilde{P}.\ee
Clearly if $\mathsf{P}$ is coefficientwise totally positive then so is $\tilde{\mathsf{P}}$, and hence so is $\scro(\mathsf{P})$ (Theorem~2.9 in~\cite{Sokal_21f}). The factorisation given above is the reason we call the factorisation in Corollary~\ref{cor.prodlikeFactor} \emph{production-like}. The \emph{quasi-production-like} factorisation in Corollary~\ref{cor.quasiprodlike} is named after the related concept of \emph{quasi-production matrices} (see~\cite{Gilmore22}) that we do not require here.

\subsection{General Abel Polynomials}

Let us now return to the $q$-forests matrix $\fm(q)$
and consider the matrix
\be\label{eq.barf}\bar{\fm}(q)\ceq(\bar{f}_{n,k}(q))_{n,k\geq 0}\;=\; \fm(q)\diag((q^{n(n-1)/2})_{n\geq 0}),\ee
the row-generating polynomials of which are defined to be
\be\bar{\fm}_0(x;q)\ceq 1\ee
and
\be \bar{\fm}_n(x;q)\ceq\sum_{k=0}^n\bar{f}_{n,k}(q)x^k\;=\;\sum_{k=0}^n\qbin{n-1}{k-1}(\qn{n})^{n-k}x^k\ee
for $n>0$. By employing the finite sum version of the $q$-binomial theorem~\cite[Theorem~3.3]{andrews84}:
\be \sum_{k=0}^{n}\qbin{n}{k}(-z)^kq^{\binom{k+1}{2}}\;=\;\prod_{i=0}^{n-1}(1-zq^i),\ee
in which $z$ has been replaced with $(-qx/\qn{n})$, we obtain the explicit formulas
\be \bar{\fm}_0(x;q)\;=\;1\ee
and
\be \bar{\fm}_n(x;q)\;=\;x\prod_{i=1}^{n-1}(xq^i+\qn{n}).\ee

For $q$ specialised to $q=1$ the polynomials $\bar{\fm}(x;1)$ reduce to the Abel polynomials
\be \bar{\fm}_n(x;1)\;=\;x(x+n)^{n-1},\ee
which are the row-generating polynomials of the forests matrix $\fm$. Sokal~\cite{Sokal_21f} has already proven that the Abel polynomials (along with further generalisations of the row-generating polynomials of the forests matrix) are coefficientwise Hankel-totally positive. It is natural to ask, then, whether this property is preserved if $q$ is left as an indeterminate in $\bar{\fm}_n(x;q)$. 

The Hankel matrix associated to $\bar{\fm}_n(x;q)$,
\be H_{\infty}((\bar{\fm}_n(x;q))_{n\geq 0})=(\pi_{n,k})_{n,k\geq0},\ee begins:
\be\left[
\begin{array}{ccc}
 x & x (q x+q+1)& \cdots\\
 x (q x+q+1) & x \left(q^2+q x+q+1\right) \left(q^2 x+q^2+q+1\right) &\cdots\\
 \vdots&\vdots&\ddots
\end{array}
\right]\ee
which contains the $2\times 2$ minor
\be\pi_{0,0}\pi_{1,1}-\pi_{1,0}\pi_{0,1}\;=\;2 q^2 x^2 + 2 q^3 x^2 + q^4 x^2 - q x^3 + 2 q^3 x^3 + q^4 x^3 - 
 q^2 x^4 + q^3 x^4\not\myge 0,\ee
so our generalisation does not preserve coefficientwise Hankel-total positivity.
By studying the minor above, however, one might suppose that shifting $q$ by $1$ (that is, replacing $q$ with $1+r$) might yield a polynomial sequence $(\bar{\fm}_n(x;1+r))_{n\geq 0}$ that is coefficientwise Hankel-totally positive; indeed, computer experiments conducted with Alan Sokal seem to confirm this for $n\leq 10$. We conjecture:

\begin{conjecture}[with Alan Sokal]\label{conj.modFq}
The polynomial sequence $(\bar{\fm}_n(x;1+r))_{n\geq0}$ is coefficientwise Hankel-totally positive jointly in $x$ and $r$.
\end{conjecture}

In the following subsections we consider further generalisations of the forests matrix to more indeterminates, and a common theme begins to emerge; in each case we have observed that modifying these matrices in a similar way to~\eqref{eq.barf} above and shifting the indeterminates by $1$ yields a set of matrices with row-generating polynomials that appear to be coefficientwise Hankel-totally positive jointly in a number of indeterminates. 

The polynomials $\bar{\fm}_{n}(x;q)$ are, in fact, specialisations of a $q$-generalisation of the \emph{general Abel polynomials} presented in~\cite{Cigler_08,Jackson_1910,Johnson_96}:
\begin{eqnarray}\label{eq.abelRotheq}a_0(x;b,h,w,q)&\ceq&1\\ a_n(x;b,h,w,q)&\ceq&(x+b)\prod_{i=1}^{n-1}(xq^i+b+\qn{i}h+\qn{n}w)\bigskip\textrm{    for }n\geq 1.\end{eqnarray}
Note that for $w=1$, $b=h=0$ the above formula gives a $q$-generalisation of the Stirling cycle polynomials, so the polynomials $a_n(x;b,h,w,q)$ interpolate between $q$-forest polynomials and a $q$-generalisation of the Stirling cycle polynomials (below we will consider a \emph{different} $q$-generalisation of the Stirling cycle polynomials).

Johnson~\cite[Section~4]{Johnson_96} showed that the polynomials $a_n(x;b,h,w,q)$ are of $q$\emph{-binomial type}, that is, they satisfy
\begin{multline} a_n(x;a+b,h,w,q)\;=\;\sum_{k=0}^n\qbin{n}{k}a_{n-k}(0;a,h+(1-q)(b+\qn{k}w),w(1-(1-q)\qn{k}))\\\times a_k(x;a,h,w,q) \end{multline}
and
\begin{multline}
a_n(x+y;a+b,0,w,q)\;=\;\sum_{k=0}^n\qbin{n}{k}a_{n-k}(yq^k;a,(1-q)(b+\qn{k}w),w(1-(1-q)\qn{k}))\\\times a_k(x;a,0,w,q).\end{multline}
The identities above are $q$-extensions of an identity that goes back to Rothe~\cite{Rothe93} in 1793 and Pfaff~\cite{Pfaff95} in 1795, usually expressed as:
\be R_n(x+y;h,w)\;=\;\sum_{k=0}^n\binom{n}{k}R_k(x;h,w)R_{n-k}(y;h,w)\ee
where
\be R_0(x;h,w)\ceq1\ee
and
\be R_n(x;h,w)\ceq x\prod_{i=1}^{n-1}(x+ih+nw)\ee
for $n>0$.

In~\cite[Section~6]{Sokal_21f} Sokal explains how these polynomials are a rewriting of the \emph{Schl\"afli-Gessel-Seo} polynomials:
\be P_0(x;a,b)\ceq 1\ee
and
\be P_n(x;a,b)\ceq x\prod_{i=1}^{n-1}(x+ai+(n-i)b)\ee
which were introduced by Schl\"afli in 1847~\cite{Schlafli47}, and resurfaced much later in a 2006 paper by Gessel and Seo~\cite{Gessel_06} who showed that 
\be P_{n,k}(a,b)\;=\;[x^k]P_n(x;a,b)\ee 
enumerates $k$-component forests on the vertex set $[n]$ by \emph{proper} and \emph{improper} vertices, and also by \emph{ascents} and \emph{descents}. Sokal's paper contains a fuller discussion of these polynomials as well as a number of interesting conjectures regarding their coefficientwise Hankel-total positivity, and we enthusiastically encourage the reader to consult the final section of~\cite{Sokal_21f} for further details.

There is, unsurprisingly, something of an overlap between this current paper and the generalisations of the forest matrix studied by Sokal~\cite{Sokal_21f}. In the next subsection we will show how our planar network $\scrn_2$ from Section~\ref{sec.coeffpos} can be modified to derive some of the results proven in~\cite{Sokal_21f}.

\subsection{Enumerating forests by proper and improper edges}

In~\cite{Sokal_21f} Sokal considers total positivity properties of generalisations of the forests matrix that differs from the $q$-generalisations studied in this paper. We will now show how making some simple modifications to the planar network $\scrn_2$ yields a path matrix that agrees with one of the generalisations studied in~\cite{Sokal_21f}.

Given two vertices $i$ and $j$ of a tree $T$ belonging to a forest, we say that $j$ is a \emph{descendant} of $i$ if the unique path from the root of $T$ to $j$ passes through $i$ (note that every vertex is a descendant of itself). Let $e=ij$ be an edge of $T$. We say that $e$ is \emph{improper} if there is a descendant of $j$ (possibly $j$ itself) that is lower-numbered than $i$, otherwise we say that $e$ is \emph{proper} (see Section~1 of~\cite{Sokal_21f}).

Consider the matrix
\be\fm(y,z)\;\ceq\;(f_{n,k}(y,z))_{n,k\geq0},\ee
the entries of which are given by
\be f_{n,k}(y,z)\;\ceq\;\sum_{m=0}^{n-k}f_{n,k,m}y^mz^{n-k-m}\ee
where $f_{n,k,m}$ is the number of forests of rooted trees on the vertex set $[n]$ that have $k$ components and $m$ improper edges. The first few rows of $\fm(y,z)$ are:
\be \left[
\begin{array}{cccccc}
 1 &  &  &  &  &\\
 0 & 1 &  &  &  &\\
 0 & y+z & 1 &  &&  \\
 0 & 3 y^2+4 y z+2 z^2 & 3 y+3 z & 1 & & \\
 0 & 25 y^2 z+15 y^3+18 y z^2+6 z^3 & 15 y^2+22 y z+11 z^2 & 6 y+6 z & 1& \\
 \vdots&\vdots&\vdots&\vdots&\vdots&\ddots\\
\end{array}
\right],\ee
each entry is a homogeneous polynomial of degree $n-k$, and under the specialisation $y=z=1$ we recover the forests matrix (that is, $\fm(1,1)=\fm$). Sokal also observes (see the remark on page~7 of~\cite{Sokal_21f}) that the polynomials $f_{n,1}(y,z)$ (which enumerate rooted trees with respect to improper edges) are homogenised versions of the \emph{Ramanujan polynomials} (see~\cite{ChenYang21,Dumont96,Guo18,Guo07,Jousat15,Lin14,Randazzo21,Shor95,Zeng99} and~\cite[A054589]{OEIS}). 

The matrix $\fm(y,z)$ is the exponential Riordan array $\scrr[F,G']$ with $F(t)=1$ and
\be G'(t)\;=\;\frac{1}{z}\left[G\left(\left(1-\frac{z}{y}+\frac{z^2}{y}t\right)e^{-\left(1-\frac{z}{y}\right)}\right)-\left(1-\frac{z}{y}\right)\right]\ee
(Section~3.2 of~\cite{Sokal_21f}) where $G(t)$ is the tree function (see~\cite{Corless96} and~\eqref{eq.treefunc} above). Sokal proved that $\fm(y,z)$ is coefficientwise totally positive in $\Z[y,z]$ using the production-matrix method; the production matrix for $\fm(y,z)$ is the lower-Hessenberg matrix\footnote{A matrix is lower-Hessenberg if all entries above the super-diagonal are $0$.} $\mathsf{P}(y,z)\ceq(p_{n,k}(y,z))_{n,k\geq0}$ with entries
\be p_{n,k}(y,z) \;=\;\frac{n!}{(k-1)!}\sum_{l=0}^{n+1-k}\frac{y^{n-l}z^l}{l!}\ee
(see~\cite[Proposition~4.7]{Sokal_21f}). The first few rows of $\mathsf{P}(y,z)$ are:
\be\left[
\begin{array}{cccccc}
 0 & 1 &  &  & & \\
 0 & y+z & y &  & & \\
 0 & 2 y^2+2 y z+z^2 & 2 y^2+2 y z & y^2 & & \\
 0 & 6 y^2 z+6 y^3+3 y z^2+z^3 & 6 y^2 z+6 y^3+3 y z^2 & 3 y^2 z+3 y^3 & y^3 &\\
 \vdots&\vdots&\vdots&\vdots&\vdots&\ddots
\end{array}
\right],\ee
and $\mathsf{P}(y,z)$ has the factorisation~\cite[Proposition~4.8]{Sokal_21f}
\be \mathsf{P}(y,z)\;=\;B_z\diag((n!)_{n\geq0})\tp(y)\diag((k!)_{k\geq0})^{-1}\Delta,\ee
where $\tp(y)$ is the Toeplitz matrix of powers of $y$ (see Subsection~\ref{subsec.MatFac} above) and $\Delta\ceq (\delta_{n+1,k})_{n,k\geq 0}$ is the lower-Hessenberg matrix with $1$ on the superdiagonal and $0$ elsewhere.

\begin{figure}\centering\includegraphics[scale=0.43]{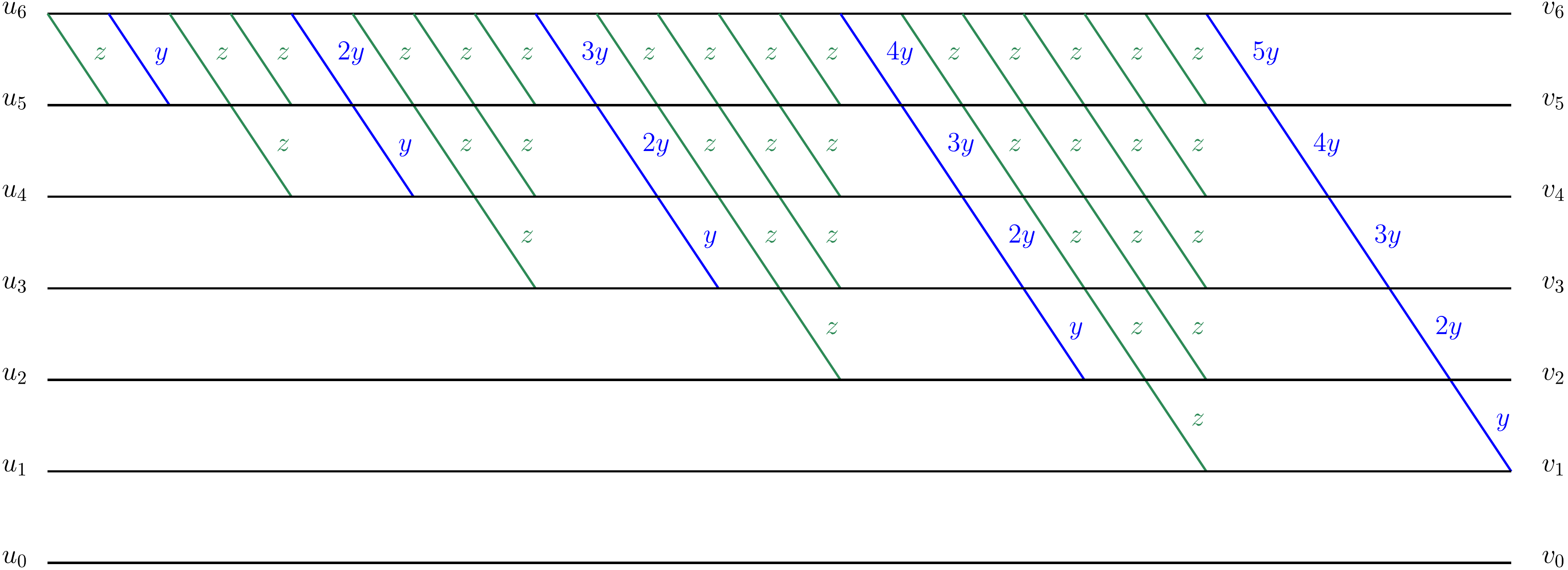}\caption{The planar network $\scrn_3$ arising from the production-matrix approach used in~\cite{Sokal_21f} (unlabelled edges have weight 1).}\label{fig.alanProd}\end{figure}

As described above (see~\eqref{eq.augProd}), the augmented production matrix
\be \tilde{\mathsf{P}}(y,z)\;\ceq\;\left[\begin{array}{c}1\,\,0\,\,0\,\,0\,\,\cdots\\\hline \mathsf{P}(y,z)\end{array}\right]\ee with $\mathsf{P}(y,z)$ defined previously yields the following factorisation of $\fm(y,z)$:

\be\label{eq.fyz} \fm(y,z)\;=\;\cdots\left[\begin{array}{c|c}
            I_3 & 0\\\hline
            0 & B_zT({\bf y}^*)\end{array}\right]\cdot\left[\begin{array}{c|c}
            I_2 & 0\\\hline
            0 & B_zT({\bf y}^*)\end{array}\right]\cdot\left[\begin{array}{c|c}
            1 & 0\\\hline
            0 & B_zT({\bf y}^*)\end{array}\right]\ee
where ${\bf y}^*=((n+1)y)_{n\geq 0}$ [equivalently, $T({\bf y}^*)=D\diag((n!)_{n\geq0})\tp(y)\diag((k!)_{k\geq 0}))^{-1}$] and $B_z=B_{z,1}$ is the weighted binomial matrix $B_{x,y}$ with $x=z$ and $y=1$. A diagram of the corresponding planar network $\scrn_3$ can be found in Figure~\ref{fig.alanProd}. 

Now let $\scrn_4$ denote the planar network in Figure~\ref{fig.zy}, the corresponding path matrix of which has the production-like factorisation
\be P_{\scrn_4}\;=\;\cdots\left[\begin{array}{c|c}
            I_3 & 0\\\hline
            0 & L(3z)T({\bf y}^*)\end{array}\right]\cdot\left[\begin{array}{c|c}
            I_2 & 0\\\hline
            0 & L(2z)T({\bf y}^*)\end{array}\right]\cdot\left[\begin{array}{c|c}
            1 & 0\\\hline
            0 & L(z)T({\bf y}^*)\end{array}\right]\ee
where $L(iz)$ denotes the lower-bidiagonal matrix with 1 on the diagonal, $iz$ on the subdiagonal, and $0$ everywhere else. Comparing Figure~\ref{fig.coeffFq} in Section~\ref{sec.coeffpos} with Figure~\ref{fig.zy} it is easy to see that we obtain $\scrn_4$ from $\scrn_2$ by specialising $q=1$ and multiplying the weights of the blue and green edges by $y$ and $z$ respectively. We make the following proposition:
\begin{figure}\centering\includegraphics[scale=0.4]{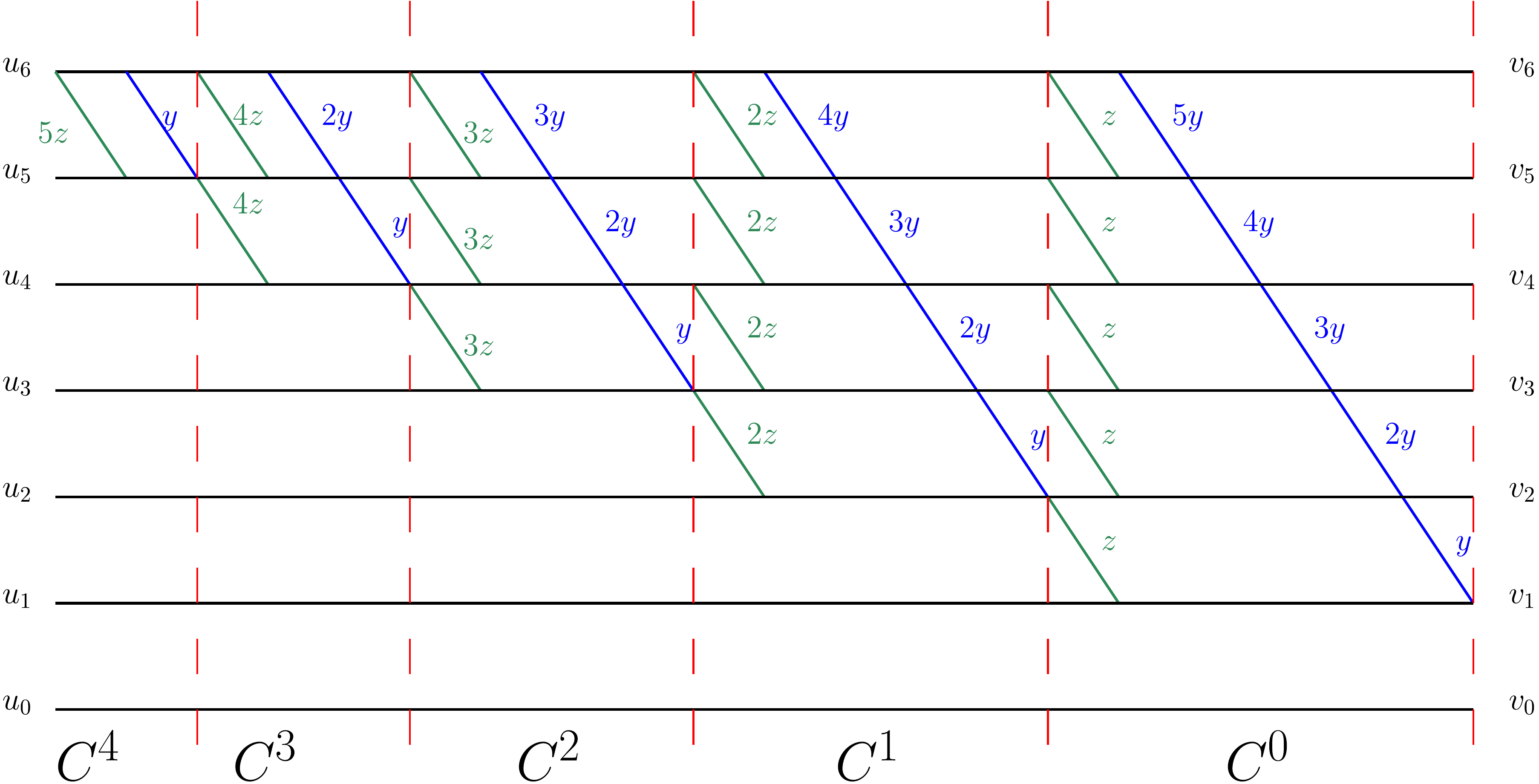}\caption{The planar network $\scrn_4$ up to source and sink vertices $u_6$ and $v_6$.}\label{fig.zy}\end{figure}
\begin{proposition}\label{prop.sokal}
The path matrix corresponding to $\scrn_4$ satisfies
\be P_{\scrn_4} \;=\;\fm(y,z).\ee
\end{proposition}

\begin{proof}
Consider the factorisation of $\fm(y,z)$ in~\eqref{eq.fyz} above:
\be\fm(y,z)\;=\;\cdots\left[\begin{array}{c|c}
            I_3 & 0\\\hline
            0 & B_zT({\bf y}^*)\end{array}\right]\cdot\left[\begin{array}{c|c}
            I_2 & 0\\\hline
            0 & B_zT({\bf y}^*)\end{array}\right]\cdot\left[\begin{array}{c|c}
            1 & 0\\\hline
            0 & B_zT({\bf y}^*)\end{array}\right].\ee
The matrix $B_z$ has the production-like factorisation
\be\label{eq.BzLz} B_z\;=\;\cdots\left[\begin{array}{c|c}
            I_3 & 0\\\hline
            0 & L(z)\end{array}\right]\cdot\left[\begin{array}{c|c}
            I_2 & 0\\\hline
            0 & L(z)\end{array}\right]\cdot\left[\begin{array}{c|c}
            1 & 0\\\hline
            0 & L(z)\end{array}\right]L(z)\;=\;\left[\begin{array}{c|c}
            1 & 0\\\hline
            0 & B_z\end{array}\right]\cdot L(z)\ee
(see Corollary~\ref{cor.binomprodLike} with $y=1$ and $x$ replaced with $z$, and also observe that the right-hand side above implies that $L(z)$ is the augmented production matrix of $B_z$). It follows that
\be\fm(y,z)\;=\;\cdots\left[\begin{array}{c|c}
            I_3 & 0\\\hline
            0 & B_zT({\bf y}^*)\end{array}\right]\cdot\left[\begin{array}{c|c}
            I_2 & 0\\\hline
            0 & B_zT({\bf y}^*)B_z\end{array}\right]\cdot\left[\begin{array}{c|c}
            1 & 0\\\hline
            0 & L(z)T({\bf y}^*)\end{array}\right].\ee
A straightforward calculation confirms that
\begin{multline} T({\bf y}^*)B_z\;=D\tp(y)D^{-1}D\tp((z^n/n!)_{n\geq0})D^{-1}\;=\;D\tp((z^n/n!)_{n\geq0})D^{-1}D\tp(y)D^{-1}\\\;=\;B_zT({\bf y}^*)\end{multline}
where $D=\diag((n!)_{n\geq 0})$, since Toeplitz matrices matrices commute: 
\be \tp(\bfa)\tp(\bfb)\;=\;\tp(\bfa*\bfb)\;=\;\tp(\bfb)\tp(\bfa)\ee
where $\bfa *\bfb$ is the convolution of the sequences $\bfa$ and $\bfb$
\be \bfa * \bfb \;=\; \left(\sum_{k=0}^na_{k}b_{n-k}\right)_{n\geq0}.\ee
Furthermore, it is easily verified that
\be B_zB_{z'}\;=\;B_{z+z'}\;=\;\left[\begin{array}{c|c}
            1 & 0\\\hline
            0 & B_{z+z'}\end{array}\right]\cdot L(z+z'),\ee
hence
\be\fm(y,z)\;=\;\cdots\left[\begin{array}{c|c}
            I_3 & 0\\\hline
            0 & B_zT({\bf y}^*)B_{2z}\end{array}\right]\cdot\left[\begin{array}{c|c}
            I_2 & 0\\\hline
            0 & L(2z)T({\bf y}^*)\end{array}\right]\cdot\left[\begin{array}{c|c}
            1 & 0\\\hline
            0 & L(z)T({\bf y}^*)\end{array}\right].\ee
Iteratively commuting $T({\bf y}^*)$ with $B_{iz}$ and writing
\be B_zB_{iz}\;=\;\left[\begin{array}{c|c}
            1 & 0\\\hline
            0 & B_{(i+1)z}\end{array}\right]L((i+1)z)\ee
for increasing $i$ in the factorisation of $\fm(y,z)$ given above completes the proof.
\end{proof}

\begin{figure}\includegraphics[scale=0.35]{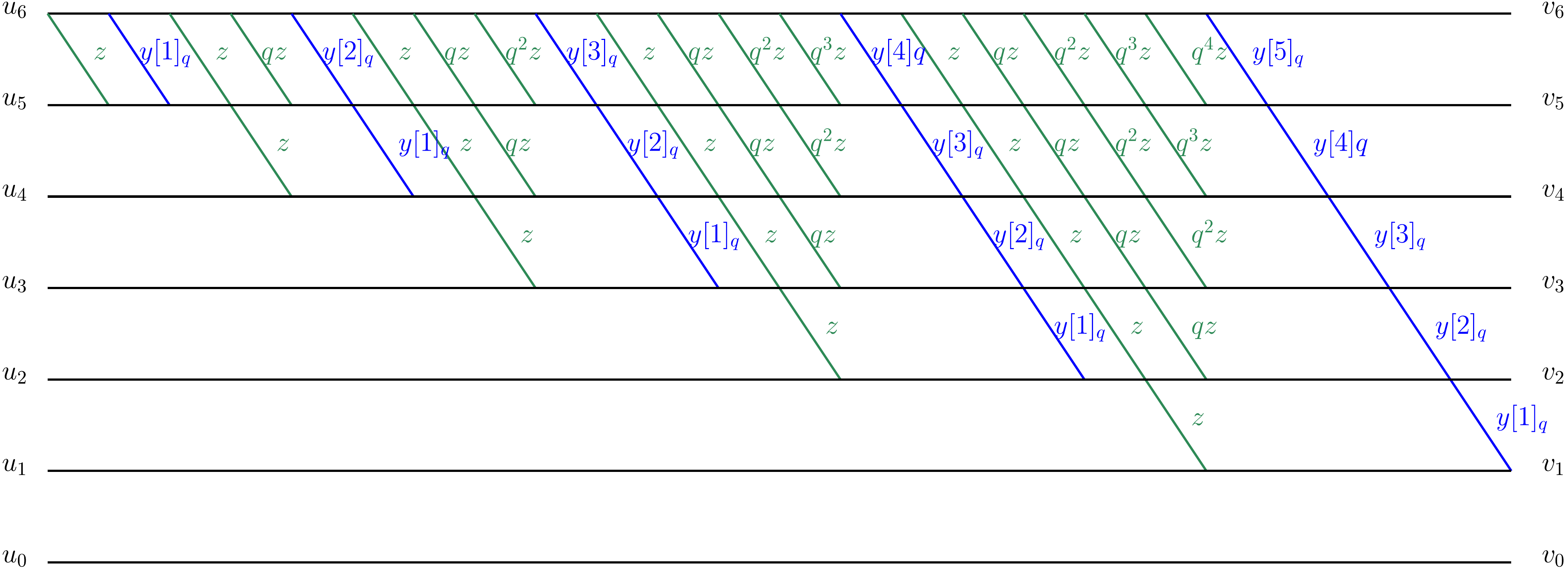}\caption{The planar network $\scrn_5$, obtained by introducing $q$-weights into the network $\scrn_3$.}\label{fig.alan2}\end{figure}

The focus of this paper has been $q$-generalisations of the forests and trees matrices, and note that we can easily introduce $q$-weights into the planar network $\scrn_3$. Consider the network $\scrn_5$ in Figure~\ref{fig.alan2}, which has the corresponding path matrix
\begin{multline} P_{\scrn_5}\;=\;\cdots\left[\begin{array}{c|c}
            I_3 & 0\\\hline
            0 & B_{q;z}T((y[n+1]_q)_{n\geq 0})\end{array}\right]\cdot\left[\begin{array}{c|c}
            I_2 & 0\\\hline
            0 & B_{q;z}T((y[n+1]_q)_{n\geq 0})\end{array}\right]\\\cdot\left[\begin{array}{c|c}
            1 & 0\\\hline
            0 & B_{q;z}T((y[n+1]_q)_{n\geq 0})\end{array}\right]\end{multline}
where $B_{q;z}$ is the $q$-binomial matrix
\be B_{q;z}\;=\;\left(\qbin{n}{k}z^{n-k}\right)_{n,k\geq 0}.\ee
Now let
\be \hat{\fm}(q,y,z)\;=\;(\hat{f}_{n,k}(q,y,z))_{n,k\geq 0}\;=\;P_{\scrn_5}\ee
and observe that $\hat{\fm}(q,y,z)$ is the output matrix of the lower-Hessenberg production matrix $\mathsf{P}'(q,y,z)=(p_{n,k}(q,y,z))_{n,k\geq 0}$ where
\be\label{eq.pd} \mathsf{P}'(q,y,z)\;=\;D_q\tp((z^n/[n]_q!)_{n\geq 0})\tp(y)D_q^{-1}\Delta\ee
where $D_q=\diag(([n]_q!)_{n\geq 0})$.

Define the row-generating polynomials of $\hat{\fm}(q,y,z)$ to be:
\be \hat{\fm}_0(x;q,y,z)\;=\;1\ee
and
\be \hat{\fm}_n(x;q,y,z)\;=\;\sum_{k=0}^n\hat{f}_{n,k}(q,y,z)x^k\ee
for $n\geq 1$. We have the following lemma:
\begin{lemma}\label{lem.alanqyz}
The sequence $(\hat{\fm}_n(x;q,y,z))_{n\geq 0}$ of row-generating polynomials of $\hat{\fm}(q,y,z)=\scro(\mathsf{P}'(q,y,z))$ is coefficientwise Hankel-totally positive jointly in $q,x,y,z$.
\end{lemma}
In order to prove Lemma~\ref{lem.alanqyz} we first recall two easy results from~\cite{Sokal_21f}. The first is:
\begin{lemma}[Lemma~2.6 in~\cite{Sokal_21f}]\label{lem.alanbpb}
Let $P=(p_{i,j})_{i,j\geq 0}$ be a row-finite matrix (with entries in a commutative ring $R$) with output matrix $A=\scro(P)$; and let $B=(b_{i,j})_{i,j\geq 0}$ be a lower-triangular matrix with invertible (in $R$) diagonal entries. Then
\be AB\;=\;b_{0,0}\scro(B^{-1}PB).\ee
That is, up to a factor $b_{0,0},$ the matrix $AB$ has production matrix $B^{-1}PB$.
\end{lemma}
The second is the following theorem:
\begin{theorem}[Theorem~2.14 in~\cite{Sokal_21f}]\label{thm.alanh0}
Let $P= (p_{i,j})_{i,j\geq0}$ be an infinite row-finite or column-finite matrix with entries in a partially ordered commutative ring $R$, and define the infinite Hankel matrix $H_{\infty}(\scro_0(P)) = ((P^{n+n'})_{0,0})_{n,n'\geq 0}$. If $P$ is totally positive of order $r$, then so is $H_{\infty}(\scro_0(P))$.\end{theorem}

Our proof of Lemma~\ref{lem.alanqyz} relies on these two useful facts, and follows along the same lines as the argument found in the proof of Lemma~2.16 of~\cite{Sokal_21f} where $D$ is replaced with $D_q$, $B_x$ is replaced with $B_{q;x}$, and
\be \bsym\phi\;=\;(\phi_n)_{n\geq0}\;=\;\sum_{l=0}^n\frac{z^ly^{n-l}}{\qn{l}!}.\ee
\begin{proof}[Proof of Lemma~\ref{lem.alanqyz}]
The row-generating polynomials of $\hat{\fm}(q,y,z)$ are the entries in the zeroth column of the matrix
\be \hat{\fm}^*(q,x,y,z)\;=\;\hat{\fm}(q,y,z)B_{q;x},\ee
and according to Lemma~\ref{lem.alanbpb} above we have
\be \hat{\fm}^*(q,x,y,z)\;=\;\scro(B_{q;x}^{-1}\mathsf{P}'(q,y,z)B_{q;x}).\ee
Since $\mathsf{P}'(q,y,z)$ is coefficientwise totally positive it suffices (thanks to Theorem~\ref{thm.alanh0}) to show that the matrix 
\be B_{q;x}^{-1}\mathsf{P}'(q,y,z)B_{q;x}\ee
is coefficientwise totally positive.
We have
\be B_{q;x}\;=\;D_q\tp((x^n/[n]_q!)_{n\geq 0})D_q^{-1},\ee
and
\be \mathsf{P}'(q,y,z)\;=\;D_q\tp(\bsym\phi)D_q^{-1}\Delta,\ee
where 
\be \bsym\phi\;=\;(\phi_n)_{n\geq0}\;=\;\sum_{l=0}^n\frac{z^ly^{n-l}}{\qn{l}!}.\ee
It follows that $B_{q;x}$ and $D_q\tp(\bsym\phi)D_q^{-1}$ commute.

Note further that according to the recurrences for $q$-binomial coefficients (see ~\eqref{eq.binomIdentities1} and~\eqref{eq.binomIdentities2} in Section~\ref{sec.combin.interp}) we have
\be \Delta B_{q;x}\;=\;B_{q;x}(\diag((q^nx)_{n\geq0})+\Delta).\ee
We therefore have
\begin{align*}B_{q;x}^{-1}\mathsf{P}'(q,y,z)B_{q;x}&=B_{q;x}^{-1}D_q\tp(\bsym\phi)D_q^{-1}\Delta B_{q;x}\\
&=B_{q;x}^{-1}D_q\tp(\bsym\phi)D_q^{-1}B_{q;x}(\diag((q^nx)_{n\geq0})+\Delta)\\
&=D_q\tp(\bsym\phi)D_q^{-1}(\diag((q^nx)_{n\geq0})+\Delta)\\
&=D_q\tp(\bsym\phi)D_q^{-1}\Delta (I+\Delta^T\diag((q^nx)_{n\geq 0}))\end{align*}
(since $\Delta\Delta^T=I$), yielding the identity
\be B^{-1}_{q;x}\mathsf{P}'(q,y,z)B_{q;x}\;=\;\mathsf{P}'(q,y,z)(I+\Delta^T\diag((q^nx)_{n\geq 0})),\ee
so $B^{-1}_{q;x}\mathsf{P}'(q,y,z)B_{q;x}$ coefficientwise totally positive. This completes the proof.
\end{proof}
My experiments suggest that there is a more general principle lurking behind Lemma~\ref{lem.alanqyz}, concerning production matrices of $q$-Riordan arrays (see ~\cite{Cheon2013}, and also~\cite{Johnson96,Gessel062}). This will be rigorously addressed in future work; for now we turn our attention to a more general planar network that yields both $\scrn_2$ and $\scrn_4$ under suitable specialisations. 
\begin{figure}\includegraphics[scale=0.6]{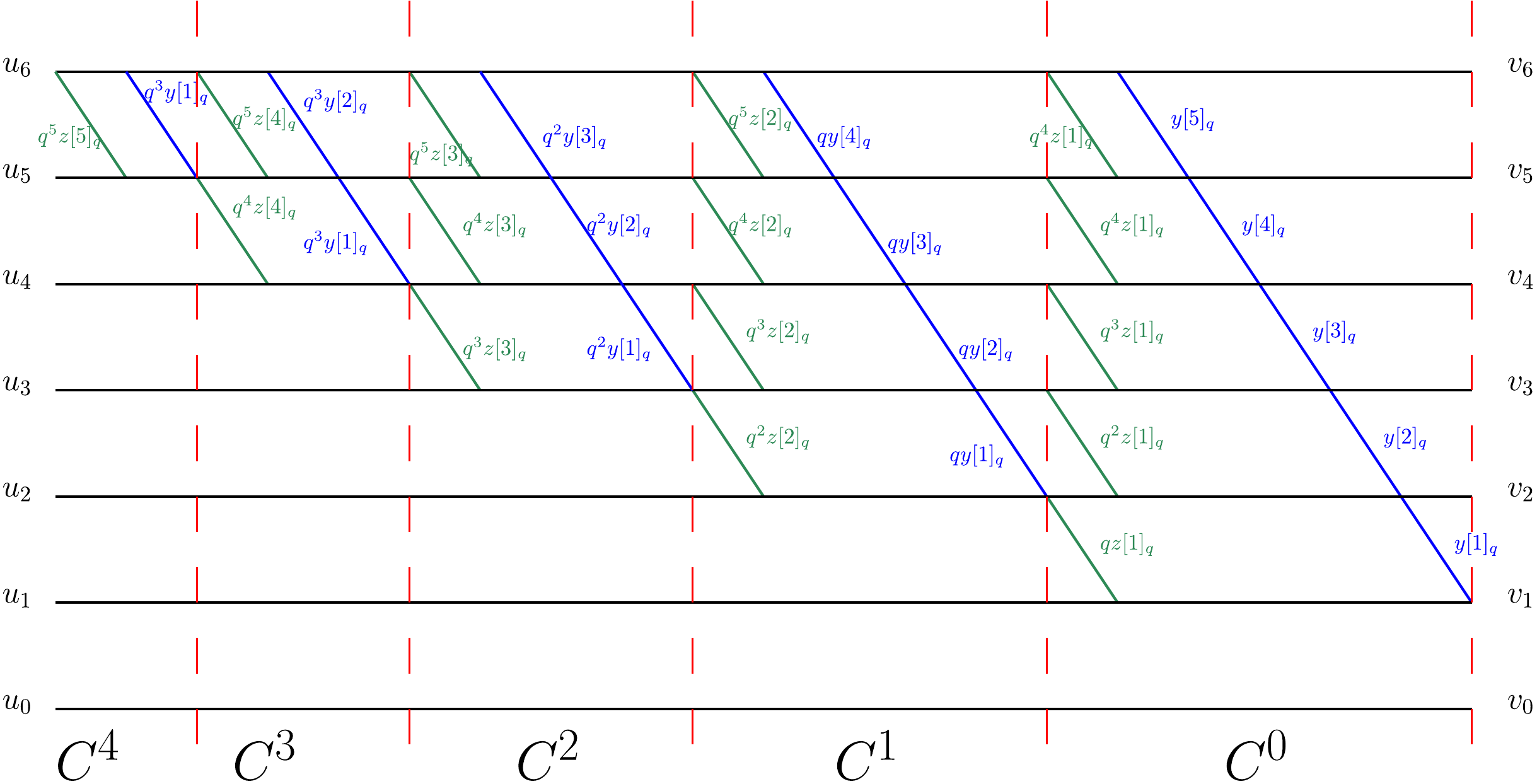}\caption{The planar network $\scrn_6$.}\label{fig.qprod_zyq}\end{figure}

Consider now the planar network $\scrn_6$ in Figure~\ref{fig.qprod_zyq} and observe that $\scrn_4$ can be obtained from $\scrn_6$ by specialising $q=1$ in $\scrn_6$ to $q=1$, and by instead specialising $y=z=1$ in $\scrn_6$ we obtain the planar network $\scrn_2$. The path matrix corresponding to $\scrn_6$ has the factorisation:

\begin{multline} P_{\scrn_6}\;=\;\cdots\left[\begin{array}{c|c}
            I_3 & 0\\\hline
            0 & L(\bfa_2(q,z))T(\bfb_2(q,y))\end{array}\right]\cdot\left[\begin{array}{c|c}
            I_2 & 0\\\hline
            0 & L(\bfa_1(q,z))T(\bfb_1(q,y))\end{array}\right]\\\cdot\left[\begin{array}{c|c}
            1 & 0\\\hline
            0 & L(\bfa_0(q,z))T(\bfb_0(q,y))\end{array}\right]\end{multline}
where 
\be\bfa_i(q,z)\;=\;(q^{n+1+i}z[i+1]_q)_{n\geq 0}\ee
and
\be\bfb_i(q,y)\;=\;(q^iy[n+1]_q)_{n\geq 0}.\ee
Let
\be \tilde{\fm}(q,y,z)\;=(\tilde{\fm}_{n,k}(q,y,z))_{n,k\geq 0}\;=\;P_{\scrn_6}\ee
and note that $\tilde{\fm}(q,y,z)$ is a generalisation of $\fm(q)$ and $\fm(y,z)$ in the sense that $\tilde{\fm}(1,y,z)=\fm(y,z)$, and $\tilde{\fm}(q,1,1)=\fm(q)$.

Sokal showed in~\cite{Sokal_21f} that the row-generating polynomials of $\fm(y,z)$:
\be\fm_0(x;y,z)\ceq1\ee
and
\be\fm_n(x;y,z)\ceq\sum_{k=0}^nf_{n,k}(y,z)x^k\ee
for $n>0$ are coefficientwise Hankel-totally positive jointly in $x,y,z$ (Theorem~1.3 of~\cite{Sokal_21f}). One might wonder, then, whether the sequence $(\tilde{\fm}_n(x;q,y,z))_{n\geq 0}$ of row-generating polynomials of $\tilde{\fm}(q,y,z)$:
\be \tilde{\fm}_0(x;q,y,z)\;=\;1\ee
and
\be \tilde{\fm}_n(x;q,y,z)\;=\;\sum_{k=0}^{n}\tilde{f}_{n,k}(q,y,z)x^k\ee
is also Hankel totally-positive. Unfortunately this is not the case, for the Hankel matrix $H_{\infty}((\tilde{\fm}_{n+n'})_{n,k\geq 0})=(\pi_{n,k})_{n,k\geq 0}$ contains the $2\times 2$ minor:
\begin{multline} \pi_{1,0}\pi_{2,1}-\pi_{2,0}\pi_{1,1}=q^3 x^2 y z+3 q^2 x^2 y z+q^4 x^2 z^2+q^3 x^2 z^2-q^2 x^2 z^2+q^3 x^3 z+2 q^2 x^3 z+2 q x^2 y^2\\-2 q x^2 y z+2 q x^3 y-2 q x^3 z-x^3 y\not\myge0.\end{multline}
However, it appears that coefficientwise total positivity can be restored by multiplying $\tilde{\fm}(q,y,z)$ by a simple diagonal matrix. Consider the matrix
\be \bar{\fm}(q,y,z)\;=(\bar{f}_{n,k}(q,y,z))_{n\geq 0}\;=\;\tilde{\fm}(q,y,z)\diag((q^{k(k-1)/2})_{k\geq0}),\ee
the row-generating polynomials of which are
\be \bar{\fm}_0(x;q,y,z)\;=\;1\ee
and
\be\bar{\fm}_n(x;q,y,z)=\sum_{k=0}^n\bar{f}_{n,k}(q,y,z)x^k\ee
for $n\geq 1$. The matrix $\bar{\fm}(q,y,z)$ is manifestly totally positive, and the sequence of row-generating polynomials of $\bar{\fm}(q,y,z)$ appear empirically (up to $5\times 5$) to be coefficientwise totally positive jointly in $q,x,y,z$. I conjecture:
\begin{conjecture}
The Hankel matrix $H_{\infty}(\bar{\fm}_{n+n'}(x;q,y,z))_{n,n'\geq 0})$ is coefficientwise totally positive jointly in $q,x,y,z$.
\end{conjecture}

But what of a combinatorial interpretation of the entries of $\bar{\fm}(q,y,z)$? As we have seen, on the one hand the entries of $\bar{\fm}(1,y,z)$ count forests on $n$ vertices with $k$ components with respect to proper and improper children, and on the other the entries of $\bar{\fm}(q,1,1)$ count forests on $n$ vertices with $k$ components with respect to the statistics defined in Section~\ref{sec.combin.interp}. I therefore pose the following problem:
\begin{problem}
Find a combinatorial interpretation of the entries of $\bar{\fm}(q,y,z).$
\end{problem}

\subsection{A generalisation with more variables}
\begin{figure}\centering\includegraphics[scale=0.42]{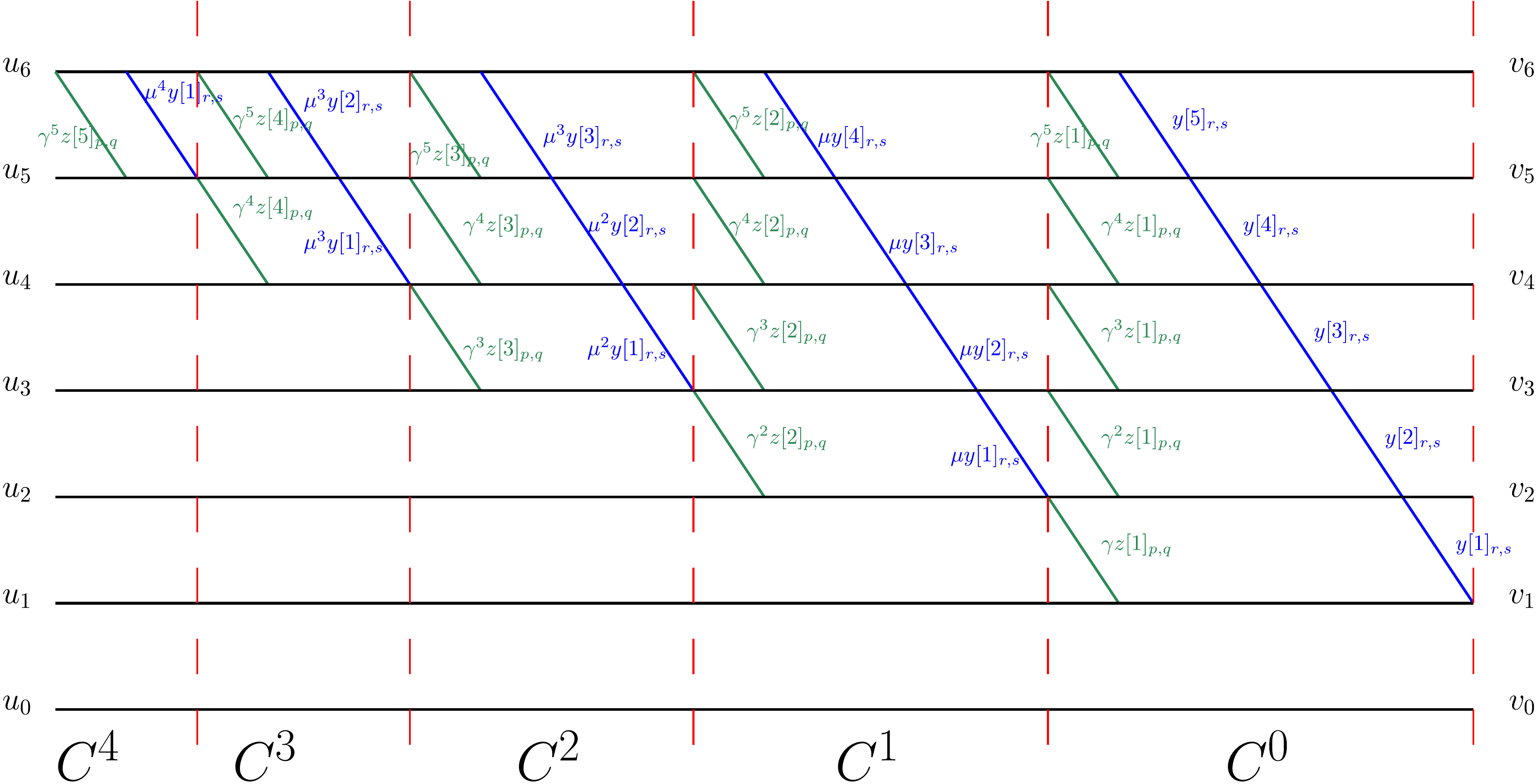}\caption{The planar network $\scrn_7$ up to source and sink vertices $u_6$ and $v_6$.}\label{fig.qSix}\end{figure}

Consider now the planar network $\scrn_7$ in Figure~\ref{fig.qSix}. The path matrix corresponding to $\scrn_7$ is
\be P_{\scrn_7}\ceq (f_{n,k}(p,q,r,s,\gamma,\mu,y,z))_{n,k\geq 0}\ee
which has the factorisation
\be P_{\scrn_7}\;=\;\cdots P_{C^2}P_{C^1}P_{C^0}\ee
where
\be P_{C^i}\ceq\left[\begin{array}{c|c}
            I_{i+1} & 0\\\hline
            0 & L(\bfa_i(p,q,\gamma,z))T(\bfb_i(r,s,\mu,y))\end{array}\right]\ee
in which $L(\bfa_i(p,q,\gamma,z))$ is the lower-bidiagonal matrix with edge sequences
\be\bfa_i(p,q,\gamma,z)\;=\;(\gamma^{i+n+1}z[i+1]_{p,q})_{n\geq 0}\ee
and $T(\bfb_i(r,s,\mu,y))$ the inverse lower-bidiagonal matrix with edge sequences
\be\bfb_i(r,s,\mu,y)\;=\;(\mu^{i}y[n+1]_{r,s})_{n\geq0}.\ee
Now let
\be \fm(p,q,r,s,\gamma,\mu,y,z)\;\ceq\; (f_{n,k}(p,q,r,s,\gamma,\mu,y,z))_{n,k\geq 0}\;=\;P_{\scrn_7},\ee
and note that $\scrn_7$ reduces to $\scrn_6$ (see Figure~\ref{fig.qprod_zyq}) when $p,r,s,\gamma,\mu$ are specialised to $p=r=1$ and $s=\gamma=\mu=q$, so $\fm(p,q,r,s,\gamma,\mu,y,z)$ is a generalisation of the forests matrix to \emph{eight} indeterminates.

The first few rows of $\fm(p,q,r,s,\gamma,\mu,y,z)$ are:
\be
\left[
\begin{array}{cccc}
 1 &  &  &  \\
 0 & 1 &  &  \\
 0 & \gamma z+y & 1 &  \\
 0 & \gamma^2 p y z+\gamma^3 p z^2+\gamma^2 q y z+\gamma^3 q z^2& \gamma^2 p z+\gamma^2 q z+\gamma^2 z&1\\
 &+\gamma^2 y z+\gamma \mu y z+\mu y^2+r y^2+s y^2 &+\mu y+r y+s y &  \\
 \vdots &\vdots &\vdots &\ddots
\end{array}
\right]\ee
Clearly this matrix is coefficientwise totally positive jointly in $p,q,r,s,\gamma,\mu$, and the entries must count forests of trees on $n$ vertices with $k$ components with respect to \emph{eight} different statistics. I present two main open problems regarding the matrix $\fm(p,q,r,s,\gamma,\mu,y,z)$. The first is:

\begin{problem}
Find a combinatorial interpretation of the entries of $\fm(p,q,r,s,\gamma,\mu,y,z)$.
\end{problem}

We know that specialising $p=r=y=z=1$ and $s=\gamma=\mu=q$ yields the $q$-forests matrix $\fm(q)$, and the entries count forests on $n$ vertices with $k$ components according to the statistics presented in Section~\ref{sec.combin.interp}. Similarly, by specialising $\gamma=\mu=p=q=r=s=1$ we obtain the matrix $\fm(y,z)$, the entries of which count $k$-component forests on $n$ vertices with respect to proper and improper edges. In the following subsections we identify other specialisations of $p,q,r,s,\gamma,\mu,y,z$ yielding matrices with entries that appear to count a variety of interesting combinatorial objects including: permutations with respect to \emph{left-right maxima}, unordered forest of increasing labelled trees, $0$-$1$ tableaux by \emph{inversions} and \emph{noninversions}, and perfect matchings with respect to \emph{crossings} and \emph{nestings}.

Now define the row-generating polynomials of $\fm(p,q,r,s,\gamma,\mu,y,z)$ to be
\be \fm_0(x;p,q,r,s,\gamma,\mu,y,z)\;\ceq\;1\ee
and
\be \fm_n(x;p,q,r,s,\gamma,\mu,y,z)\;\ceq\;\sum_{l=0}^nf_{n,l}(p,q,r,s,\gamma,\mu,y,z)x^l\ee
for $n>0$, and consider the Hankel matrix
\be H_{\infty}((\fm_n(x;p,q,r,s,\gamma,\mu,y,z))_{n\geq0})\ceq(\pi_{n,k})_{n,k\geq0}\ee
where
\be\pi_{n,k}\;=\;\fm_{n+k}(x;p,q,r,s,\gamma,\mu,y,z).\ee

The first few $\fm_n(x;p,q,r,s,\gamma,\mu,y,z)$ are:
\be
\begin{array}{c|c}
 n & \fm_n(x;p,q,r,s,\gamma,\mu,y,z)\\\hline
 0 & 1 \\
 1 & x \\
 2 & x(x+y+\gamma z) \\
 3 & x (x^2 + (\mu + r + s) x y + (\mu + r + s) y^2 + 
   \gamma^2 (1 + p + q) x z \\
   &+ \gamma (\mu + \gamma (1 + p + q)) y z + 
   \gamma^3 (p + q) z^2)\\
   \vdots&\vdots
\end{array}
\ee
so the associated Hankel matrix contains the $2\times 2$ minor:
\begin{multline}\pi_{1,0}\pi_{2,1}-\pi_{2,0}\pi_{1,1}\;=\;\gamma^2 p x^2 y z+\gamma^3 p x^2 z^2+\gamma^2 p x^3 z+\gamma^2 q x^2 y z+\gamma^3 q x^2 z^2+\gamma^2 q x^3 z\\+\gamma^2 x^2 y z-\gamma^2 x^2 z^2+\gamma^2 x^3 z+\gamma \mu x^2 y z-2 \gamma x^2 y z-2 \gamma x^3 z\\+\mu x^2 y^2+\mu x^3 y+r x^2 y^2+r x^3 y+s x^2 y^2+s x^3 y-x^2 y^2-2 x^3 y\not\myge 0.\end{multline}
The polynomial sequence $(\fm_n(x;p,q,r,s,\gamma,\mu,y,z))_{n\geq0}$ is evidently not coefficientwise Hankel-totally positive. However, consider instead the matrix 
\begin{multline}\label{eq.tildefm}\tilde{\fm}(p,q,r,s,\gamma,\mu,y,z)\ceq(\tilde{f}_{n,k}(p,q,r,s,\gamma,\mu,y,z))_{n,k\geq 0}\\\;=\; \fm(p,q,r,s,\gamma,\mu,y,z)\diag(((pqrs\gamma\mu)^{n(n-1)/2})_{n\geq 0}),\end{multline}
the row-generating polynomials of which are
\be \tilde{\fm}_0(x;p,q,r,s,\gamma,\mu,y,z)\ceq1\ee
and
\be \tilde{\fm}_n(x;p,q,r,s,\gamma,\mu,y,z)\ceq\sum_{k=0}^n\tilde{f}_{n,k}(p,q,r,s,\gamma,\mu)x^k\ee
for $n>0$.

The sequence $(\tilde{\fm}_n(x;p,q,r,s,\gamma,\mu,y,z))_{n\geq 0}$ is \emph{not} coefficientwise Hankel-totally positive, however, my computations suggest that shifting the indeterminates $p,q,r,s,\gamma,\mu$ by $1$ yields a polynomial sequence that appears to be coefficientwise Hankel-totally positive jointly in \emph{nine} indeterminates! I make the following conjecture:

\begin{conjecture}
The polynomial sequence 
\be(\tilde{\fm}_n(x;1+p',1+q',1+r',1+s',1+\gamma',1+\mu',y,z))_{n\geq 0}\ee
is coefficientwise Hankel-totally positive jointly in $x,y,z,p',q',r',s',\gamma',\mu'$.
\end{conjecture}

I have verified the above conjecture to $5\times 5$, and I do not yet understand what these polynomials enumerate. The conjectured coefficientwise Hankel-total positivity of the shifted row-generating polynomials of $\tilde{\fm}(p,q,r,s,\gamma,\mu,y,z)$ is in a sense \emph{stronger} than the coefficientwise total positivity of $\tilde{\fm}(p,q,r,s,\gamma,\mu,y,z)$ alone, and in light of this a combinatorial interpretation of the entries of $\tilde{\fm}(p,q,r,s,\gamma,\mu,y,z)$ is much desired. Clearly each $(n,k)$-entry enumerates $k$-component forests on $n$ vertices with respect to eight statistics; interpreting these statistics is left as an open problem:
\begin{problem}
Find a combinatorial interpretation of the entries of~ $\tilde{\fm}(p,q,r,s,\gamma,\mu,y,z)$.
\end{problem}

\subsection{\texorpdfstring{$(p,q)$}{TEXT}-Stirling polynomials}

The network $\scrn_7$ is an interlacing of two different planar networks (one with green diagonal edges, one with blue). By setting $r=s=\mu=y=0$ we effectively remove the blue edges and are left with a binomial-like planar network we shall denote $\scrn_7'$. Specialising $\gamma =z=1$ in $\scrn_7'$ we obtain the binomial-like network with edge weights specialised to $\alpha_{i,l}=1$ and 
\be \beta_{i,l}=[i-1]_{p,q},\ee with corresponding path matrix
\be G(p,q)\ceq (g_{n,k}(p,q))_{n,k\geq 0}\;=\;P_{\scrn_7'}\;=\;\fm(p,q,0,0,1,0,0,1)\ee
The weights of the network $\scrn_7'$ are purely $i$-dependent so the entries of $G(p,q)$ satisfy the recurrence:
\be\label{eq.srecurr} g_{n,k}(p,q)\;=\;g_{n-1,k-1}(p,q)\;+\;[n-1]_{p,q}g_{n-1,k-1}(p,q)\ee
for $n\geq 1$ with initial condition $g_{0,k}=\delta_{0k}$. Moreover we have
\be\label{eq.stirling.symm} g_{n,k}(p,q)\;=\;\bfe_{n-k}([1]_{p,q},[2]_{p,q},\ldots,[n-1]_{p,q})\ee
by way of~\eqref{eq.eSymm} in Section~\ref{sec.prelim}. The recurrence~\eqref{eq.srecurr} can be found in~\cite[Proposition~2.2]{DeMedicis_93} and~\cite{Zeng95}, and the expression in terms of elementary functions~\eqref{eq.stirling.symm} can be found in Proposition~2.3 (part (a)) of~\cite{DeMedicis_93}. 
The entries of $G(p,q)$ are thus the \emph{$(p,q)$-Stirling cycle numbers} (or \emph{unsigned $(p,q)$-Stirling numbers of the first kind}, denoted $c_{p,q}[n,k]$) and the row-generating polynomials of $G(p,q)$ are
\be c_n(x;p,q)\ceq x\prod_{i=1}^{n-1}(x+y[i]_{p,q})\ee
where $y=1$ (see Proposition~2.3 part (b) of~\cite{DeMedicis_93}). Note that the authors of~\cite{DeMedicis_93} have a combinatorial interpretation of $c_{p,q}[n,k]$, but it has nothing to do with forests and trees! Intead they consider $0$-$1$ tableaux: a $0$-$1$ tableau is a pair $\psi=(\lambda,f)$ where $\lambda=(\lambda_1\geq\lambda_2\geq\cdots\geq\lambda_k)$ is a partition of an integer $m=|\lambda|$ and $f=(f_{i,j})_{1\leq j\leq\lambda_i}$ a filling of of the corresponding Ferrer's diagram of $\lambda$ with $0$s and $1$s such that there is exactly one $1$ in each column. In~\cite{DeMedicis_93} de M\'{e}dicis and Leroux define two statistics on $\psi$: the \emph{inversion number}, $\inv(\psi)$, is the number of $0$s below a $1$ in $\psi$, and the \emph{noninversion number}, $\noninv(\psi)$, is the number of $0$s above a $1$ in $\psi$. They then show that
\be c_{p,q}[n,k]\;=\;\sum_{\psi\in Td(n-1,n-k)}p^{\noninv(\psi)}q^{\inv(\psi)}\ee
where $Td(k,r)$ denotes the set of $0$-$1$ tableaux with $k$ rows and columns of length at most $r$, where the lengths of the columns are distinct.

Permutations, however, are perhaps arguably more natural objects to study in relation to generalisations of the Stirling cycle numbers (the polynomial $c_n(x;1,1)$ is the generating function for permutations counted with respect to cycles, after all). For $p$ specialised to $p=1$ Zeng~\cite[Section~3, where \texorpdfstring{$s_q(n,k) = q^{n-k}c_{1,q}[n,k]$}{TEXT}]{Zeng95} showed that
\be c_{n}(x;1,q)\;=\;\sum_{\sigma\in\mathfrak{S}_n}x^{\textrm{lrm}(\sigma)}q^{\inv(\sigma)}\ee
where lrm$(\sigma)$ denotes the number of left-right maxima in a permutation $\sigma$ (that is, the number of $\sigma(i)$ such $\sigma(i)>\sigma(j)$ for all $j<i$), and $\inv(\sigma)$ the number of inversions of $\sigma$ (that is, the number of pairs $i<j$ such that $\sigma(i)>\sigma(j)$), hence
\be g_{n,k}(1,q)\;=\;\sum_{\sigma\in\mathfrak{S}_n(k)}q^{\inv(\sigma)+k-n},\ee
where $\mathfrak{S}_n(k)$ denotes the set of permutations on $[n]$ with $k$ left-right maxima.

Clearly $G(p,q)$ is totally positive, and for $p$ specialised to $p=1$ we can express the row-generating polynomials of $G(1,q)$ as the $S$-type (or \emph{Stieltjes-type}) continued fraction\footnote{The study of continued fractions goes back to Euler, and was reinvigorated more recently by Flajolet in~\cite{Flajolet80}.} 
\be 1+ \sum_{n,k\geq 1}c_{1,q}[n,k]x^kt^n\;=\cfrac{1}{1-\cfrac{\lambda_1t}{1-\cfrac{\lambda_2t}{\cfrac{\cdots}{1-\cfrac{\lambda_nt}{\cdots}}}}}
\ee
where $\lambda_{2n-1}\ceq (x+\qn{n})q^{n-1}$ and $\lambda_{2n}=\qn{n}q^{n+1}$ for $n\geq 1$ (Lemma~3 of Zeng 1993, where we have replaced $x$ with $tq$, and $a$ with $x/q$). If a sequence of polynomials has an $S$-fraction with nonnegative coefficients then the corresponding Hankel matrix is coefficientwise totally positive (see~\cite{Sokal22} and~\cite[Section~9]{petreolle2020}), so we have:

\begin{corollary}
The polynomial sequence $(c_n(x;1,q))_{n\geq 0}$ is coefficientwise Hankel-totally positive in $x,q$.
\end{corollary}

The Hankel matrix of $c_n(x;p,q)$, where $p$ and $q$ are indeterminates, however, begins
\be \left[
\begin{array}{cccc}
 1 & x & x (x+1)&\cdots \\
 x & x (x+1) & x (x+1) (p+q+x)&\cdots \\
 x (x+1) & x (x+1) (p+q+x) & x (x+1) (p+q+x) \left(p^2+p q+q^2+x\right)&\cdots \\
\vdots&\vdots &\vdots&\ddots
\end{array}
\right]\ee
which contains the $2\times2$ minor
\be \pi_{0,1}\pi_{1,2}-\pi_{1,1}\pi_{0,2}\;=\;p x^3+p x^2+q x^3+q x^2-x^3-x^2\not\myge 0\ee
so is not even TP$_2$ in $\Z[p,q]$. Again it appears, however, that right-multiplying $G(p,q)$ by a simple diagonal matrix yields a matrix whose row-generating polynomials appear to be coefficientwise Hankel-totally positive when $p,q$ are shifted by $1$.

Consider the matrix 
\be G'(p,q,\gamma,z)\;=\;(g'_{n,k}(p,q,\gamma,z))_{n,k\geq 0}\;=\;\fm(p,q,0,0,\gamma,0,0,z)\diag(((pq\gamma)^{n(n-1)/2})_{n\geq 0})\ee
and define its row-generating polynomials to be
\be G'_n(x;p,q,\gamma,z)\ceq\sum_{k=0}^ng'_{n,k}(p,q,\gamma,z)x^k.\ee
Clearly $G'(p,q,\gamma,z)$ is coefficientwise totally positive jointly in $p,q,\gamma,z$ (this follows from Cauchy-Binet and the LGV lemma). I make the following conjecture:
\begin{conjecture}
The polynomial sequence $(G'_n(x;p+1,q+1,\gamma+1,z))_{n\geq 0}$ is coefficientwise Hankel-totally positive jointly in $p,q,\gamma,z,x$.
\end{conjecture}
I have verified this up to $7\times 7$, and am yet to discover a suitable combinatorial interpretation of these polynomials. I therefore pose the following problem:
\begin{problem}
Find a combinatorial interpretation of the polynomials $G'_n(x;p,q,\gamma,z)$.
\end{problem}

\subsection{Generalised Bessel polynomials}

Now consider the planar network $\scrn_7''$ obtained by specialising $p=q=\gamma=z=0$ in $\scrn_7$ (effectively removing all diagonal blue edges from the network), and let
\be \mm(r,s,\mu,y)\ceq(m_{n,k}(r,s,\mu,y))_{n,k\geq 0}\;=\;P_{\scrn_7''} \;=\;\fm(0,0,r,s,0,\mu,y,0)\ee
be the corresponding path matrix. The matrix $\mm(r,s,\mu,y)$ has a factorisation in terms of downshifted inverse bidiagonal matrices:
\be \mm(r,s,\mu,y)\;=\;\cdots\left[\begin{array}{c|c}
          I_3 & 0 \\
          \hline
          0 & T(\bfa_2(r,s,\mu,y))
          \end{array}\right]\cdot\left[\begin{array}{c|c}
          I_2 & 0 \\
          \hline
          0 & T(\bfa_1(r,s,\mu,y))
          \end{array}\right]\cdot \left[\begin{array}{c|c}
          1 & 0 \\
          \hline
          0 & T(\bfa_0(r,s,\mu,y))
          \end{array}\right]\ee
where $\bfa_i(r,s,\mu,y)=(\mu^iy[n+1]_{r,s})_{n\geq 0}$.

Each matrix $T(\bfa_i(r,s,\mu,y))$ has $(n,k)$-entry:
\be \mu^{i(n-k)}y\frac{[n]_{r,s}!}{[k]_{r,s}!}\ee
and specialising $r,s,\mu$ to $r=s=\mu=1$ we obtain
\be\label{eq.q1Tarray} \left[\begin{array}{c|c}
          1 & 0 \\
          \hline
          0 & T(\bfa_i(1,1,1,y))
          \end{array}\right]\;=\;\left[
\begin{array}{ccccccc}
 1&&&&&&\\
 0& 1 &  &  &  & & \\
 0& y & 1 &  &  & & \\
 0& 2y^2 & 2y & 1 &  & & \\
 0&6y^3 & 6y^2 & 3y & 1 & & \\
 0& 24y^4 & 24y^3 & 12y^2 & 4y & 1 & \\
\vdots&\vdots&\vdots&\vdots&\vdots&\vdots&\ddots
\end{array}
\right]\ee
(this array appears in~\cite[A094587]{OEIS}).
Compare the above with the production matrix $\mathsf{P}_{\bsym\phi}\ceq (p_{n,k})_{n,k\geq 0}$ defined in Proposition~1.4 of~\cite{Petreolle_2021}, where
\be p_{n,k}\;=\;\begin{cases} 0&\textrm{if }k=0,\\
                              \displaystyle\frac{n!}{(k-1)!}\phi_{n-k+1}&\textrm{if } n\geq k-1.\end{cases}\ee
We have
\be\label{eq.phiMat}\mathsf{P}_{\bsym\phi}\;=\;\left[
\begin{array}{ccccccc}
 0& 1 &  &  &  & & \\
 0& \phi_1 & 1 &  &  & & \\
 0& 2\phi_2 & 2\phi_1 & 1 &  & & \\
 0&6\phi_3 & 6\phi_2 & 3\phi_1 & 1 & & \\
 0& 24\phi_4 & 24\phi_3 & 12\phi_2 & 4\phi_1 & 1 & \\
\vdots&\vdots&\vdots&\vdots&\vdots&\vdots&\ddots
\end{array}
\right],\ee
(where $\bsym\phi=(\phi_n)_{n\geq 0}$ is a sequence of indeterminates, and for ease of reading we set $\phi_0=1$). It is easy to see that~\eqref{eq.q1Tarray} is the augmented production matrix obtained by setting $\phi_i=y^i$ in $\mathsf{P}_{\bsym\phi}$.

The output matrix of $\mathsf{P}_{\bsym\phi}$ is the \emph{generic Lah triangle}:
\be\label{eq.phioutput} \mathsf{L}\ceq(L_{n,k}(\bsym\phi))_{n,k\geq0}\;=\;
\left[
\begin{array}{ccccccc}
 1&&&&&&\\
 0& 1 &  &  &  &  \\
 0& \phi_1 & 1 &  &  &  \\
 0& \phi_1^2+2\phi_2 & 3\phi_1 & 1 &  &  \\
 0& \phi_1^3+8\phi_1\phi_2+6\phi_3 & 7\phi_1^2+8\phi_2 & 6\phi_1 & 1 & \\
\vdots&\vdots&\vdots&\vdots&\vdots&\ddots
\end{array}
\right],\ee
the row-generating polynomials of which are
\be L_n(\bsym\phi,y)\ceq\sum_{k=0}^nL_{n,k}(\bsym\phi)y^k.\ee
The entries of $\mathsf{L}$ are generating polynomials for \emph{unordered forests of increasing ordered trees} on the vertex set $[n]$ having $k$ components, in which each vertex with $i$ children is assigned a weight $\phi_i$.\footnote{An \emph{ordered tree} is a rooted tree in which the children of each vertex are linearly ordered. An unordered forest of ordered trees is an unordered collection of ordered trees. An \emph{increasing ordered tree} is an ordered tree in which the vertices carry distinct labels from a linearly ordered set (usually some set of integers) in such a way that the label of each child is greater than the label of its parent; otherwise put, the labels increase along every path downwards from the root. An \emph{unordered forest of increasing ordered trees} is an unordered forest of ordered trees with the same type of labeling.}
Setting $\bsym\phi=(1)_{n\geq 0}$ in~\eqref{eq.phioutput} we obtain the \emph{reverse Bessel triangle}:

\be \mm(1,1,1,1)\;=\;\left[
\begin{array}{ccccccc}
 1 &  &  &  &  &  &\\
 0 & 1 &  &  &  &  &\\
 0 & 1 & 1 &  &  &  &\\
 0 & 3 & 3 & 1 &  &  &\\
 0 & 15 & 15 & 6 & 1 &  &\\
 0 & 105 & 105 & 45 & 10 & 1 &\\
 \vdots&\vdots&\vdots&\vdots&\vdots&\vdots&\ddots
\end{array}
\right]\ee
(see~\cite[A001497]{OEIS}). The row-generating polynomials of $\mm(1,1,1,1)$ are the \emph{reverse Bessel polynomials} (see~\cite{Grosswald78})
\be\theta_n(x)\;\ceq\;\sum_{k=0}^n\frac{(n+k)!}{(n-k)!k!}\left(\frac{x^{n-k}}{2^k}\right),\ee
which is an orthogonal sequence of polynomials related to the \emph{Bessel polynomials} (see~\cite{Krall49,Grosswald78}):
\be y_{n}(x)\ceq\sum_{k=0}^n\frac{(n+k)!}{(n-k)!k!}\left(\frac{x}{2}\right)^k\ee
via the identity
\be \theta_n(x)\;=\;x^ny_n(1/x).\ee
It follows from~\cite{Petreolle_2021} that the sequence of polynomials $\theta_n(x)$ is coefficientwise Hankel-totally positive. 

Consider now the row-generating polynomials of $\mm(r,s,\mu,y)$, defined to be
\be \mm_0(x;r,s,\mu,y)\;\ceq\;1\ee
and
\be \mm_n(x;r,s,\mu,y)\ceq\sum_{k=0}^nm_{n,k}(r,s,\mu,y)x^k\ee
for $n>0$, and note that the polynomials $\mm_n(x;r,s,\mu,y)$ are generalisations of the reverse Bessel polynomials in the sense that $\mm_n(x;1,1,1,1)=\theta_n(x)$. 

It follows from~\cite{Petreolle_2021} that the polynomial sequence $(\mm_n(x;1,1,1,y))_{n\geq 0}$ is coefficientwise totally positive jointly in $x$ and $y$. Again, it is natural to ask whether the generalisation where $r,s,\mu,y$ are left as indeterminates preserves coefficientwise Hankel-totally positivity. The answer is no, since the Hankel matrix \be H_{\infty}((\mm_n(x;r,s,\mu,y))_{n\geq 0})\;=\;(\pi_{n,k})_{n,k\geq0},\ee
where
\be \pi_{n,k}\;=\;\mm_{n+k}(x;r,s,\mu,y)\ee
 contains the $2\times 2$ minor
\be \pi_{1,0}\pi_{2,1}-\pi_{2,0}\pi_{1,1}\;=\;\mu x^2 y^2+\mu x^3 y+r x^2 y^2+r x^3 y+s x^2 y^2+s x^3 y-x^2 y^2-2 x^3 y\not\myge0.\ee
Once more, however, it appears that Hankel-total positivity could well be restored by right-multiplying $\mm(r,s,\mu,y)$ by a simple diagonal matrix and shifting $r,s,\mu$ by $1$. 

Consider the matrix
\be \mm'(r,s,\mu,y)\ceq (m'_{n,k}(r,s,\mu,y))_{n,k\geq 0}\;=\;\mm(r,s,\mu,y)\diag(((rs\mu)^{n(n-1)/2})_{n\geq 0})\ee
and let 
\be  \mm'_0(x;r,s,\mu,y)\ceq 1\ee
and
\be  \mm'_n(x;r,s,\mu,y)\ceq \sum_{k=0}^nm'_{n,k}(r,s,\mu,y)x^k\ee
for $n>0$ be the corresponding row-generating polynomials of $\mm'(r,s,\mu,y)$. I conjecture:
\begin{conjecture}
The sequence of polynomials $(\mm'_n(x;r+1,s+1,\mu+1,y))_{n\geq 0}$ is coefficientwise Hankel-totally positive jointly in $r,s,\mu,x,y$.
\end{conjecture}
I have verified this up to $5\times 5$, and am yet to discover a suitable combinatorial interpretation of the entries of $\mm(r,s,\mu,y)$ or indeed $\mm'(r,s,\mu,y)$, and therefore pose the following problem:
\begin{problem}
Find combinatorial interpretations of $\mm'(r,s,\mu,y)$ and $\mm(r,s,\mu,y)$.
\end{problem}

There is one final observation to make concerning perfect matchings that may help in tackling the foregoing problem.
A \emph{perfect matching} $M$ on $[2n]$ is a partition of $1,2,\ldots,2n$ into $n$ pairs $(i,j)$ where $i<j$. We denote the set of all such perfect matchings $\scrm_{2n}$, and it is well-known that
\be |\scrm_{2n}|\;=\;(2n-1)!!.\ee
Two pairs $(i_r,j_r)$ and $(i_s,j_s)$ belonging to $M\in\scrm_{2n}$ form a \emph{crossing} if $i_r<i_s<j_r<j_s$, and form a \emph{nesting} if $i_r<i_s<j_s<j_r$. The number of crossings in $M$ is denoted $\cross(M)$ and the number of nestings is denoted $\nest(M)$.

If $\mu$ and $y$ are specialised to $\mu=y=1$ then the first few rows of $\mm(r,s,1,1)$ are:
 \be  {\small\left[
\begin{array}{cccccc}
 1 &  &  &  &  &\\
 0 & 1 &  &  &  &\\
 0 & 1 & 1 &  &  &\\
 0 & r+s+1 & r+s+1 & 1 &&  \\
 0& 2r^2 s+r^3+r^2+2 r s^2+2 r s\hfill&2 r^2 s+r^3+r^2+2 r s^2+2 r s\hfill&r^2+r s+r\hfill&1&\\
  & \hfill +2 r+ s^3+s^2+2 s+1& \hfill+2 r+s^3+ s^2+2 s+1& \hfill +s^2+s+1 &  &\\
 \vdots&\vdots &\vdots&\vdots&\vdots&\ddots\\
\end{array}
\right]}.\ee
For $s$ specialised to $s=1$ and $0<n\leq 20$ I have verified empirically that
\be m_{n,1}(r,1,1,1)\;=\;(1-r)^{1-n}\sum_{k=0}^{n-1}(-1)^kr^{k(k+1)/2}\frac{2k+1}{2n-1}\binom{2n-1}{n+k}.\ee
The formula above (where $n$ is replaced with $n+1$) was found implicitly by Touchard~\cite{Touchard52}, and explicitly by Riordan~\cite{Riordan75} (see also~\cite{Blitvic12,Kim13,Penaud95,Prodinger12,Read79}), and counts perfect matchings with respect to crossings (or nestings):
\be\sum_{M\in\scrm_{2n}}r^{\cross(\scrm)}\;=\;\sum_{M\in\scrm_{2n}}r^{\nest(\scrm)}\;=\;(1-r)^{-n}\sum_{k=0}^{n}(-1)^kr^{k(k+1)/2}\frac{2k+1}{2n+1}\binom{2n+1}{n+k+1}.\ee

However, more may indeed be true, since it looks like the $k=1$ column of $\mm(r,s,1,1)$ counts perfect matchings with respect to crossings \emph{and} nestings. For $0<n\leq 20$ my computations have confirmed that
\be m_{n,1}(r,s,1,1)\;=\;[t^n]S(t;r,s)\ee
where $S(t;r,s)$ is the $S$-type continued fraction
\be S(t;r,s)\ceq\cfrac{1}{1-\cfrac{[1]_{r,s}t}{1-\cfrac{[2]_{r,s}t}{\cfrac{\cdots}{1-\cfrac{[n]_{r,s}t}{\cdots}}}}}.
\ee
Since $S(t;r,s)$ is the master $S$-fraction for perfect matchings first given by Kasraoui and Zeng in~\cite[Proposition~4.4, with $p,q$ replaced by $r,s$ respectively]{Kasraoui06} (see also~\cite[Theorem~4.4]{sokal2020multivariate} specialised to $x=y=u=v=1$, $p_+=p_-=r$, and $q_+=q_-=s$), I make the following conjecture:
\begin{conjecture}
For $n>1$
\be m_{n,1}(r,s,1,1)\;=\;\sum_{M\in\scrm_{2(n-1)}}r^{\cross(M)}s^{\nest(M)}.\ee
\end{conjecture}
As described above, there exists an explicit formula for counting perfect matchings with respect to crossings \emph{or} nestings given by Touchard and Riordan, however, an explicit closed formula that counts perfect matchings with respect to crossings \emph{and} nestings is yet to be found. I therefore pose the following problem:
\begin{problem}
Find a closed form formula for $m_{n,1}(x;r,s,1,1)$.
\end{problem}

\section*{Acknowledgements}
I would like to thank Shaoshi Chen, Xi Chen, Bishal Deb, Alan Sokal, and Alexander Dyachenko for many helpful and illuminating conversations.

I would also like to thank the many volunteers who maintain and expand the Online Encyclopedia of Integer Sequences~\cite{OEIS}, and Neil Sloane for founding it in the first place. It is a fabulous resource that I make use of almost daily.

This research was supported by Engineering and Physical Sciences Research Council grant EP/N025636/1.


\begin{thebibliography}{10}

\bibitem{Aigner18}
M.~Aigner and G.~M. Ziegler, \emph{Proofs from {T}he {B}ook}, 6th ed.,
  Springer-Verlag, Berlin, 2018.

\bibitem{andrews84}
G.~E. Andrews, \emph{The {T}heory of {P}artitions}, Encyclopedia of
  {M}athematics and its {A}pplications, Cambridge University Press, 1984.

\bibitem{Avron2016}
A.~Avro and N.~Dershowitz, \emph{{{C}}ayley's formula: A page from {T}he
  {B}ook}, The American Mathematical Monthly \textbf{123} (2016), no.~7,
  pp.~699 -- 700.

\bibitem{Blitvic12}
N.~Blitvi\'c, \emph{The $(q,t)$-{{G}}aussian process}, J.~Funct.~Anal.
  \textbf{263} (2012), pp.~3270 -- 3305.

\bibitem{Brenti_95}
F.~Brenti, \emph{Combinatorics and total positivity}, J.~Combin.~Theory A
  \textbf{71} (1995), pp.~175 -- 218.

\bibitem{Brenti_96}
\bysame, \emph{The applications of total positivity to combinatorics, and
  conversely.}, Total Positivity and its Applications (M.~Gasca and C.~A.
  Micchelli, eds.), Kluwer, Dordrecht, 1996, pp.~451 -- 473.

\bibitem{Chauve99}
C.~Chauve, S.~Dulucq, and O.~Guibert, \emph{Enumeration of some labelled
  trees}, Research Report RR-1226-99 (1999), LaBRI, Universite Bordeaux I.

\bibitem{Chauve00}
\bysame, \emph{Enumeration of some labelled trees}, Formal Power Series and
  Algebraic Combinatorics (FPSAC 2000), Springer-Verlag, Berlin, 2000, pp.~146
  -- 157.

\bibitem{ChenYang21}
W.~Y.~C. Chen and H.~R.~L. Yang, \emph{A context-free grammar for the
  {{R}}amanujan–-{{S}}hor polynomials}, Adv.~Appl.~Math. \textbf{126} (2021),
  no.~101908, p.~24.

\bibitem{Gilmore22}
Xi. Chen, B.~Deb, A.~Dyachenko, T.~Gilmore, and A.~D. Sokal,
  \emph{Coefficientwise total positivity of some matrices defined by linear
  recurrences}, in preparation.

\bibitem{Gilmore20}
\bysame, \emph{Coefficientwise total positivity of some matrices defined by
  linear recurrences},  (2020), to appear in the proceedings of FPSAC 2021.

\bibitem{Cheon2013}
G.-S. Cheon, J.-H. Jung, and Y.~Lim, \emph{A $q$-analogue of the {{R}}iordan
  group}, Linear Algebra and its Applications \textbf{439} (2013), no.~12,
  pp.~4119 -- 4129.

\bibitem{Cigler_08}
J.~Cigler, \emph{$q$-{{A}}bel polynomials},  (2008), arXiv:0802.2886.

\bibitem{Clarke58}
L.~E. Clarke, \emph{On {{C}}ayley's {{F}}ormula for {{C}}ounting {{T}}rees},
  Journal of the London Mathematical Society \textbf{1-33} (1958), no.~4,
  pp.~471--474.

\bibitem{Comtet_74}
L.~Comtet, \emph{Advanced {{C}}ombinatorics: The {{A}}rt of {{F}}inite and
  {{I}}nfinite expansions}, Reidel, Dordrecht--Boston, 1974, translated from
  the French by J.W.~Nienhuys. [French original: {\em Analyse Combinatoire}\/,
  tomes~I et II, Presses Universitaires de France, Paris, 1970.].

\bibitem{Corless96}
R.~M. Corless, G.~H. Gonnet, D.~E.~G. Hare, D.~J. Jeffrey, and D.~E. Knuth,
  \emph{On the {{L}}ambert {{W}} function}, Advances in Computational
  Mathematics \textbf{5} (1996), no.~1, pp.~329 -- 359.

\bibitem{DeMedicis_93}
A.~de~M{\'e}dicis and P.~Leroux, \emph{A unified combinatorial approach for
  $q$- (and $p,q$-) {{S}}tirling numbers}, Journal of Statistical Planning and
  Inference \textbf{34} (1993), no.~1, pp.~89 -- 105.

\bibitem{Deutsch09}
E.~Deutsch, L.~Ferrari, and S.~Rinaldi, \emph{Production {{M}}atrices and
  {{R}}iordan {{A}}rrays}, Annals of Combinatorics \textbf{13} (2009), no.~1,
  pp.~65 -- 85.

\bibitem{Deutsch04}
E.~Deutsch and L.~Shapiro, \emph{Exponential {{R}}iordan arrays}, handwritten
  lecture notes, Nankai University, available online at
  \url{http://www.combinatorics.net/ppt2004/Louis%20W.%20Shapiro/shapiro.pdf},
  26 February 2004.

\bibitem{Dumont96}
D.~Dumont and A.~Ramamonjiosa, \emph{Grammaire de {{R}}amanujan et arbres de
  {{C}}ayley}, Electron.~J.~Combin. \textbf{3} (1996), no.~R17.

\bibitem{Egecioglu86}
{\"O}.~E{\v{g}}ecio{\v{g}}lu and J.~B. Remmel, \emph{Bijections for {{C}}ayley
  trees, spanning trees, and their $q$-analogues}, J.~Combin.~Theory A
  \textbf{42} (1986), no.~1, pp.~15 -- 30.

\bibitem{Fallat_11}
S.~M. Fallat and C.~R. Johnson, \emph{Totally {{N}}onnegative {{M}}atrices},
  Princeton University Press, Princeton NJ, 2011.

\bibitem{Flajolet80}
F.~Flajolet, \emph{Combinatorial aspects of continued fractions}, Discrete
  Math. (1980), pp.~125 -- 161.

\bibitem{Fomin_00}
S.~Fomin and A.~Zelevinsky, \emph{Total positivity: tests and
  parametrizations}, Math.~Intelligencer \textbf{22} (2000), pp.~23 -- 33.

\bibitem{Francon_74}
Jean Fran\c{c}on, \emph{Preuves combinatoires des identites d'{{A}}bel},
  Discrete Mathematics \textbf{8} (1974), no.~4, pp.~331 -- 343.

\bibitem{Gantmacher_02}
F.~R. Gantmacher and M.~G. Krein, \emph{Oscillation {{M}}atrices and
  {{K}}ernels and {{S}}mall {{V}}ibrations of {{M}}echanical {{S}}ystems}, AMS
  Chelsea Publishing, Providence RI, 2002, based on the second Russian edition,
  1950.

\bibitem{Gasca87b}
G.~Gasca and G.~M{\"u}hlbach, \emph{A test for strict total positivity via
  {{N}}eville elimination}, Current {{T}}rends in {{M}}atrix {{T}}heory
  (F.~Uhlig and R.~Groue, eds.), Elsevier Science, 1987, pp.~225 -- 232.

\bibitem{Gasca87}
M.~Gasca and E.~Lebr{\`o}n, \emph{Elimination techniques and interpolation},
  J.~Comput.~Appl.~Math. \textbf{19} (1987), pp.~125 -- 132.

\bibitem{Gasca87a}
M.~Gasca and G.~M{\"u}hlbach, \emph{Generalized {{S}}chur complements and a
  test for total positivity}, Appl.~Numer.~Math. \textbf{3} (1987), pp.~215 --
  232.

\bibitem{Gasca92}
M.~Gasca and J.~M. Pe{\~n}a, \emph{Total positivity and {{N}}eville
  elimination}, Linear Algebra and its Applications \textbf{165} (1992), pp.~25
  -- 44.

\bibitem{Gessel_06}
I.~M Gessel and S.~Seo, \emph{A refinement of {{C}}ayley's formula for trees},
  Elec.~J.~Combin. \textbf{11} (2006), no.~R.27.

\bibitem{Gessel062}
Ira~M. Gessel, \emph{A q-analog of the exponential formula}, Discrete
  Mathematics \textbf{306} (2006), no.~10, pp.~1022 -- 1031, 35th Special
  Anniversary Issue.

\bibitem{Grosswald78}
E.~Grosswald, \emph{Bessel {{P}}olynomials}, Lecture Notes Math. \textbf{698}
  (1978).

\bibitem{Guo17}
S.~Guo and J.~W. Guo, \emph{A recursive algorithm for trees and forests},
  Discrete Mathematics \textbf{340} (2017), no.~4, pp.~695 -- 703.

\bibitem{Guo18}
V.~J.~W. Guo, \emph{A bijective proof of the {{S}}hor recurrence}, European J.
  Combin. \textbf{70} (2018), pp.~92 -- 98.

\bibitem{Guo07}
V.~J.~W. Guo and J.~Zeng, \emph{A generalization of the {{R}}amanujan
  polynomials and plane trees}, Avd.~Appl.~Math. \textbf{39} (2007), pp.~96 --
  115.

\bibitem{Jackson_1910}
F.~H. Jackson, \emph{A $q$-generalization of {{A}}bel's series}, Rendiconti
  Palermo \textbf{29} (1910), pp.~340 -- 346.

\bibitem{Johnson_96}
W.~P. Johnson, \emph{$q$-extensions of identities of {{A}}bel-{{R}}othe type},
  Discrete Mathematics \textbf{159} (1996), no.~1, pp.~161 -- 177.

\bibitem{Johnson96}
Warren~P. Johnson, \emph{Some applications of the q-exponential formula},
  Discrete Mathematics \textbf{157} (1996), no.~1, pp.~207 -- 225.

\bibitem{Jousat15}
M.~Josuat-Verg{\`e}s, \emph{Derivatives of the tree function}, Ramanujan J.
  \textbf{38} (2015), pp.~1 -- 15.

\bibitem{Kim13}
M.~Josuat-Verg\`es and J.~S. Kim, \emph{{{T}}ouchard–-{{R}}iordan formulas,
  {{T}}--fractions, and {{J}}acobi's triple product identity}, Ramanujan J.
  \textbf{30} (2013), pp.~341 -- 378.

\bibitem{Karlin_68}
S.~Karlin, \emph{Total {{P}}ositivity}, Stanford University Press, Stanford CA,
  1968.

\bibitem{Kasraoui06}
A.~Kasraoui and J.~Zeng, \emph{Distribution of crossings, nestings and
  alignments of two edges in matchings and partitions}, Elec.~J.~Comb.
  \textbf{13} (2006), no.~1.

\bibitem{Krall49}
H.~L. Krall and O.~Frink, \emph{A new class of orthogonal polynomials: The
  {{B}}essel polynomials}, Trans.~Amer.~Math.~Soc. \textbf{63} (1949), pp.~100
  -- 115.

\bibitem{Lin14}
Z.~Lin and J.~Zeng, \emph{Positivity properties of {{J}}acobi–{{S}}tirling
  numbers and generalized {{R}}amanujan polynomials}, Adv.~Appl.~Math
  \textbf{53} (2014), pp.~12 -- 27.

\bibitem{Moon70}
J.~W. Moon, \emph{Counting labelled trees}, Canadian Mathematical Congress,
  Montreal, 1970.

\bibitem{Penaud95}
J.-G. Penaud, \emph{Une preuve bijective d'une formule de
  {{T}}ouchard–-{{R}}iordan}, Discrete Math. \textbf{139} (1995), pp.~347 --
  360.

\bibitem{Petreolle_2021}
M.~P{\'e}tr{\'e}olle and A.~D. Sokal, \emph{Lattice paths and branched
  continued fractions ii. {{M}}ultivariate {{L}}ah polynomials and {{L}}ah
  symmetric functions}, European Journal of Combinatorics \textbf{92} (2021),
  pp.~103 -- 235.

\bibitem{petreolle2020}
M.~P{\'e}tr{\'e}olle, A.~D. Sokal, and B.~Zhu, \emph{Lattice paths and branched
  continued fractions: An infinite sequence of generalizations of the
  {{S}}tieltjes--{{R}}ogers and {{T}}hron--{{R}}ogers polynomials, with
  coefficientwise {{H}}ankel-total positivity},  (2020), accepted for
  publications in Memoirs of the AMS.

\bibitem{Pfaff95}
J.~F. Pfaff, \emph{Allgemeine {{S}}ummation einer {{R}}eihe, worinn h{\"o}here
  differenziale vorkommen}, Archiv der reinen und angewandten Mathematik
  \textbf{1} (1795), pp.~337 -- 347.

\bibitem{Pinkus_10}
A.~Pinkus, \emph{Totally {{P}}ositive {{M}}atrices}, Cambridge University
  Press, Cambridge, 2010.

\bibitem{Pitman_02}
J.~Pitman, \emph{Forest {{V}}olume {{D}}ecompositions and
  {{A}}bel–-{{C}}ayley-–{{H}}urwitz {{M}}ultinomial {{E}}xpansions},
  J.~Combin.~Theory, Series A \textbf{98} (2002), no.~1, pp.~175 -- 191.

\bibitem{Prodinger12}
H.~Prodinger, \emph{On {{T}}ouchard's continued fraction and extensions:
  combinatorics-free, self-contained proofs}, Quaestiones Math. \textbf{35}
  (2012), pp.~431 -- 445.

\bibitem{Prufer18}
H.~Pr{\"u}fer, \emph{Neuer {{B}}eweis eines {{S}}atzes \"{u}ber
  {{P}}ermutationen}, Arch.~Math.~Phys. \textbf{27} (1918), pp.~742 -- 744.

\bibitem{Randazzo21}
L.~Randazzo, \emph{Arboretum for a generalisation of {{R}}amanujan
  polynomials}, Ramanujan J. \textbf{54} (2021), pp.~591 -- 604.

\bibitem{Read79}
R.~C. Read, \emph{The chord intersection problem}, Second International
  conference on Combinatorial Mathematics: Annals of the New York Academy of
  Sciences No.~319 (A.~Gerwitz and L.~V Quintas, eds.), New York Academy of
  Sciences, New York, 1979, pp.~pp.~444 -- 454.

\bibitem{Riordan68}
J.~Riordan, \emph{Forests of labeled trees}, J.~Combin.~Theory \textbf{5}
  (1968), pp.~90 -- 103.

\bibitem{Riordan75}
\bysame, \emph{The distribution of crossings of chords joining pairs of $2n$
  points on a circle}, Math.~Comp. \textbf{29} (1975), pp.~215 -- 222.

\bibitem{Rothe93}
H.~A. Rothe, \emph{Formulae de serierum reversione demonstratio universalis
  signislocalibus combinatorioanalyticorum vicariis exhibita}, Litteris
  Sommeriis, Leipzig, 1793.

\bibitem{Sagan_83}
Bruce~E. Sagan, \emph{A note on {{A}}bel polynomials and rooted labeled
  forests}, Discrete Mathematics \textbf{44} (1983), no.~3, pp.~293 -- 298.

\bibitem{Schlafli47}
L.~Schl{\"a}fli, \emph{Bemerkungen {\"u}ber die {{L}}ambertische {{R}}eihe},
  Archiv der Mathematikund Physik \textbf{10} (1847), pp.~332 –-- 340,
  Reprinted in L. Schl{\"a}fli, \emph{Gesammelte Mathematische Abhandlungen},
  Band I (Springer, Basel, 1950), pp. 38–45.

\bibitem{Shor95}
P.~W. Shor, \emph{A new proof of {{C}}ayley's formula for counting labeled
  trees}, J.~Combin.~Theory A \textbf{71} (1995), pp.~154 -- 158.

\bibitem{OEIS}
N.~Sloane, \emph{The {{O}}nline {{E}}ncyclopedia of {{I}}nteger {{S}}equences},
  available at \url{http://oeis.org}.

\bibitem{Sokal22}
A.~D. Sokal, \emph{Coefficientwise total positivity (via continued fractions)
  for some hankel matrices of combinatorial polynomials}, in preparation.

\bibitem{Sokal_20}
\bysame, \emph{A remark on the enumeration of rooted labeled trees}, Discrete
  Mathematics \textbf{343} (2020), no.~7.

\bibitem{Sokal_21f}
A.~D. Sokal, \emph{Total {{P}}ositivity of some polynomial matrices that
  enumerate labelled trees and forests: {{I}}. {{F}}orests of rooted labelled
  trees},  (2021), arXiv:2105.05583.

\bibitem{Sokal_21t}
A.~D. Sokal and X.~Chen, \emph{Total positivity of some polynomial matrices
  that enumerate labelled trees and forests: {{II}}. {{R}}ooted labelled
  trees},  (2021), in preparation.

\bibitem{sokal2020multivariate}
A.~D. Sokal and J.~Zeng, \emph{Some multivariate master polynomials for
  permutations, set partitions, and perfect matchings, and their continued
  fractions}, 2020, arXiv:2003.08192.

\bibitem{Stanley_86}
R.~P. Stanley, \emph{Enumerative {{C}}ombinatorics}, vol.~2, Wadsworth \&
  Brooks/Cole, Monterey CA, 1986, reprinted by Cambridge University Press,
  1999.

\bibitem{Stieltjes94}
T.~J. Stieltjes, \emph{Recherches sur les fractions continues},
  Ann.~Fac.~Sci.~Toulouse \textbf{8} (1894), no.~J1 -- J122.

\bibitem{Stieltjes95}
\bysame, \emph{Recherches sur les fractions continues}, Ann.~Fac.~Sci.~Toulouse
  \textbf{9} (1895), no.~A1 -- A47.

\bibitem{Stieltjes93}
\bysame, \emph{Euvres {{C}}ompl{\`e}te/{{C}}ollected {{P}}apers}, vol.~II,
  Springer-Verlag, Berlin, 1993.

\bibitem{Touchard52}
J.~Touchard, \emph{Sur un probl{\`e}me de configurations et sur les fractions
  continues}, Canad.~J.~Math \textbf{4} (1952), pp.~2 -- 25.

\bibitem{Zeng95}
J.~Zeng, \emph{The $q$-{{S}}tirling numbers, continued fractions and the
  $q$-{{C}}harlier and $q$-{{L}}aguerre polynomials}, Journal of Computational
  and Applied Mathematics \textbf{57} (1995), no.~3, pp.~413 -- 424.

\bibitem{Zeng99}
\bysame, \emph{A {{R}}amanujan sequence that refines the {{C}}ayley formula for
  trees}, Ramanujan J. \textbf{3} (1999), pp.~45 -- 54.

\end{thebibliography}
\end{document}